\documentclass[12pt]{article}
\usepackage{setspace}
\usepackage[utf8]{inputenc}\usepackage{authblk}
\usepackage{graphicx}
\usepackage{amsmath}
\usepackage{amssymb,amsthm}
\usepackage{geometry}
\usepackage{cases}
\usepackage[hidelinks]{hyperref}
\usepackage{cleveref}

\graphicspath{ {images/} }
\newtheorem{theorem}{Theorem}
\newtheorem{lemma}{Lemma}

\newtheorem{prop}{Proposition}

\theoremstyle{definition}
\newtheorem{rk}{Remark}

\DeclareMathOperator*{\supess}{ess\,sup}
\DeclareMathOperator*{\supp}{supp}

\begin{document}

%\tableofcontents
%\printbibliography
\def\R{\mathbb{R}}
\def\H{\mathbb{H}}
\def\C{\mathbb{C}}
\def\g{\rm{\bf g}}
\def\a{\alpha}
\def\e{\epsilon}\def\s{\sigma} 
\def\si{\sigma}
\def\ga{\gamma}
\def\la{\lambda}
\def\k{\kappa}
\def\ds{\displaystyle}
\def\ni{\noindent}
\def\avint{\mathop{\,\rlap{\bf--}\!\!\int}\nolimits}
\def\ba{\begin{aligned}}
\def\ea{\end{aligned}}
\def\ol{\overline}
\def\wt{\widetilde}
\def\b{\beta}\def\ab{{\overline\a}}
\def\avgdj{\avint_{(D_{j+1}\setminus D_{j})\setminus E_\tau}}
 \def\avgdJ{\avint_{(D_{J+1}\setminus D_{J})\setminus E_\tau}}
 \def\avgdJo{\avint_{(D_{J_1+1}\setminus D_{J_1})\setminus E_\tau}}
 \def\dJ{\int_{(D_{J+1}\setminus D_{J})\setminus E_\tau}}
 \def\dj{\int_{(D_{j+1}\setminus D_{j})\setminus E_\tau}}
 \def\wphi{\wt{\phi}_{\e,r}}
 \def\mc{\mathcal}\def\inj{{\text {inj}}}
 \def\wh{\widehat}
 \def\AA{\max\{A_0,A_\infty\}}
\def\be{\begin{equation}}\def\ee{\end{equation}}\def\nt{\notag}\def\bc{\begin{cases}}\def\ec{\end{cases}}\def\ba{\begin{aligned}}
\def\ea{\end{aligned}}
\def\qed{{\setlength{\fboxsep}{0pt}\setlength{\fboxrule}{0.2pt}\fbox{\rule[0pt]{0pt}{1.3ex}\rule[0pt]{1.3ex}{0pt}}}}
\def\QEDopen{{\setlength{\fboxsep}{0pt}\setlength{\fboxrule}{0.2pt}\fbox{\rule[0pt]{0pt}{1.3ex}\rule[0pt]{1.3ex}{0pt}}}}
\def\na{\sum_{k=0}^{2N-1}}
\def\ia{\sum_{k=2N}^{\infty}}
\def\aa{\sum_{k=0}^{\infty}}\def\Tfc{(Tf)^\circ}
\def\udl{\underline}\def\llceil{\left\lceil}\def\rrceil{\right\rceil}
\def\bfe{{\bf e_1}}\newcommand\fH[1]{\sbox0{#1}\dimen0=\ht0 \advance\dimen0 -1ex
  \sbox2{\'{}}\sbox2{\raise\dimen0\box2}%
  {\ooalign{\hidewidth\kern.1em\copy2\kern-.5\wd2\box2\hidewidth\cr\box0\crcr}}}
\def\S{\mathbb{S}}\def\SS{\mathcal{S}}
\centerline{\bf \Large{Sharp Adams inequalities with exact growth conditions }}
\centerline{\bf\Large{on metric measure spaces and applications}}
\font\twelvemi=cmmi12 at 13pt\font\elevenmi=cmmi11 at 9 pt
\renewcommand{\chi}{\raisebox{.13\baselineskip}{\hbox{\twelvemi\char31}}}
\newcommand{\gam}{\gamma}
\newcommand{\sgam}{\gamma}
\def\lan{\big\langle}\def\ran{\big\rangle}
\vskip 0.1in
\centerline{\bf{Carlo Morpurgo,\ Liuyu Qin\footnote{The second author was supported by the National Natural Science Foundation of China (12201197) \;\;\;\; 

\smallskip MSC Class: 46E36; 26D10}}}

%
%
%
%%
%% \chapter*{Acknowledgements}
\begin{abstract} 
Adams  inequalities with exact growth conditions are derived for Riesz-like potentials on metric measure spaces. The results extend and improve                                                                                                                                                                                                                                                                                                                                                                              those obtained recently on $\R^n$ by the second author, for Riesz-like convolution operators. As a consequence, we will obtain new sharp Moser-Trudinger inequalities with exact growth conditions on~$\R^n$, the Heisenberg group, and  Hadamard  manifolds. On $\R^n$ such inequalities will be used to prove the existence of radial ground states solutions for a class of quasilinear elliptic equations.
\end{abstract}
\numberwithin{equation}{section}

\

 \section{Introduction}

\def\lc{\left\lceil}\def\rc{\right\rceil}\def\na{\frac n{\a}}\def\nna{\frac{n}{n-\a}}

There are several extensions of the classical Moser-Trudinger inequality (MTI) for functions in $W_0^{\a,n/\a}(\Omega)$, with $\Omega\subset \R^n$ and $\a$ an integer with $0<\a<n$. The most important types can be  categorized using  the single inequality 
 \be \int_{\{|u|\ge 1\}} \frac{\exp\left[\gamma\,|u|^{\nna}\right]}{|u|^\lambda}\,dx\le C\big(1-\e+\e\|u\|_{n/\a}^{n/\a}\big),\qquad \kappa\|u\|_{n/\a}^{n/\a}+\|\nabla^\a u\|_{n/\a}^{n/\a}\le1\label{MS1}\ee
 where  $\nabla^\a=(-\Delta)^{\frac{\a}2}$ if  $\a$ is even, and $\nabla^\a=\nabla(-\Delta)^{\frac{\a-1}2}$ if  $\a$ is odd.
 
 \def\la{\lambda}
 The {\it classical MTI} corresponds  to the case $\la=\e=\kappa=0$, and it holds for $\gamma$ up to and including an explicit constant $\gamma_{n,\a}$ (see \eqref{mf4}) and not for a larger $\gamma$, provided that $\Omega$ is an open set with finite measure. This is the classical result by Adams [A1]. More generally, in [FM3] it was shown that the classical MTI  holds if $\Omega$ is Riesz-subcritical, i.e., roughly speaking, if it misses enough dimensions at infinity.
 
 \smallskip
 
 The {\it Ruf MTI} corresponds to $\la=\e=0$ and $\kappa=1$, and it  holds on $W_0^{\a,n/\a}(\R^n)$ for $\gamma\le \gamma_{n,\a}$ and it fails for $\gamma>\gamma_{n,\a}$.  The first such inequality was derived by Ruf in [Ruf] for the gradient in $\R^2$, later extended to all dimensions in [LR], the Laplacian in [LL2], and to arbitrary $\alpha$ in [FM2]. The inequality fails at $\gamma=\gamma_{n,\a}$ (and it's true for $\gamma<\gamma_{n,\a}$) if the Ruf norm condition is replaced by the weaker condition $\max\big\{\|u\|_{n/\a}^{n/\a},\|\nabla^\a u\|_{n/\a}^{n/\a}\big\}\le1$.

\smallskip
The {\it Adachi-Tanaka MTI} is  the case $\la=\kappa=0$ and $\e=1$, and it is true only for $\gamma<\gamma_{n,\a}$ and not true for $\gamma=\gamma_{n,\a}$ regardless of whether or not $\Omega$ has infinite measure [AT], [FM2]. In fact, such an inequality with $\gamma=\gamma_{n,\a}$ fails even if $0<\la<\nna$. The Adachi-Tanaka inequality is implied by  the Ruf MTI.

\eject

\smallskip
The {\it MTI  with exact growth condition}, which we will also call {\it Masmoudi-Sani inequality (MSI)}, corresponds to $\la=\nna,\; \e=1,\;\kappa=0$, and it holds on any $\Omega\subseteq \R^n$  for $\gamma\le \gamma_{n,\a}$ while it fails for $\gamma>\gamma_{n,\a}$. This result was derived first in [IMN] for $n=2,\; \a=1$, followed by [MS1] for $n=4\;,\a=2$. For general $n$,   [LTZ] settled the case $\alpha=2$, [MS2] the case $\alpha=1$, and  [MS3] any integer $\alpha<n$. When $\Omega=\R^n$ the inequality yields, in particular,  a uniform bound of the exponential integral with $\gamma=\gamma_{n,\a}$ under the condition $\max\big\{\|u\|_{n/\a}^{n/\a},\|\nabla^\a u\|_{n/\a}^{n/\a}\big\}\le1$ (which is not possible if $\la<\nna$), which is weaker than the Ruf norm condition.  It is in fact possible to derive the Ruf MTI from the MSI, even though both inequalities are in effect giving optimal growth conditions under different Sobolev norms.

\smallskip

When $\Omega$ has infinite measure, or when $\e=1$, all of the above inequalities  can be stated in equivalent forms with the integrals  taken  over the entire $\Omega$, provided the exponential is suitably regularized and the denominator is replaced by $1+|u|^\la$.

\smallskip

As a side note, we mention that  in $\R^2$ the limiting case $\a=2$ in $L^1$ has been studied in [CRT] and [FM4], in the context of the so-called ``reduced Sobolev spaces",  where the Moser type inequalities are often called ``Brezis-Merle inequalities". The case where no boundary condition is assumed, i.e. $W^{\a,n/\a}(\Omega)$, has been settled completely by Cianchi in [Ci1] for $\alpha=1$, whereas when $\alpha=2$ sharp results have been derived  in [FM5] when $\Omega$ is a ball. We also mention  [Ci2], where Cianchi derives  the  several sharp Moser-Trudinger inequalities 
with respect to Frostman (or upper  Ahlfors-regular) measures, on a wide class of bounded domains of $\R^n$, with and without zero boundary conditions.

\bigskip

In this paper we shall call  ``Adams inequality" (AI)  any of the versions of \eqref{MS1} applied to potentials rather than Sobolev functions, namely 

  \be \int_{\{|Tf|\ge 1\}} \frac{\exp\left[\gamma\,|Tf|^{\nna}\right]}{|Tf|^\la}\,dx\le C\big(1-\e+\e\|Tf\|_{n/\a}^{n/\a}\big),\qquad \kappa\|Tf\|_{n/\a}^{n/\a}+\|f\|_{n/\a}^{n/\a}\le1\label{MS2}\ee
 where $f$ is a compactly supported function in $L^{n/\a}(\Omega)$, and $Tf$ is an integral operator with a kernel $k$ which behaves like a Riesz potential, in a suitable sense. The classical AI was first derived by Adams for the Riesz potential in [A1]. In a series of papers ([FM1], [FM2], [FM3]) the first author and Fontana refined the Adams machinery and used it  to derive Adams inequalities in several  settings. In particular, they proved that if $Tf=K*f$, where $K(x)\sim g(x^*)|x|^{\a-n}$ for small $x$, where $x^*=x/|x|$, then the classical AI holds for  $\gamma\le \gamma_g$, where $\gamma_g$ is given explicitly in terms of $g$, provided $\Omega$ has finite measure, or more generally if $K$ satisfies a critical integrability condition at infinity ([FM1], [FM3]). In [FM2] it was shown that for such $K$ the Ruf AI holds for $\gamma\le \gamma_g$ when $\Omega=\R^n$. In all these cases the constant $\gamma_g$ is sharp (largest) provided $K$ is regular enough. These results applied to the Riesz potential imply the corresponding MT inequalities mentioned above.
 
 The Masmoudi-Sani Adams inequality was recently settled by the second author in [Qin2],  for convolution operators $Tf=K*f$ where $K(x)\sim g(x^*)|x|^{\a-n}$ for small $x$ and for large $x$, with the same sharp exponential constant $\gamma_g$. As a consequence Qin derived sharp  MTI with exact growth condition for arbitrary fractional powers of the Laplacian and for general homogeneous  differential operators with constant coefficients. 
 
 \medskip
 \eject
 In this paper we push the techniques used in [Qin2] to general metric measure spaces, and derive sharp and improved Masmoudi-Sani and Ruf's Adams inequalities for Riesz-like potentials in that setting. So far only the classical Adams inequalities are known in the context of measure spaces [FM1], [FM3]; the metric structure was not needed there, but we could not avoid it in this paper.
 
 As a consequence, we will obtain Moser-Trudinger inequalities with exact growth conditions on the Heisenberg group and on Hadamard manifolds. We will also give an application concerning the existence of a ground state solution of a nonlinear PDE in the same spirit as in [MS2], but for a slightly more general equation. Other applications in the context of Riemannian manifolds with nonnegative curvature will appear in a forthcoming paper.
 
\smallskip
The outline of the paper is as follows.

\smallskip
In Section \ref{main} we will first introduce the assumptions we need on our metric measure spaces. We will only require that the measure of balls is a finite continuous function of the radius, and we will not require the doubling property. 

Next, we will introduce the Riesz-like kernels. There are several papers dealing with Riesz kernels on a metric measure space $X$, with metric $d$ and Borel measure $\mu.$ Normally such kernels $k(x,y)$ are  comparable to $\mu\big(B(x,d(x,y)\big)^{-\beta}$ ($\beta>0$),  where $B(x,r)$ is the ball centered at $x$ and with radius $r$, or in other cases they are defined as $d(x,y)^{\a}\mu\big(B(x,d(x,y)\big)^{-1}$. If $\mu$ is Ahlfors-regular of order $d$ (i.e. balls of radius $r$ have measures comparable to $r^d$) these definitions amount to assuming that $k$ is essentially comparable to $d(x,y)^{\a-d}$ (see for example [ARSW], [HK]). The key step in the derivation of sharp MTI using Adams inequalities, is  to represent a Sobolev function in terms of a (pseudo-)differential operator of order $\alpha$, and to  have precise information on the behavior of its fundamental solution (Green's function)  near the diagonal. In general we cannot expect  the Green function behaves like a constant multiple of $d(x,y)^{\a-d}$, or   $\mu\big(B(x,d(x,y)\big)^{-{\frac{d-\a}d}}, $ when $d(x,y)$ is small. On the Heisenberg group $\H^n$, for example, powers higher than 2 of the subLaplacian have a Green function of type $g\big((xy^{-1})^*\big) |xy^{-1}|^{\a-Q}$, where $Q=2n+2$, and $x^*=x/|x|$, for some nonconstant function $g$ on the Heisenberg sphere.

For our purpose we will need sharp asymptotic Riesz-like behavior of the kernel both for small and large distances, which will be captured  by integral conditions on the kernel, rather than pointwise conditions.

The main inequalities of the paper are stated in Theorem 1, the sharpness statements are in Theorem 2. Besides being in the context of measure spaces, our MSI result contains two other  substantial improvements  from the known versions.  The first improvement is the explicit dependence of sharp exponential constant from  either the asymptotic behavior of the kernel along the diagonal, or the one at infinity. 
It has been already observed by Qin [Qin2, Remark 3] that the MSI for the Riesz kernel $|x-y|^{\a-n}$ in $\R^n$ modified to be $2|x-y|^{\a-n}$ for $|x-y|\ge 2$ cannot hold with the same sharp exponential constant for the unmodified  Riesz kernel.

Secondly,  we will allow the power of $|Tf|$ in the denominator to vary, as well as  the norm on the right-hand side, accordingly. As consequence of this result we obtain the following improvement of the  original MSI on $\R^n$:

  \be \int_{\{|u|\ge 1\}} \frac{\exp\left[\gamma_{n,\a}\,|u|^{\nna}\right]}{|u|^{\frac{p\a}{n-\a}}}\,dx\le C\|u\|_p^p,\qquad \|\nabla^\a u\|_{n/\a}^{n/\a}\le1\label{MS3}\ee
  
  \smallskip
   valid for all $p\ge 1$, whereas the  known results in the literature have $p=n/\a.$  The natural space of functions associated to such an inequality is a Sobolev space with mixed norms.

   \medskip
   In Sections \ref{pfthm1},\ref{lemma2},\ref{prop1}  we will prove Theorem 1. The proof follows the same lines as the one in [Qin2], however there are many changes and improvements, due to the different background space, the more general assumptions on the kernels, and the different conclusions of the theorem. 
 The key ideas used in [Qin2] was to use an improved O'Neil's Lemma combined with a suitable potential version of the ``optimal descending growth condition''. Such condition was introduced in [IMN], [MS2] for radial functions in $W^{1,n}(\R^n)$, and states that for $u\in W^{1,n}(\R^n)$,  radially decreasing and such that $u(r)\ge 1$ and $\|(\nabla u)\chi_{\{|x|\ge r\}}\|_n\le 1$ we have
 \be \frac{\exp\left[\gamma_{n,1}\,u(r)^{\frac{n}{n-1}}\right]}{u(r)^{\frac{ n}{n-1}}}\, r^n\le C\|u\chi_{\{|x|\ge r\}} \|_n^n\label{odgc1}\ee
(for simplicity we let $u(x)=u(|x|)$), where $\gamma_{n,1}=n\omega_{n-1}^{\frac1{n-1}}$ is the sharp constant in the MTI for the gradient. This estimate quantifies optimally the growth of radially decreasing Sobolev functions under the sole condition on their gradients, from which the authors in [IMN], [MS2] were able to derive the MSI for the case $\alpha=1$. A version of \eqref{odgc1} for the case $\alpha=2$, using Lorentz norms of the Laplacian, was obtained in [MS1], [MS3], [LTZ].
 
 Following [Qin2] a version of \eqref{odgc1} will be obtained for general potentials, see \eqref{odgc} in Remark \ref{E1}, where the function $u(r)$ above is replaced by 
 \be(Tf)^\circ(\tau)=\supess_{x\in E_\tau}\big|T\big(f\chi_{F_\tau^c\cap B_\tau(x)^c}\big)(x)\big|
 \label{tftau} \ee
where $E_\tau$, $F_\tau$ are level sets of $Tf$ and $f$ resp., with volume $\tau$, and where $B_\tau(x)$ is a ball centered at $x$,  with volume $\tau.$ When $T$ is the Riesz potential in $\R^n$ and $f$  is radially decreasing, then $(T f)^\circ$
can be compared from above or from below by the decreasing rearrangement $(T f)^*$ in some instances, and it actually coincides with it when $Tf=\nabla|y|^{2-n}* f$, if $|f|$ is radially decreasing (see Appendix, Prop. \ref{circstar}).

The potential version of \eqref{odgc1} is stated in Proposition \ref{le102}, in a slightly improved form, together with a version under the Ruf condition. Both of these results are key tools from which the main inequalities will be obtained.
  The proof of Proposition \ref{le102} is the most challenging one, and will be effected through an improved version of a  discretization procedure introduced in [Qin2], which was  itself adapted from the original papers by Masmoudi et al.
  
  \medskip
  In Section \ref{sharpness} we will prove Theorem 2, which gives sufficient conditions for the exponential constants in the inequalities of Theorem 1 to be sharp.
   
\eject
 In Section \ref{pdeapp} we will apply inequality \eqref{MS3} in the case $\a=1$ to show the existence of a radial ground state solution of the nonlinear equation 
 
 \be -\Delta_n u+V_0|u|^{p-2}u=f(u)\label {qlequation}\ee
 where $\Delta_n$ is the $n-$Laplacian,  $V_0$ is a positive constant, $p>1$,  $\|u\|_p<\infty$, $\|\nabla u\|_n<\infty$, and where $f$ satisfies certain critical exponential growth conditions. This result extends those in [MS2], which were proved for $p=n.$ 
 
Section \ref{hsbg} is devoted to the Heisenberg group, where we obtain, as an application of Theorems 1 and 2, the sharp MSI  for the horizontal gradient and for the powers of the subLaplacian, both of which were previously unknown.  We also obtain the sharp Ruf MTI for such operators, which were previously known only  for the horizontal gradient [LL1]. 

Finally, in Section \ref{mf} we apply Theorem 1 and obtain sharp MSI on complete, simply connected Riemannian manifolds with sectional curvature bounded above and below by two negative constants. In the hyperbolic space the case $\alpha=1$ was settled in [LT], and the case $\alpha=2$ in [NN]. Recently, Bertrand and Sandeep [BS] proved that the classical MTI in such manifold holds, extending previous results on the Hyperbolic space by Fontana and Morpurgo in [FM3]. Both papers exploit the fact that the Green kernel for $\nabla^\a$ on simply connected manifolds with constant and negative curvatures have a nice exponential decay, and satisfy a critical integrability condition at infinity.

\medskip

Some  results in this paper are included in Qin's PhD thesis [Qin1], specifically: 1) 
 the Adams MSI with varying powers in the denominator, for Riesz-like convolution potentials on $\R^n$;  2)  the application of  \eqref{MS3} to the quasilinear equation  \eqref{qlequation}; 3) the MSI inequality on $\H^n$, in particular \eqref{hs2b} of Theorem 5,  and Theorem 6 below (but without varying powers in the denominator). The latter results were obtained by adapting the proof in $\R^n$ in the context of the Heisenberg group.

\bigskip
\medskip\noindent{\bf Acknowledgment.} The authors would like to thank Luigi Fontana, Stefano Pigola, and Giona Veronelli  (Universit\`a di Milano-Bicocca) for helpful discussions and comments.

\bigskip

\section{Adams inequalities on metric measure spaces}\label{main}

Let $(M,d,\mu)$ be a metric measure space, that is a set $M$ endowed with a distance function $d$ and a Borel measure $\mu$.
 Define 
$$B(x,r)=\{y:\; d(x,y)\le r\},\qquad V_x(r)=\mu\big(B(x,r)\big).$$

We assume that for all $x\in M$
$$i)\; \; \forall r>0,  V_x(r)<\infty,\qquad ii)  \; r\to V_x(r) \;{\rm continuous.}$$
Such hypothesis easily implies that the boundary of any ball has zero measure, and that for any  $x\in M$ and any  $E$ measurable the function $r\to \mu\big(B(x,r)\cap E\big)$ is continuous on $[0,\infty)$. In particular, the measure $\mu$ is nonatomic.

A measurable function $k:M\times M\to \R$ will be called  a  {\it Riesz-like kernel of order $\beta>1$ and normalization constants $A_0>0$ and $A_\infty\ge0$} if there exists $B\ge0$ such that for all $x\in M$

\be\mathop\int\limits_{r_1\le d(x,y)\le r_2} |k(x,y)|^\b d\mu(y)\le  A_0\log{\frac{V_x(r_2)}{V_x(r_1)}}+B,  \quad  0<V_x(r_1)<V_x(r_2)\le1,\label{k1}\tag{K1}\ee
 
  \be\mathop\int\limits_{r_1\le d(x,y)\le r_2} |k(x,y)|^\b d\mu(y)\le  A_\infty\log{\frac{V_x(r_2)}{V_x(r_1)}}+B,  \quad  1\le V_x(r_1)<V_x(r_2)\label{k2}\tag{K2}\ee

for all $x,y\in M$
\be |k(x,y)|\le BV_x\big(d(x,y)\big)^{-1/\b}\qquad |k(x,y)|\le B V_y\big(d(x,y)\big)^{-1/\b}\label{k3}\tag{K3}\ee

\medskip
for each $\delta>0$  there is $B_\delta>0$ such that for all $x\in M$

\be\mathop\int\limits_{ d(x,y)> R} |k(x',y)-k(x,y)|^\b d\mu(y)\le B_\delta,\quad V_x(R)\ge (1+\delta) V_x(r),\quad \forall x'\in B(x,r).\label{k4}\tag{K4}\ee

The classical Riesz kernel on $\R^n$, i.e. $|x-y|^{\a-n}$, satisfies (K1)-(K4) with $\b=\frac{n}{n-\a}$ and $A_0=A_\infty=|B_1|$, the volume of the unit ball. More generally, the convolution kernels treated in [FM2] or [Qin2]  are of the above type.
When $\beta=1$ kernels satisfying conditions similar to  \eqref{k3} and \eqref{k4} are known as ``standard kernels" in the literature involving singular integrals (see e.g. [Ch]).

It's easy to check that the continuity hypothesis on $V_x$ implies that 
\be \left(\frac{1}{V_x(d(x,\cdot))}\right)^*(t)=\frac{1}{t}\qquad x\in M,\;  t>0\label {vstar}\ee
where the function on the left is the decreasing rearrangement of $V_x(d(x,\cdot))\big)^{-1}$ for each fixed $x$ (for the definition of decreasing rearrangement see \eqref{star}  below).  Using \eqref{vstar} we get that  \eqref{k3} implies 
\be\mathop\int\limits_{r_1\le d(x,y)\le r_2} |k(x,y)|^\b d\mu(y)\le B^\b \log{\frac{V_x(r_2)}{V_x(r_1)}},\qquad 0<V(r_1)<V(r_2)<\infty.\label {1}\ee
 In essence, conditions \eqref{k1} and \eqref{k2} quantify the asymptotic behavior of $k(x,y)$ when $V_x(d(x,y))$ is either small or  large, in terms of  specific  constants $A_0, A_\infty$.
 
 \smallskip
 The regularity condition \eqref{k4}  is a consequence of the following pointwise regularity estimate
  \be |k(x',y)-k(x,y)|\le C_\delta V_{x}(d(x,x'))^\eta\,  V_{x}(d(x,y))^{-\eta-1/\b}\label{xreg}\ee
   valid for some $\eta>0$, and for all $x\in M$, all $x'\in B(x,r)$, and all $y\notin B(x,R)$, whenever $V_{x}(R)\ge (1+\delta)V_{x}(r)$.

 We also  note that conditions \eqref{k1},\eqref{k2} can be expressed in terms of the decreasing rearrangement of $k$. For $t>0$, let 
\be k_1^*(t)=\sup_{x\in M}\big(k(x,\cdot)\big)^*(t),\qquad
 k_2^*(t)=\sup_{y\in M}\big(k(\cdot,y)\big)^*(t)\label{k12star}\ee
where $\big(k(x,\cdot)\big)^*(t)$ and $\big(k(\cdot,y)\big)^*(t)$ are the decreasing rearrangements of $k(x,y)$ for fixed $x$ and fixed $y$, respectively.

The first observation is that    condition \eqref{k3} implies 
 \be  k_1^*(t)\le C t^{-1/\b},\qquad k_2^*(t)\le Ct^{-1/\b},\qquad t>0\label{kk12}\tag{K3'}\ee 
 (with $C=B^{1/\b}$). It is also true (and it will be used later) that if $A_\infty=0$ and $k$ is symmetric, then the validity of \eqref{k3} for  $V_x(d(x,y))\le 1$ is enough to imply  \eqref{kk12} (see Remark \ref{k1prime} in Appendix A).

Secondly, one can see that assuming \eqref{k3} conditions  \eqref{k1}, \eqref{k2}  are equivalent to

\be\mathop\int_{t_1}^{t_2} \big(k^*(x,u)\big)^\b du\le  A_0\log{\frac{t_2}{t_1}}+B,  \quad  0<t_1<t_2\le1,\;\;\forall x\in M\label{k1p}\tag{K1'}\ee
 
  \be\mathop\int_{t_1}^{t_2} \big(k^*(x,u)\big)^\b du\le  A_\infty\log{\frac{t_2}{t_1}}+B,  \quad  1\le t_1<t_2,\;\;\forall x\in M\label{k2p}\tag{K2'}\ee
  where for simplicity we let 
  \be k^*(x,t)=\big(k(x,\cdot)\big)^*(t)\ee
(see Appendix A for a proof). 

In turn, the integral conditions above are implied by following pointwise asymptotic conditions on $k_1^*$:

 \be k_1^*(t)\le A_0^{1/\b}t^{-1/\b}+C(t^{-1/\b+\delta_1})\qquad 0<t\le 1\label{k1pp}\tag{K1''}\ee
 
 \be k_1^*(t)\le A_\infty^{1/\b}t^{-1/\b}+C(t^{-1/\b-\delta_2})\qquad t\ge 1\label{k2pp}\tag{K2''}\ee
for some $\delta_1,\delta_2>0$. Condition \eqref{k1pp}, together with  the second condition in \eqref{kk12}, was first introduced in \cite{fm1}, as a main hypothesis  in order to derive Adams inequalities on measure spaces of finite measure, for potentials with kernel $k$. We note here, and this will be also used later, that the main result in [FM1], Theorem 1, is true under the more general conditions \eqref{k1} and \eqref{k3}, or  \eqref{k1p} and \eqref{kk12}.

The special case $A_\infty=0$ in (K2') corresponds to what in [FM3] was labeled as a ``critical integrability" condition, namely 
\be\supess_{x\in M}\int_1^\infty \big( k^*(x,t)\big)^\b dt<\infty.\label{crit}\ee
It was proven in [FM3] that the above condition  together with (K1'') and (K3') guarantee a classical Adams inequality with exponential constant $1/A_0$, on general measure spaces. Even in this case the proof goes though if one only assumes \eqref{k1p} and (K3').

A Riesz-like potential is an integral operator
 \be Tf(x)=\int_M k(x,y)f(y)d\mu(y).\label{de3}\ee
 where $k(x,y)$ is  Riesz-like. Such operator is well defined for $f\in L^{\b'}$ and with compact support.  Here and from now on $\b'$ will denote the exponent conjugate to $\b$:
 $$\frac {1}\b+\frac 1{\b'}=1.$$

 \begin{theorem}\label{m1} If $k(x,y)$ is a Riesz-like kernel of order $\b>1$ and normalization  constants $A_0, A_\infty$, and if  
$p\ge 1$, 
 then there is a constant $C$ such that for any measurable function $f$ supported in a ball with 

 \be\|f\|_{\b'}\le 1\label{1b}\ee we have
 
 \be \int_M \frac{\exp_{\llceil p/\b-1\rrceil}\bigg[\dfrac{1} {\max\{{A_0},{A_\infty}\}}|Tf|^{\b}\bigg]}{1+|Tf|^{p\b/\b'}}\,d\mu(x)\le C\|Tf\|_{p}^{p}. \label{1a}\ee

\def\k{\kappa}
Moreover, if $A_\infty>0$ and given any  $\k>0$ there exists $C$ (depending on $\k$) such that for any measurable  $f$ supported in a ball with 
\be \|f\|_{\b'}^{\b'}+\k \|Tf\|_{\b'}^{\b'}\le 1\label{1ba}\ee 
we have 
 \be \int_M \frac{\exp_{\llceil p/\b-1\rrceil}\bigg[\dfrac1 {A_0}|Tf|^{\b}\bigg]}{1+|Tf|^{p\b/\b'}}\,d\mu(x)\le C\|Tf\|_{p}^{p}, \label{1a1}\ee
and
 \be \int_M \exp_{\llceil \b'-2\rrceil}\bigg[\dfrac1 {A_0}|Tf|^{\b}\bigg]\,d\mu(x)\le  {C}. \label{1a2}\ee

\ni If $A_\infty=0$ (critical integrability) then \eqref{1a1}, \eqref{1a2} holds under \eqref{1b}, and the same is true if the first estimate in  (K3) holds  for $V_x\big(d(x,y)\big)\le 1$, (K4) holds for $V_x(r)\le 1$, and  \eqref{kk12} holds.

 \end{theorem}

    \begin{rk}\label{rk1} If $M=\R^n$ $\mu=$Lebesgue measure and $k(x,y)=|x-y|^{\a-n}$, the Riesz kernel, then it is easy to see by dilation that the constant $C$ in \eqref{1a1} is independent of $\kappa$, whereas in \eqref{1a2} it is of type $C/\kappa$ for any $\k>0$. In general, it is possible to see   that the current proof only yields this behavior for large $\kappa$.\end{rk}

 \begin{rk}\label{rkvector} The conclusions of Theorem 1 hold if $k$ and $f$ are vector-valued. Specifically, if $k=(k_1,...,k_m)$ is defined on $M\times M$ and measurable, and if $f=(f_1,...,f_m)$ is measurable on $M$, then let $Tf(x)=\int_M k(x,y)\cdot f(y)d\mu$ where the $``\,\cdot\, ”$ denotes the standard Euclidean inner product on $\R^m.$ If $k$ satisfies \eqref{k1}-\eqref{k4}, and if $|k|=(k\cdot k)^{1/2},\ |f|=(f\cdot f)^{1/2}$, $\;\|f\|_{\b'}=\big(\int_M|f|^{\b'}\big)^{1/\b'}$,   then the conclusions of Theorem 1 hold (just apply
 the inequality $|k\cdot f|\le |k||f|$ in \eqref{mtest2}, \eqref{fla8}, \eqref{la3}, \eqref{lem2}, \eqref{lc4}, and \eqref{lem6} in the proof).\end{rk}
 
 % to the real-valued operator $\int_M |k(x,y)|\,|f(y)|d\mu\ge |Tf(x)|$).\end{rk}
    \begin{rk}\label{AA} If $A_\infty\le A_0$, then the exponential constants in \eqref{1a} and \eqref{1a2} are the same, and it would be possible to deduce \eqref{1a2} from \eqref{1a} in this case,  using a method similar to the one used in [Qin2, Proof of Corollary 3]. However, if $A_\infty>A_0$ this cannot be done.
    \end{rk}

\begin{rk}\label{modified} We will see in Section \ref{mf} an instance where  $A_\infty=0$ but the estimates in \eqref{k3} is only available for balls with volumes less than or equal 1.

\end{rk}

\begin{rk}\label{finite measure}  If $\mu(M)<\infty$ then condition (K3) implies condition (K2) with $A_\infty=0$. Hence, if $k$ is a Riesz-like kernel satisfying (K1), (K3) and (K4) on a space with finite measure, we can deduce \eqref{1a1}, \eqref{1a2} under $\|f\|_{\b'}\le 1$. Note that the validity of \eqref{1a2} under critical integrability is guaranteed  by results in [FM3].
\end{rk}

The next theorem deals with the sharpness of the exponential constants in Theorem 1. The essence is that if \eqref{k1} or \eqref{k2} are sharp at some $x_0$, then  
under further regularity conditions on $k$ the exponential constants in Theorem 1 are sharp. 

We say that a Riesz-like kernel with normalization constants $A_0, A_\infty$ is {\it proper on the diagonal} if
 there exists $x_0\in M$ such that 
 \be\mathop\int\limits_{r_1\le d(x_0,y)\le r_2} |k(x_0,y)|^\b d\mu(y)\geq  A_0\log{\frac{V_{x_0}(r_2)}{V_{x_0}(r_1)}}-C,  \quad  0<V_{x_0}(r_1)<V_{x_0}(r_2)\le1,\label{k5}\ee
and that is {\it proper at infinity} if $\mu(M)=+\infty$, $A_\infty>0$, and there is $x_0\in M$ such that  
 \be\mathop\int\limits_{r_1\le d(x_0,y)\le r_2} |k(x_0,y)|^\b d\mu(y)\geq  A_\infty\log{\frac{V_{x_0}(r_2)}{V_{x_0}(r_1)}}-C,  \quad  1\le V_{x_0}(r_1)<V_{x_0}(r_2).\label{k7}\ee
 
 It is easy to check that, due to the $x-$regularity in the form \eqref{k4}, for any given $r_1, r_2>0$ if \eqref{k5}  holds, then  there is $r'<r_1$ (depending on $r_1,r_2$) so that  \eqref{k5} holds with $x_0$ replaced by any $x\in B(x_0,r')$, and the same is true for \eqref{k7}.

 We will also require a slightly stronger  $x-$regularity condition than \eqref{k4}, and an even higher  regularity of $k(x,y)$ in the $y$ variable. For simplicity we will assume pointwise regularity in the same spirit as \eqref{xreg}.

 We will say that a Riesz-like kernel has a {\it Taylor formula with volume remainder  of order $\eta\ge 0$ in the $y$ variable at $x_0$} if there is an integer $m\ge1$ and measurable functions $k_j(x,x_0), p_j(y,x_0)$, $j=0,1,...,m-1$, where the $\{p_j\}_0^{m-1}$ are bounded and linearly  independent on $B(x_0,r)$ for any $r>0$,  and such that

 \be \bigg|k(x,y)-\sum_{j=0}^{m-1} k_j(x,x_0)p_j(y,x_0)\bigg|\le C_\delta V_{x_0}(d(y,x_0))^{\eta} \; V_{x_0}(d(x,x_0))^{-\eta-1/\b}\label{yreg2}\ee
  whenever $V_x(R)\ge (1+\delta)V_x(r)$, $\;d(y,x_0)\le r,\;d(x,x_0)\ge R$.
 We will also require that if $\{v_j(y,x_0)\}_0^{m-1}$ is an orthonormal basis of the space spanned by the $p_j$ restricted to the ball $B(x_0,r)$,  then 
 \be|v_j(y,x_0)|\le\frac {C}{ \sqrt{V_{x_0}(r)}},\qquad y\in B(x_0,r),\qquad j=0,1,..,m-1.\label{yreg3}\ee
  
In concrete cases where $k(x,y)$ is a differentiable function one can typically  choose the $p_j$ so that for $j=0,1,...,m-1$
\be |k_j(x,x_0)|\le C_\delta  V_{x_0}(d(x,x_0))^{-\eta_j-1/\b},\qquad |p_j(y,x_0)|\le C_\delta V_{x_0}(d(y,x_0))^{\eta_j}\ee
with $\eta_0=0$, $k_0=p_0=0$ (which is (K3)), and $k_1(x,x_0)=k(x,x_0)$, $p_1(y,y_0)=1,$ which gives the same regularity condition as in \eqref{xreg} but in the $y$ variable.

The Riesz kernel in $\R^n$ and the $m-$regular Riesz-like kernels defined in [FM2], and [Qin2] are obvious examples of kernels of this type, being approximated by their classical Taylor polynomials of order $m-1$, with  volume remainders of order  $m/n$ (see for example eq. (95) in [FM2]).

  \begin{theorem}\label{s} Under the same assumptions of Theorem 1, assume further that there is $x_0$ such that $k(x,y)$ satisfies \eqref{xreg} with  $x=x_0$, for some $\eta>0$, and has a Taylor formula with volume remainder of order $\eta\ge0 $ in the $y$ variable at $x_0$. If $p>(1+p/\b')(1+\eta)^{-1}$ then the following hold:
 \smallskip
 
  \ni (a) If $k(x,y)$  is proper on the diagonal, then   the exponential constant $A_0^{-1}$  in \eqref{1a1}, \eqref{1a2} is sharp, and it is also sharp in \eqref{1a}  when $A_0\ge A_\infty$. The power of the denominators $p\b/\b'$ in \eqref{1a} and \eqref{1a1} is also sharp.
 
 \smallskip
 
 \ni (b)  If $k(x,y)$  is proper at infinity,  then the exponential constant $A_\infty^{-1}$ in \eqref{1a} is sharp when $A_0\le A_\infty$. The power of the denominator $p\b/\b'$ in \eqref{1a} is also sharp. 
 
 \end{theorem}
 
 Here the meaning of ``sharp'' is to be intended as follows:  the ratios of the left-hand sides and right-hand sides of \eqref{1a}, \eqref{1a1}, \eqref{1a2} can be made arbitrarily large if either the exponential constants are larger or the powers of the denominators are smaller. There is only one case where we can say more, namely in case (a) the left-hand sides of   \eqref{1a}, \eqref{1a1} can be arbitrarily large while the right-hand sides stay small, if the exponential constants are larger (see Proof of Theorem 2).\medskip
 
 \ni{\bf Notes.} 1. The conclusions of the above theorem remain valid under slightly weaker regularity conditions, namely \eqref{k6}, \eqref{k8}.
 
 \ni2. If $p>1+p/\b'$ then we can take $m=0$, $\eta_0=0$ and just use (K3) to derive sharpness. If one can choose $m$ high enough with  $\eta\ge p/\b'$, then sharpness holds for all $\b,p>1$. For example, for the Riesz kernel  in $\R^n$, and $p=\b'$ then $\eta=1$ is obtained by taking the Taylor expansion of  order $n-1$ 
of  $|x-y|^{\a-n}$.

%#########################################################

%\input{3-thm1-pf}

\section{Proof of Theorem \ref{m1}}\label{pfthm1}
We will prove Theorem 1 under the hypothesis that (K1)-(K4) hold and $A_\infty>0$. In Section~\ref{Ainfinity} we will outline the changes of the proof in order to deal  with the case   $A_\infty=0$, under the  weaker conditions stated at the end of Theorem 1.  We will also assume throughout the proof, WLOG,  that $B\ge 2^{1/\b'}.$

Let us recall the definition of decreasing rearrangement. Given a measurable function $f\ :\ M\rightarrow \ [-\infty,\infty]$ its distribution function is defined as
\be m_f(s)=\mu(\{x\in M:\ |f(x)|>s\}),\ \ \ s\geq 0.\nt\ee  
Assuming that the distribution function of $f$ is finite for all $s>0 $, the 
 decreasing rearrangement of $f$ is defined as
\be f^*(t)=\inf\{s\geq0:\ m(f,s)\leq t\},\ \ \ t>0,\label{star}\ee
and we also define
\be f^{**}(t)=\frac{1}{t}\int_0^tf^*(u)du,\ \ \ t>0.\nt\ee

\par
From now on let us assume that 
\be  \|f\|_{\b'}\le 1,\qquad \supp f:=\{x:\;f(x)\neq0\}\subseteq B(x_0,R),\label{hypf}\ee
for some $x_0,\in M$, $R>0$ (depending on $f$) and after possibly redefining $f$ to be 0 on a set of 0 measure outside $B(x_0,R)$.

Under these circumstances $T f(x)$ is well-defined and finite for a.e. $x$ in $M$. This follows from the weak type estimate
\be \big(m_{T'|f|}(s)\big)^{1/\beta}\le \frac{\b' B\|f\|_1}{s}\ee
where $T'$ is the operator with kernel $|k(x,y)|$ (see \cite{qin2}, estimate (A.1)).

\smallskip 
By the ``Exponential Regularization Lemma''  (see \cite{fm2}, Lemma 9, \cite{qin2}, Lemma A)  \eqref{1a} and \eqref{1a1} are equivalent to  
 \be \int_{|Tf|\geq 1} \frac{\exp{\bigg[\dfrac{1}{A}|Tf|^{\b}\bigg]}}{1+|Tf|^{p\b/\b'}}d\mu(x)\le C\|Tf\|_{p}^{p} \label{1aa}\ee
 and \eqref{1a2} is equivalent to 
 \be \int_{\{|Tf|\geq 1\}}\exp\bigg[\frac{1}{A_0}|Tf|^{\b}\bigg]d\mu\leq C,  \label{T4a}\ee
 where
 \be A=\bc \max\{A_0,A_\infty\}\ \ \ &\text{given}\ \ \eqref{1b}\cr A_0 &\text{given}\ \ \eqref{1ba}.\ec\label{adefi}\ee
% for part (1) and 
%  \be  \int_{|Tf|\geq 1} \frac{\exp{\bigg[\dfrac{1}{A_1}|Tf|^{\b}\bigg]}}{1+|Tf|^{\si\b}}d\mu(x)\le C||Tf||_{\si\b'}^{\si\b'} \label{1a1a}\ee
%  for part (2).\par
  Let \be t_0=\mu\{x\in M\ :\ |Tf(x)|\geq 1\}\label{t0}\ee so that \eqref{1aa} is equivalent to 
 \be \int_{0}^{t_0}\frac{\exp\bigg[\dfrac{1}{A}\big((Tf)^*(t)\big)^{\b}\bigg]}{1+\big((Tf)^*(t)\big)^{p\b/\b'}}dt \leq C\|Tf\|_{p}^{p} ,\label{1dd1}\ee
 and \eqref{T4a} is equivalent to 
  \be \int_{0}^{t_0}{\exp\bigg[\dfrac{1}{A_0}\big((Tf)^*(t)\big)^{\b}\bigg]}dt \leq C.\label{rufnew8}\ee
%  and \eqref{1a1a} is equivalent to 
% \be \int_{0}^{t_0}\frac{\exp\bigg[\dfrac{1}{A_1}\big((Tf)^*(t)\big)^{\b}\bigg]}{1+\big((Tf)^*(t)\big)^{\si\b}}dt \leq C||Tf||_{\si\b'}^{\si\b'}. \label{1dd2}\ee
 
%By looking at the exponential constants in the above inequalities, it is clear that WLOG we can assume $A_0,\ A_\infty>0$. Also in part (2) it is enough to assume that $A_\infty>A_0$, since otherwise it can be reduced to the case in part (1) and hence \eqref{1a1} follows.
% By remark 1 we only need to consider the case $\mu(M)=\infty$. \par
A first  estimate for $(Tf)^*(t)$ is given via the {\it O'Neil functional},  defined as 
\be Uf(t)=Ct^{-\frac{1}{\b}} \int_0^{t}{f}^*(u)du+\int_t^{\infty}k_1^*(u){f}^*(u)du.  \nt\ee
The O'Neil Lemma on measure spaces  (see \cite{fm1}, \cite{fm3}) gives $(Tf)^*(t)\le (Tf)^{**}(t)\le Uf(t)$ for any $t>0$. 

The basic Adams inequality on measure spaces with finite measure can be stated in terms of $Uf$ as follows: 
\be\int_0^\tau\exp\bigg[\frac1{A_0} \big(Uf(t)\big)^\b\bigg]dt\le C\big(\tau+\mu(\supp f)\big),\qquad \tau>0\label{Adams}\ee
under conditions \eqref{k1}, \eqref{k3}, where $C$ is independent of $f, T$. Such an inequality has been stated and used in [Qin2] (see Thm.$\;\!\!$E) under  more restrictive pointwise conditions such as \eqref{k1pp}, \eqref{kk12}. Its proof is basically contained in the proof of [FM3, Corollary 2], and the proof under the more general conditions \eqref{k1},\; \eqref{k3} follows by modifying the function $g(x,\xi,\eta)$ in [FM3, p. \hskip-.2em10], letting it equal to $k_\tau(x,e^{-\xi})e^{-{\xi\over\beta}}$, when $\xi\in (-\infty,\eta)$, and by applying the integral condition in the form \eqref{k1p} to obtain the estimate on [FM3 p. 29, line 3].

\smallskip
The  key starting point  in \cite{qin2}, in the case of convolution Riesz-like kernels in $\R^n$, was to split $f$ on a suitable level set of given measure $\tau$.

To this end, for each $\tau>0$ we consider a measurable  set $F_\tau$ satisfying
\be
\bc
\{x\ :\ |f(x)|>f^*(\tau)\}\subseteq F_{\tau}\subseteq \{x\ :\ |f(x)|\geq f^*(\tau)\}\cr \mu(F_{\tau})=\tau
  \ec\label{d3}\ee
%\be \mu(E_{\tau})=\tau,\qquad
%\{x\ :\ |Tf(x)|>(Tf)^*(\tau)\}\subseteq E_{\tau}\subseteq \{x\ :\ |Tf(x)|\geq(Tf)^*(\tau)\}
 %\label{d2}\ee 

\ni the existence of which is guaranteed by  the continuity assumption of  $\mu(B(x,r))$ in $r$, for all $x\in M$. Such a set may not be unique, for example if $f^*$ is constant around $0$, however a prescription can be given to identify it uniquely,  once a preferred point $x_0\in M$ is fixed. Namely (see \cite{qin2}) if $V_1=\{x\ :\ |f(x)|>f^*(\tau)\}$ and $V_2=\{x\ :\ |f(x)|\ge f^*(\tau)\}$, then $F_\tau$ can be found of type $V_1\cup \big(B(x_0,r)\cap(V_2\setminus V_1)\big)$, for some $r\in[0,\infty]$, uniquely up to a set of zero measure.
  Note that in $\R^n$ if $f$ is radially decreasing and $\tau=|B(0,r)|$ then $F_\tau$ either the open or the closed ball of center 0 and radius~$\,r.$
  
  The set $F_\tau$ with $\mu(F_\tau)=\tau$ satisfies the identity
  \be\int_{F_\tau}\Phi(|f(x)|)d\mu(x)=\int_0^\tau \Phi\big(f^*(u)\big)du\label{etauid}
  \ee
  where $\Phi$ is any non-negative, measurable function on $[0,\infty)$.  
See [BeSh, Lemma 2.,5], for the existence of such set on finite nonatomic measure spaces, and where the above identity is given for $\Phi(u)=u$, and which  can be easily extended to arbitrary $\Phi$ (see also [Qin2, (3.10)]).

\medskip 
  Likewise, we define a corresponding level set $E_\tau$  relative to $Tf$:
  
 \be \bc 
\{x\ :\ |Tf(x)|>(Tf)^*(\tau)\}\subseteq E_{\tau}\subseteq \{x\ :\ |Tf(x)|\geq(Tf)^*(\tau)\}\\ \mu(E_{\tau})=\tau.\qquad\ec
 \label{d2}\ee 

Let \be f_{\tau}=f\chi^{}_{F_{\tau}},\ \ \ f_\tau'=f\chi^{}_{F^c_{\tau}}=f-f_\tau \nt\ee
\ni and let 
\be r_x(\tau)=\min\{r>0: V_x(r)=\tau\}.\label{rtau}\ee 

We will use for simplicity the notation
\be B_\tau(x)=B(x,r_x(\tau))\label {Btau}\ee
to denote the smallest ball centered at $x$ and with volume $\tau$.

As in [Qin2], the crucial quantity we will use to  control  $f$ outside $F_\tau$ is defined as follows:

%\be W_\tau=\supess_{x\in E_\tau}\int_\tau^{2\tau}k_1^*(u)(f_\tau'\chi^{}_{B(x,r_x(\tau))})^*(u-\tau)du \label{d7} \ee
\be \Tfc(\tau):=\supess_{x\in E_\tau}|T(f_\tau'\chi^{}_{B_\tau(x)^c})(x)|=\supess_{x\in E_\tau}\bigg|\int_{B_\tau(x)^c\cap F_\tau^c } k(x,y)f(y)d\mu(y)\bigg|, \label{d8}\ee

which is finite for every $\tau$ from H\"{o}lder's inequality, \eqref{k3}, and \eqref{hypf}, which  imply
\be \Tfc(\tau)\le \Big(\frac{BV_{x_0}(R)}{\tau}\Big)^{1/\b}.\label{circ0}\ee

\begin {lemma} There is $z\in E_\tau$ such that for $0<t\le\tau$ we have 
\be (Tf)^{*}(t)\leq (Tf)^{**}(t)\leq Uf_\tau(t)+B\|f_\tau'\chi_{B_\tau(z)}\|_{\b'}+\Tfc(\tau).\label{e1}\ee
\end{lemma}

\ni{\bf Proof.}  Arguing exactly as in [Qin2], using the improved O'Neil Lemma ([Qin2, Lemma 1])  we obtain 
\be (Tf)^{**}(t)\leq Uf_\tau(t)+\supess_{x\in E_\tau}\int_\tau^{2\tau}k_1^*(u)(f_\tau'\chi^{}_{B_\tau(x)})^*(u-\tau)du+\Tfc(\tau).\label{e2}\ee

Note that by \eqref{kk12} we have, for each $x\in E_\tau$,
 \be  \ba &\int_\tau^{2\tau}k_1^*(u)(f_\tau'\chi^{}_{B_\tau(x)})^*(u-\tau)du
  \leq B^{1/\b}\int_\tau^{2\tau}u^{-\frac{1}{\b}}(f_\tau'\chi^{}_{B_\tau(x)})^*(u-\tau)du \leq B^{1/\b}\|(f_\tau'\chi^{}_{B_\tau(x)})^*\|_{\b'}\cr&= B^{1/\b}\|f_\tau'\chi^{}_{B_\tau(x)}\|_{\b'}\le  B^{1/\b}\supess_{x\in E_\tau} \|f_\tau'\chi^{}_{B_\tau(x)}\|_{\b'}\le 2B^{1/\b}\|f_\tau'\chi^{}_{B_\tau(z)}\|_{\b'}\le B \|f_\tau'\chi^{}_{B_\tau(z)}\|_{\b'}\ea\ee
 for some $z\in E_\tau $ (note that $\|f_\tau'\chi^{}_{B_\tau(x)}\|_{\b'}$ is continuous in $x$, and that we also assumed $B\ge 2^{1/\b'}$).
\hfill\QEDopen\\

\bigskip

%Let $z\in E_\tau$ be such that 
%\be B_1\|f_\tau'\chi^{}_{B(z,r_x(\tau))}\|_{\b'}\geq \frac{1}{2}W_\tau,\label{w2}\ee
%and define 
%\be \a'=\|f_\tau'\chi^{}_{B(z,r_x(\tau))}\|_{\b'}.\label{w7}\ee Note that $\a'$ depends on $\tau,\ z$.

%We write \eqref{e1} as
%\be (Tf)^{*}(t)\leq (Tf)^{**}(t)\leq Uf_\tau(t)+2B_1\a'+\Tfc(\tau)\ \ \ \ \ \textup{for}\ 0<t\leq\tau,\label{w3}\ee
Using the above lemma we can estimate the integrand in \eqref{1dd1}. For   every $\e\in (0,1)$, using $(a+b)^\beta\le\e^{1-\b}a^\b+(1-\e)^{1-\b} b^\b$ we have
\be\ba 
&\frac{\exp\bigg[\dfrac{1}{A}\big((Tf)^*(t)\big)^{\b}\bigg]}{1+\big((Tf)^*(t)\big)^{p\b/\b'}}\le \cr&\quad\leq C\,{\exp\bigg[\dfrac{\e^{1-\b}}{A}\big(Uf_\tau(t)\big)^{\b}\bigg]}\cdot \frac{\exp\bigg[\dfrac{(1-\e)^{1-\b}}{A}\big(\Tfc(\tau)+B\|f_\tau'\chi_{B_\tau(z)}\|_{\b'}\big)^{\b}\bigg]}{1+\big(\Tfc(\tau)+B\|f_\tau'\chi_{B_\tau(z)}\|_{\b'}\big)^{p\b/\b'}}\ea
\label{main1}\ee
and
\be\ba 
&\exp\bigg[\dfrac{1}{A_0}\big((Tf)^*(t)\big)^{\b}\bigg]\le \cr&\quad\leq C\,{\exp\bigg[\dfrac{\e^{1-\b}}{A_0}\big(Uf_\tau(t)\big)^{\b}\bigg]}\cdot \exp\bigg[\dfrac{(1-\e)^{1-\b}}{A_0}\big(\Tfc(\tau)+B\|f_\tau'\chi_{B_\tau(z)}\|_{\b'}\big)^{\b}\bigg]\ea
\label{rufnew7}\ee
for some $C>0$ independent of $f,\tau,\e$.

\bigskip
The following is an easy consequence of \eqref{Adams}:

\begin{prop}\label{lemmaI21}
There is $C>0$ such  that for any $\e>0$ 
\be \int_0^\tau {\exp\bigg[\dfrac{\e^{1-\b}}{A_0}\big(Uf_\tau(t)\big)^{\b}\bigg]}\le C\tau,\qquad \tau>0 \label{c1bbb1}\ee 
provided that  \be {\|f_\tau\|_{\b'}^{\b'}}\le \e. \label{etau1}\ee 
\end{prop}
\ni {\bf Proof of Proposition \ref{lemmaI21}.} Let  \be \wt{f}:=\frac{f_{\tau}}{\e^{1/\b'}}, \nt\ee then we have that $\wt{f}$ has measure of support $\mu(\textup{supp} \wt{f})\leq \tau$ with \be\|\wt{f}\|_{\b'}\leq1.\nt\ee
%We also have 
%\be \frac{B_1\|f_\tau'\chi_{B_\tau(z)}\|_{\b'}}{\e^{1/\b'}}\le B_1\bigg(\frac{\e-\|f_\tau\|^{\b'}}{b\e}\bigg)^{1/{\b'}}=B_1b^{-1/\b'}(1-\|\wt{f}\|^{\b'}_{\b'})^{1/\b'}.\label{w4}\ee
 Therefore \eqref{c1bbb1} follows from the (improved) Adams inequality in the form \eqref{Adams}.\hfill\QEDopen\\
\par

\begin{prop}\label{le102}  There exists a constant $C^*>0$ such that for  $\kappa\ge0$, $\e\in[0,1)$ and for any measurable $f$ supported in a ball  satisfying 
\be \|f\|_{\b'}^{\b'}+\k \|Tf\|_{\b'}^{\b'}\le 1\label{1b1}\ee
\be \|f_\tau'\|_{\b'}^{\b'}+\kappa\|(Tf)\chi_{E_\tau^c}\|_{\b'}^{\b'}\le 1-\e \label{lambda11}\ee
\be \Tfc(\tau)> C^*\label{cstar}\ee
we have
\be \frac{\exp\bigg[\dfrac{(1-\e)^{1-\b}}{A}\big(\Tfc(\tau)+B\|f_\tau'\chi_{B_\tau(z)}\|_{\b'}\big)^{\b}\bigg]}{1+\big(\Tfc(\tau)+B\|f_\tau'\chi_{B_\tau(z)}\|_{\b'}\big)^{p\b/\b'}}\leq \frac{C}{\tau(1-\e)^{p\b/\b'}}\|(Tf)\chi^{}_{E_\tau^c}\|^{p}_{p}\qquad {\hbox{ if }} \;\kappa\geq 0, \label{l1a02} \ee
where $A$ is the same as in \eqref{adefi}, and 
\be {\exp\bigg[\dfrac{(1-\e)^{1-\b}}{A_0}\big(\Tfc(\tau)+B\|f_\tau'\chi_{B_\tau(z)}\|_{\b'}\big)^{\b}\bigg]}\le \frac{C}{\tau} \qquad {\hbox{ if }} \;\kappa>0,\label{rufnew2}\ee
where $C$ is a constant independent of $f,\e.$

\end{prop}

\begin{rk}\label{E}
The set $E_\tau$ in \eqref{lambda11} and \eqref{l1a02} can also be replaced by any measurable set $E$ with measure less than or equal $\tau$. 
\end{rk}

\smallskip

\begin{rk}\label{E1}
An immediate consequence of  \eqref{l1a02} is that under $\|f\|_{\b'}\le 1$, $\|f_\tau'\|_{\b'}\le 1$ and $\Tfc(\tau)> C^*$ we have the following form of the ``optimal descending growth condition'':

\be \frac{\exp\bigg[\dfrac{1}{\max\{A_0,A_\infty\}}\big(\Tfc(\tau)\big)^{\b}\bigg]}{1+\big(\Tfc(\tau)\big)^{\b}}\leq \frac{C}{\tau} \|(Tf)\chi^{}_{E_\tau^c}\|^{\b'}_{\b'}, \label{odgc}\ee
which includes the corresponding condition  \eqref{odgc1} derived  in [IMN] and [MS2] in $\R^n$, when  $k(x,y)=c_2\nabla|x-y|^{2-n}$ and $f=\nabla u$, where $u  $ positive and radially decreasing,  $-\partial_\rho u$ decreasing, and $u(R)>1$. In this case we indeed have $(Tf)^\circ(\tau)=u^*(R)$, if $\tau=|B(0,R)|$ (see Appendix \ref{tf}).
\end{rk}

\medskip

 The proof of Proposition \ref{le102} is quite involved and will be given in  Section \ref{lemma2}. Assuming Proposition 2 we will now derive the main inequalities of Theorem 1.
 
 \bigskip 
 
\ni {\underline{\bf {Proof of \eqref{1a} and \eqref{1a1}}}.}  For notational convenience let
\be  I_1(\tau,\e)=\frac{\exp\bigg[\dfrac{(1-\e)^{1-\b}}{A}\big(\Tfc(\tau)+B\|f_\tau'\chi_{B_\tau(z)}\|_{\b'}\big)^{\b}\bigg]}{1+\big(\Tfc(\tau)+B\|f_\tau'\chi_{B_\tau(z)}\|_{\b'}\big)^{\b}\big)^{p\b/\b'}}\ee
\be I_2(\tau,t,\e)= {\exp\bigg[\dfrac{\e^{1-\b}}{A_0}\big(Uf_\tau(t)\big)^{\b}\bigg]}\ee

First note that  if 
$ \Tfc(\tau)\le C^*$, 
then by direct computation and the fact that $(Tf)^*(t)\geq~1$ for $t<t_0$ we have 
\be I_1(\tau,\e)\le {e^{{(1-\e)^{1-\b}}{A^{-1}}(C^*+B)^{\b}}}{\tau^{-1}}\|Tf\|^{p}_{p}. \label{lambda3}\ee
\smallskip

Let
\be\e_\tau=\min\bigg\{ \dfrac{\|f_\tau\|_{\b'}^{\b'}}{\|f\|_{\b'}^{\b'}+\kappa\|Tf\|_{\b'}^{\b'}},\;\dfrac{1}{4}\bigg\}\label{dtau1}, \ee
then 
\be3/4\le 1-\e_\tau<1.\label{w6} \ee
For any $\kappa\ge 0$, using condition  \eqref{1b1} and \eqref{dtau1} we get

\be\ba &\|f_\tau'\|_{\b'}^{\b'}+\kappa\|(Tf)\chi^{}_{E_\tau^c}\|_{\b'}^{\b'}\le  \|f\|_{\b'}^{\b'}-\|f_\tau\|_{\b'}^{\b'}+\k\|Tf\|_{\b'}^{\b'}\cr
&=\big(\|f\|_{\b'}^{\b'}+\k\|Tf\|_{\b'}^{\b'}\big)\bigg(1-\frac{\|f_\tau\|_{\b'}^{\b'}}{\|f\|_{\b'}^{\b'}+\k\|Tf\|_{\b'}^{\b'}}\bigg) \le \big(\|f\|^{\b'}_{\b'}+\k\|Tf\|_{\b'}^{\b'}\big) (1-\e_\tau)\le 1-\e_\tau. \ea\label{lambda1}\ee
Taking $\e=\e_\tau$ in \eqref{etau1} and \eqref{lambda11},  the proof of \eqref{1a}, \eqref{1a1} follows in the same manner as in [Qin2, proof of (3.34)]. Let us define

\be\tau_0=\max\Big\{\tau\in[0,t_0]: \;\dfrac{\|f_\tau\|_{\b'}^{\b'}}{\|f\|_{\b'}^{\b'}+\kappa\|Tf\|_{\b'}^{\b'}}\le\frac{1}{4}\Big\}.\nt\ee
 It is clear that by the definition of $\e_\tau$ that $\e_{\tau_0}\ge \|f_\tau\|_{\b'}^{\b'}$, hence from Propositions 1, 2 and \eqref{lambda3} we get
\be I_1({\tau_0},\e_{{\tau_0}})\leq \frac{C}{{\tau_0}}\|Tf\|^{p}_{p}\qquad\text{and}\qquad \int_0^{{\tau_0}}I_2({\tau_0},t,\e_{{\tau_0}})dt \leq C{\tau_0}. \label{main4} \ee
Therefore, using \eqref{main1} it is immediate that
\be \int_{0}^{\tau_0}\frac{\exp\bigg[\dfrac{1}{A}\big((Tf)^*(t)\big)^{\b}\bigg]}{1+\big((Tf)^*(t)\big)^{p\b/\b'}}dt\leq\int_0^{\tau_0} I_1(\tau_0,\e_{\tau_0})I_2(\tau_0,t,\e_{\tau_0})dt\le C\|Tf\|^{p}_{p}.\label{main6}\ee

Next, we take $\tau=t$ for $\tau_0\le t\le t_0$, and $\e=1/8$ in $I_2$. Then by definition of O'Neil's operator and the fact that the support $f_t$ has measure less than or equal $t$
\be Uf_t(t)=C_0t^{-\frac{1}{\b}}\int_0^t f_t^*(u)du \leq C\bigg(\int_0^t (f_t^*)^{\b'}\bigg)^{1/\b'}\!\!=C\|f_t\|_{\b'}\leq C.\nt\ee
So we have \be I_2\bigg(t,t,\frac{1}{8}\bigg)\le C.\label{main702}\ee
Since $\tau_0\le t\le t_0$, by definition we have $\e_t=1/4$. Take $\theta=(\frac{1-1/4}{1-1/8})^{\b-1}<1$ we get
\be \ba I_1\bigg(t,\frac{1}{8}\bigg)&= \frac{\exp\bigg[ \dfrac{(7/8)^{1-\b} }{A}\big(\Tfc(t)+B\|f_t'\chi_{B_t(z)}\|_{\b'}\big)^{\b}\bigg]}{1+\big(\Tfc(t)+B\|f_t'\chi_{B_t(z)}\|_{\b'}\big)^{p\b/\b'}}\cr
&\le \left(\frac{\exp\bigg[ \dfrac{(3/4)^{1-\b} }{A}\big(\Tfc(t)+B\|f_t'\chi_{B_t(z)}\|_{\b'}\big)^{\b}\bigg]}{1+\big(\Tfc(t)+B\|f_t'\chi_{B_t(z)}\|_{\b'}\big)^{p\b/\b'}}\right)^{\!\!\theta}\cr 
&\le \left(\frac{\exp\bigg[ \dfrac{(1-\e_t)^{1-\b} }{A}\big(\Tfc(t)+B\|f_t'\chi_{B_t(z)}\|_{\b'}\big)^{\b}\bigg]}{1+\big(\Tfc(t)+B\|f_t'\chi_{B_t(z)}\|_{\b'}\big)^{p\b/\b'}}\right)^{\!\!\theta}
=I_1^{\theta}(t,\e_t)\le \frac{C}{t^\theta}\|Tf\|^{\theta p}_{p}.
\ea
\label{main802}\ee
Using \eqref{main702}, \eqref{main802} and that $t_0\le \|Tf\|_{p}^{p}$ (see \eqref{t0}), 
\be \ba
\int_{\tau_0}^{t_0}I_1\bigg(t,\frac{1}{8}\bigg)I_2\bigg(t,t,\frac{1}{8}\bigg)dt&\leq C\int_{\tau_0}^{t_0}\frac{1}{t^\theta}\|Tf\|^{\theta p}_{p}dt\le Ct_0^{1-\theta}\|Tf\|^{\theta p}_{p}\cr &\leq C\|Tf\|^{(1-\theta){p}}_{p}\|Tf\|^{\theta p}_{p}=C\|Tf\|_{p}^{p}.
\ea
 \label{c2e02}\ee
 Combining \eqref{main6} and \eqref{c2e02} we have inequality \eqref{1a1}.

\bigskip
\ni{\underline{\bf Proof of \eqref{1a2}}.}
The proof is almost the same as above. Let $k>0$ and define 
\be {I_3(\tau,\e)=\exp\bigg[\dfrac{(1-\e)^{1-\b}}{A_0}\big(\Tfc(\tau)+B\|f_\tau'\chi_{B_\tau(z)}\|_{\b'}\big)^{\b}\bigg]}\le \frac{C}{\tau},\ee
from \eqref{rufnew2}.

We take $\e_\tau$ as in \eqref{dtau1} and estimate similarly as in \eqref{main6}, \eqref{c2e02}. By using Lemmas \ref{lemmaI21}, \ref{le102} we get 
\be \int_{0}^{\tau_0}{\exp\bigg[\dfrac{1}{A_0}\big((Tf)^*(t)\big)^{\b}\bigg]}dt\leq\int_0^{\tau_0} I_3(\tau_0,\e_{\tau_0})I_2(\tau_0,t,\e_{\tau_0})dt\le C,\label{rufnew9}\ee
and 
\be \ba
\int_{\tau_0}^{t_0}I_3\bigg(t,\frac{1}{8}\bigg)I_2\bigg(t,t,\frac{1}{8}\bigg)dt\leq C\int_{\tau_0}^{t_0}\frac{1}{t^\theta}dt\le Ct_0^{1-\theta}\leq C\|Tf\|^{(1-\theta){\b'}}_{\b'}\le C.
\ea
 \label{rufnew10}\ee

\hfill\QEDopen\\

%###########################################################
%\input{4-prop2-pf}

\section{Proof of Proposition \ref{le102}}\label{lemma2} Let us assume throughout this section  that \eqref{1b1}, \eqref{lambda11} hold for some  $\kappa\ge0$, i.e. \be \|f\|_{\b'}^{\b'}+\k \|Tf\|_{\b'}^{\b'}\le 1,\qquad \|f_\tau'\|_{\b'}^{\b'}+\kappa\|(Tf)\chi_{E_\tau^c}\|_{\b'}^{\b'}\le 1-\e\nt\ee

\ni {\underline{\bf Proof of \eqref{l1a02}}.}

\bigskip

For $0<\tau<\mu(M)$ and $x\in E_\tau$ define 
%\be W_\tau(x)=\int_\tau^{2\tau}k_1^*(u)(f_\tau'\chi^{}_{B(x,r_x(\tau))})^*(u-\tau)du  \ee
\be \Tfc(\tau,x)=|T(f_\tau'\chi^{}_{B^c(x,r_x(\tau))})(x)|, \label{mtest1}\ee
so that $\Tfc(\tau)=\supess_{x\in E_\tau}\Tfc(\tau,x).$
Recall that $r_x(\tau)$ is the smallest $r$ such that $\mu(B(x,r))=\tau$.

\smallskip

To prove \eqref{l1a02} we will show that there are constants $C,C^*>0$ such that  for $0<\tau< t_0$ and for all $x\in E_\tau$ with $\Tfc(\tau,x)> C^*$
\be \frac{\exp\bigg[\dfrac{(1-\e)^{1-\b}}{A}\big(\Tfc(\tau,x)+B\|f_\tau'\chi_{B_\tau(z)}\|_{\b'}\big)^{\b}\bigg]}{1+\big(\Tfc(\tau,x)+B\|f_\tau'\chi_{B_\tau(z)}\|_{\b'}\big)^{p\b/\b'}}\leq C\frac{\|(Tf)\chi^{}_{E_\tau^c}\|_{p}^{p}}{\tau(1-\e)^{p\b/\b'}},\label{le1aaa}\ee
where $z\in E_\tau$ is as in Lemma 1.

Now let us state a more general version of the exact  growth condition for sequences given  in [IMN], [MS1]-[MS3], [LT], [LTZ], the proof of which is postponed to the Appendix.

Given any sequence $\ds a=\{a_k\}_{k\geq 0}$ we will let 
\be  \|a\|_{q}=\big(\aa |a_k|^{q}\big)^{1/q}\label{lo1}\ee 
and given a sequence $\ds{\lambda=\{\lambda_k\}_{k\ge0}}$ we will let
$\lambda a=\{\lambda_ka_k\}_{k\ge0}.$

\begin{lemma}\label{lo} For any $\b>1,\, p>0$ and given a sequence $\lambda=\{\lambda_k\}$ such that $\lambda_k\in(0,1]$ and $\lambda_k<1$ for at most  finitely many $k$, define, for $h>0$, 
 \be \mu_d(h)=\inf\Big \{\aa |a_k|^{p}e^{\b'k}: \ \|a\|_1=h,\ \|\lambda a\|_{\b'} \leq 1\Big\}.\label{muh1}\ee 
Then, for each $h_0>0$  there exist positive constants $C=C(p,\b,h_0),C'=C'(p,\b,h_0)$ such that \be Ce^{-L\b'}\,\frac{\exp\big[{\b'h^{\b}}\big]}{h^{p\b/\b'}}\leq \mu_d(h)\leq C' \frac{\exp\big[{\b'h^{\b}}\big]}{h^{p\beta/\b' }}\qquad \forall h\ge h_0.\label{lo2}\ee
where \be  L_\lambda:=\sum_{k=0}^\infty \Big(\lceil\lambda_k^{-\beta}\rceil-1\Big).\label {lambda}\ee
\end{lemma}

   \begin{rk}\label{rk3} The original lemma in [IMN],[MS1]-[MS3],[LT],[LTZ] (and which was also used in [Qin2])  had $\lambda_k=1$ for all $k$, and $p=\b'$. The original result for $\b'=2$ was given in [IMN] and [MS1], and was obtained for $\b'>2$ independently in [MS2] and [LT], while in [LTZ] the authors settled the case $1<\b'<2$. Our proof of Lemma \ref{lo}, given in the Appendix, is a bit shorter than the existing ones, and works for all $\b'>1$. \end{rk}
   
   \medskip
   %\ni  \ni{\bf Remark 3:} It is possible to prove a version the above result to the case where $\lambda_k<1$ for infinitely many $k$ and  $\sum_{0}^\infty \big(\lambda_k^{-\beta}-1\big)<\infty$.
   
\bigskip
From the above lemma we have that given a sequence $\lambda$ satisfying $L_\lambda\le C$, and given   any $\b'>1, \,\mu>0,\,h>1$, there is $C$ such that 
for any sequence $\{a_k\}$ satisfying 
\be\sum_{k=0}^\infty |a_k|=h \qquad \aa \lambda^{\b'}_{k}|a_k|^{\b'}\le \mu \label{le1a01}\ee we have
\be \frac{\exp\big[{\b'\mu^{1-\b}h^{\b}}\big]}{h^{p\b/\b'}}\le Ce^{L_\lambda\b'}\mu^{-p\b/\b'}\sum_{k=0}^\infty |a_k|^{p} e^{\b'k}.\label{le1b01}\ee
We will apply this version of Lemma 3 to prove Proposition 2, by  finding  a number $h_1$, a sequence $\{a_k\}$, both depending on  $f$ and $x$, and a sequence $\{\lambda_k\}$ also depending on $f$ and $x$ but satisfying  $L_\lambda\le C$ independent of $f,x$, such that
  \be {\b'}^{-\frac 1{\b}}A^{-\frac1 {\b}}\big(\Tfc(\tau,x)+B\|f_\tau'\chi_{B_\tau(z)}\|_{\b'}\big)\le h_1\label{le1c}\ee 
\be \sum_{k=0}^{\infty}|a_k|=h_1,\qquad\quad \aa \lambda^{\b'}_{k}|a_k|^{\b'}\le 1-\e,\label{le1d01} \ee 
and \be \sum_{k=0}^{\infty}|a_k|^{p}e^{\b'k}\le \frac {C}{\tau (1-\e)^{p\b/\b'}}\|(Tf)\chi^{}_{E_\tau^c}\|_{p}^{p}.\label{le1d02} \ee 
%Estimate  \eqref{l1a02} in Lemma \ref{le102} will then  follow from \eqref{le1a01}-\eqref{le1d02}, with $\mu=1-\e_\tau\geq 3/4>0$, and $h=h_1$. 
%Similarly \eqref{l1a03} in Lemma \ref{le102} follows with $A_0=A_1,\  \mu=1-\e_\tau\geq 1-\frac{A_1}{C_0}$.\par

\bigskip
Fix $x\in M, 0<\tau< t_0$. Recall that $t_0=\mu\{x\ :\ |Tf|\geq 1\}\le\mu(M)$.  Fix $R>0$ such that 
\be \supp f\subseteq B(x,R).  \nt\ee

% $r_0'=r_{y}(\tau)$. 
Let 
\be N=\bc0&\text{if}\ V_x(R)\le\tau \cr 1&\text{if}\ \tau<V_x(R)\le \tau e^{2\b'}\cr\Big\lceil \dfrac{1}{\b'}\log\dfrac{V_{x}(R)}{\tau}\Big\rceil-1\ \ \ &\text{if}\ V_x(R)>\tau e^{2\b'}\cr\ec\label{nn}\ee
 \be R_j=r_{x}(\tau e^{\b'j}),\qquad j=0,1,2,....\label{Rj}\ee 
 so that 
 \be e^{-\b'}V_x(R) \le V_x(R_N)=\tau e^{\b'N}< V_x(R)\le \tau e^{\b'(N+1)}.\label{RN}\ee
 %and if $N\neq 0$ let also $r_N=R$. If $N\geq 2$, for $j=1,...,N-1$ define
Define

 \be \bc r_0=R_0=r_x(\tau)\cr r_N=R & \text{if} \ N\ge1\ec \ee
 and if $N\ge2$ let for $j\le N-1$
 \be r_{j}= \bc\sup\bigg\{ r\le R_j\ :\ \displaystyle\int_{r_{j-1}\le d(x,y)\le r} |k(x,y)|^\b d\mu(y)\le  \b'A_0\bigg\},\ \  &  \text{if}\ r_{j-1}< r_x(1)\cr
     \sup\bigg\{ r\le R_j\ :\ \displaystyle\int_{r_{j-1}\le d(x,y)\le r} |k(x,y)|^\b d\mu(y)\le  \b'A_\infty\bigg\},\ \  &  \text{if}\ r_{j-1}\ge r_x(1)            \ec
                         \label{defr}\ee
and let 
\be j_1=\bc\min\big\{j\ :\ V_{x}(r_j)\geq 1\big\}\  & \text{if}\ \ \exists j\le N: V_x(r_j)\ge 1\cr N  & \text{if} \ \ \forall j\le N, V_x(r_j)<1.\ec \label{j1}\ee
 Note that with this notation we have, for $N\ge1$, 
 \be r_j\le R_j,\qquad j=0,1,...,N-1,\qquad r_N=R>R_N\label{rN}\ee
 where the last inequality is from \eqref{RN} and the fact that $V_x(r)$ is increasing.
 
 \smallskip
\begin{rk}\label{rk4} From the continuity  in $r$  of the integral inside \eqref{defr} we have that the $\sup$ is actually a $\max$, and when $r_j<R_j$ there is equality
 in the integral condition inside \eqref{defr}.\end{rk}
 
 \smallskip   

 It turns out that $B(x,r_j)$ and $B(x,R_j)$ have comparable volumes:
 If $N\ge 2$ then  we have
 \be e^{-\frac B{A_0}-\frac B{A_\infty}}\,V(R_j)\le V_x(r_j)\le e^{\b'} V_x(R_j),\qquad 0\le j\le N.\label {vols}\ee
\be Q_1\;V(R_j)\le V_x(r_{j+1})- V_x(r_j)\le e^{2\b'} V_x(R_j),\qquad 0\le j\le N-1,\label{vols1}\ee
for some constant $Q_1$ depending only on $\b',B,A_\infty,A_0.$ 
 The proof of the above estimates will be given in the Appendix.
 
  Later we will need the following estimate, valid for all $x\in M:$
 \be \int_{r_{j-1}\le d(x,y)<r_j} |k(x,y)|^\b d\mu(y)\le Q_2,\qquad j=1,...,N,\label{intk}\ee
which follows easily from (K1), (K2), \eqref{RN}, and \eqref{vols}.

 By \eqref{Rj}, \eqref{vols}, \eqref{vols1}, there exists an integer $m$ depending only on $k,\b$ such that if $N\ge m+1$ then
 
 \be \bc V_x(r_j)-V_x(r_0)\geq 2\tau\cr V_x(r_{j+1})-V_x(r_j)\geq 2\tau\ec \qquad m\le j\le N-1.\label{vols2}\ee
 
 On the other hand, if $N\le m$, by \eqref{intk} we  get
 \be\ba \Tfc(\tau,x)&=|T(f_\tau'\chi^{}_{B^c(x,r_x(\tau))})(x)|=\bigg| \int_{r_x(\tau)\le d(x,y)<R}k(x,y)f_\tau'(y)d\mu(y)\bigg|\cr&\le\bigg(\int_{r_0\le d(x,y)\le r_N} |k(x,y)|^\b d\mu(y)\bigg)^{1/\b}\le (N Q_2)^{1/\b}\le (m Q_2)^{1/\b}\le Q_2,\ea \label{mtest2}\ee
where we assumed WLOG that $Q_2>m^{\frac{1}{\b-1}}.$

From now on we assume $\Tfc(\tau,x)> Q_2$
 which then implies that $N\ge m+1\ge2$ and that \eqref{vols2} holds. 
% hence it is enough to assume 
% \be N\geq m+1.\nt\ee

\medskip
For $x\in M$  we will let
\be \bc D_j=B(x,r_j)\ \qquad  \text{for }\; j=0,1,...,N\cr D_0'=B(z,r_z(\tau))\ec\label{Dj}\ee
where the $r_j$ are defined as above (depending on $x, f,\tau$) and where $z$ is defined as in \eqref{e1}. 
 
 \smallskip
\emph{In what follows we will assume, for notational  simplicity, that $m=1$, hence $N\ge 2$. For general $m$, there are a few modifications to make. One is to replace, in the proof below, $D_j$  by $D_{j+m-1}$ for $1\le j\le N-m+1$, and  $D_{1}\setminus D_0$ by $D_m\setminus D_0$. The range of $j,\ J$ should also be changed to $0,...,N-m+1$ accordingly.  Note that by definition $\mu(D_{j+m}\setminus D_{j+m-1})$ is comparable to $V_x(r_{j+1})-V_x(r_j)$ for $1\le j\le N-m$, and $\mu(D_m\setminus D_0)$ is comparable to $V_x(r_{1})-V_x(r_0)$. This fact is used later in \eqref{cl6}, \eqref{lem2} and \eqref{lem6}. }
\emph{The other modification is to define $j_1$ as follows
\be j_1=\bc 0\ \ \ &\text{if} \ V_x(r_0)\ge 1\cr \min\{j\ :\ V_x(r_{j+m-1})\geq 1\}\ \ \ &\text{if}\ V_x(r_0)<1\ \text{and} \ \exists j\le N-m+1: V_x(r_{j+m-1})\geq 1\cr
N-m+1 &\text{if}\ \forall j\le N-m+1 : V_x(r_{j+m-1})< 1\ec\label{j11}\ee
It is clear that this definition coincides with the previous definition of $j_1$ for $m=1$ in \eqref{j1}. Note that if $j_1\neq N-m+1$, we have $\tau e^{\b'j_1}=e^{-\b'(m-1)}\tau e^{\b'(j_1+m-1)}\geq e^{-\b'(m-1)}$, which is needed to get \eqref{le1fc}  and \eqref{rufnew6}. }

 \medskip
 
Let us   split  $f'_\tau=f{\chi_{F_\tau^c}}^{}$ as follows
 \be f_\tau'=\sum_{j=0}^{N-1}f_\tau'\chi_{(D_{j+1}\setminus D_j)\setminus {D_0'}}^{}+f_\tau'\chi_{ D_0\setminus {D_0'}}^{}+f_\tau'\chi_{ {D_0'}}^{},\label{fde1}\ee
and  define 
%    \be r_j= r(\tau) e^{\frac{q}{n}j},\qquad D_j=B(x,\ r_j)  \nt\ee 
  \be \a_j=\|f_\tau'\chi_{(D_{j+1}\setminus D_{j})\setminus D_0'}^{}\|_{\b'},\qquad \a_{-1}=\|f_\tau'\chi^{}_{D_0\setminus D_0'}\|_{\b'},\qquad  \a'=\|f_\tau'\chi^{}_{ D_0'}\|_{\b'},\nt\ee
  
   \be \ab_j=\max{\{\a_{-1},\a_0,...,\a_j,\a'\}},\quad  \b_j=\|f_\tau'\chi_{D^c_{j}\setminus D_0'}^{}\|_{\b'},\quad \b_{-1}=\|f_\tau'\|_{\b'}. \nt \ee 

Notice that for any $j$
\be\a_j\le\b_j\le1,\qquad \b_{j-1}\le\b_j+\a_{j-1}.\label{abj}\ee

 Clearly $\b_j$ is decreasing, and it vanishes when $j=N$, since 
   \be \supp f\subseteq D_N=B(x,r_N).  \nt\ee

\begin{prop}\label{cl}
There exist constants $C_1,C_2>0$ and $C^*>Q_2$ independent of $f$   such that, for any $f$ supported in a ball with $\|f\|_{\b'} \le 1$ and such that $\Tfc(\tau,x)> C^*$, there is an integer  $J\le N-1$ (depending on $f,x$)
with 
\be A_0^{-\frac1 {\b}} (\b')^{-\frac 1{\b}}\big(\Tfc(\tau,x)+B\|f_\tau'\chi_{B_\tau(z)}\|_{\b'}\big)\le\bc  \displaystyle{\!\sum_{j=0}^{J}}\a_j+C_1\ab_J+C_1\b_J, &\hskip -2em \text{if}\ J\le j_1-1,\cr
\displaystyle{\!\sum_{j=0}^{j_1-1}\a_j+\dfrac{A_\infty^{\frac1 {\b}}}{A_0^{\frac1 {\b}}}\!\!\sum_{j=j_1}^{J}}\a_j+C_1\ab_J+C_1\b_J,  &\hskip -.3em \text{if}\ J\ge j_1.\ec
 \label{cl1}\ee
and
\be \sum_{j=0}^{J} \a_j^{p}e^{\b'j}+C_1\ab_J^{p}+C_1\b_J^{p}\leq \frac{C_2}{\tau}\|(Tf)\chi^{}_{E_\tau^c}\|_{p}^{p}.\label{cl2}\ee
\end{prop}

\bigskip

The proof of Proposition \ref{cl} will be given later. Assuming the proposition, our goal is to find a number $h_1$ and a sequence $a=\{a_k\}$ that satisfies \eqref{le1c},\eqref{le1d01} and \eqref{le1d02}.

\bigskip
\ni \underline{Proof of \eqref{le1aaa} in the case $\kappa=0$}.

\smallskip
\ni Assume $\|f\|_{\b'}\le1$, $\|f_\tau'\|_{\b'}^{\b'}\le 1-\e$, and $(Tf)^\circ(\tau,x)>C^*$, where $C^*$ is as in Proposition ~\ref{cl}. Let 
 \be h_1=\sum_{j=0}^J \a_j+C_1\,\ab _J+C_1\b_J.\label{le1e}\ee  If  $A=\max\{A_0, A_\infty\}$, then by Proposition \ref{cl}, we have \be A^{-\frac1 {\b}}(\b')^{-\frac 1{\b}} \big(\Tfc(\tau,x)+B\|f_\tau'\chi_{B_\tau(z)}\|_{\b'}\big)\le h_1.\label{le1f}\ee
Let 
\be a_k=\bc C_1\ab_J &\text{if}\ k=0\cr C_1\b_J &\text{if}\ k=1\cr \a_{k-2} &\text{if}\ k=2,...,J+2.\ec\label{seq1}\ee
Then it is clear that we have 
\be \sum_{k=0}^{J+2}|a_k|=\sum_{k=0}^{J+2}a_k=h_1,\nt\ee
and 
\be\sum_{k=0}^{J+2}|a_k|^{p}e^{\b'k}=C\ab_J^{p}+Ce^{\b'}\b_J^{p}+e^{2\b'}\sum_{j=0}^{J}\a_j^{p}e^{\b'j}\le \frac{C}{\tau}\|(Tf)\chi^{}_{E_\tau^c}\|_{p}^{p},\label{growth1}\ee
where in the last inequality we used  \eqref{cl2} in Proposition \ref{cl}.\par
Next, we show that the second inequality in \eqref{le1d01} holds for such sequence $\{a_k\}$. If there exists $j^*\in \{0,...,J-1\}$ such that $\ab_J=\a_{j^*}$, we take 
%\be \lambda_0=\la_1=\la_{j^*+2}=\la_{J+2}=(1+C_2^{\b'})^{-1/\b'}\nt\e1
\be \la_k=\bc(1+C_1^{\b'})^{-1/\b'} \ \ &\text{if}\ k=0,1,j^*+2,J+2\cr 1&\text{otherwise}\ec\nt\ee
so that $L_\lambda=\sum_0^\infty(\lceil\lambda_k^{-\b}\rceil-1)\le C$ independent of $f,x,z$ and 
\be\ba \aa \lambda^{\b'}_{k}|a_k|^{\b'}&=\sum_{j=0,j\neq j^*,J}^{J-1}\a_j^{\b'}+\frac{C_1^{\b'}+1}{(1+C_1^{\b'})}\a^{\b'}_{j^*}+\frac{1}{1+C_1^{\b'}}\a_J^{\b'}+\frac{C_1^{\b'}}{1+C_1^{\b'}}\b_J^{\b'}\cr
&\le \sum_{j=0}^{J-1}\a_j^{\b'} + \b_J^{\b'}=\|f_\tau'\chi^{}_{D_0^c\setminus D_0'}\|_{\b'}^{\b'}\le 1-\e.\label{lamk}\ea\ee
\par

In the case where $\ab_J=\a_{J}$ and $\a_{-1}$, it is obvious that by taking slightly different $\lambda_k$ we still get \eqref{le1d01}, so we omit it here. Therefore we find a sequence that satisfies \eqref{le1c} to \eqref{le1d02} and hence \eqref{l1a02} in Proposition \ref{le102} follows.

\bigskip\ni
\underline{Proof of \eqref{le1aaa} in the case $\kappa>0$.} 

\smallskip

\ni Recall that in this case $A=A_0$. Also, it is enough to assume $A_\infty>A_0$, otherwise we can apply the previous case, since \eqref{1b1} implies \eqref{1b}.  Let \be h_1=\bc \displaystyle{\sum_{j=0}^{J}\a_j}+C_1\ab_J+C_1\b_J, &\text{if}\ J\le j_1-1,\cr
\displaystyle{\sum_{j=0}^{j_1-1}\a_j}+\dfrac{A_\infty^{\frac1 {\b}}}{A_0^{\frac1 {\b}}}\displaystyle{\sum_{j=j_1}^{J}\a_j}+C_1\ab_J+C_1\b_J,  &\text{if}\ J\ge j_1.\ec
 \label{h1}\ee
 Then clearly \be A_0^{-\frac1 {\b}}\b'^{-\frac 1{\b}} \big(\Tfc(\tau,x)+B\|f_\tau'\chi_{B_\tau(z)}\|_{\b'}\big)\le h_1.\label{h11}\ee
 If $J\le j_1-1$, then the sequence defined in \eqref{seq1} will work by the same argument as in the proof of \eqref{le1aaa} under condition \eqref{1b}. Now suppose that $J\ge j_1$ and  define
\be a_0=C_1\ab_J,\qquad a_1=C_1\b_J,\qquad\text{and}\qquad a_k=\bc \a_{k-2}\ &\text{for}\ k=2,...,j_1+1,\cr \dfrac{A_\infty^{\frac1 {\b}}}{A_0^{\frac1 {\b}}}\a_{k-2} &\text{for}\ k=j_1+2,...,J+2.\ec\label{seq2}\ee 
Then  the $\{a_k\}$ satisfies the first identity in \eqref{le1d01}, and the estimate   \eqref{le1d02}, using the same argument as the one under condition \eqref{1b}. 

To deal with $\|\lambda a\|_{\b'}$, fix any $n_0\in{\mathbb N}$ and suppose that 
\be 0\le J-j_1\le n_0.\label{seq4}\ee
In view of \eqref{seq2}, that we can find suitable $\{\la_k\}$, where at most $n_0+4$ of the $\la_k$ are less than $1$, and the value of each $\la_k$ only depends on $C_1,A_\infty,A_0,\b$, so that $L_\lambda\le C n_0$. Arguing as  in the previous case we see that the second inequality in \eqref{le1d01} holds. For example, in the case that $\ab_J=\a_{j^*}$ with $j^*\in\{0,...,j_1-1\}$, we can take
\be \la_k=\bc(1+C_1^{\b'})^{-1/\b'} \ \ &\text{if}\ k=0,j^*+2\cr (\frac{A_\infty}{A_0}+C_1^{\b'})^{-1/\b'} \ \ &\text{if}\ k=1,j_1+2,...,J+2\cr1&\text{otherwise}\ec\label{lambda1}\ee
Therefore, arguing as in \eqref{lamk}, we have the following
\be\ba\aa &\lambda^{\b'}_{k}|a_k|^{\b'}\le\|f_\tau'\|_{\b'}^{\b'}+\k\|(Tf)\chi^{}_{E_\tau^c}\|_{\b'}^{\b'} \le 1-\e.
\ea\nt\ee
 On the other hand, assume that \eqref{seq4} does not hold, that is 
\be J-j_1>n_0.\label{seq5}\ee

Note that by \eqref{cl2} with $p=\b'$, we have that for any $j\le J$,
 \be \tau e^{\b'j}\sum_{k=j}^{J}\a_k^{\b'}\le C_2\|(Tf)\chi^{}_{E_\tau^c}\|^{\b'}_{\b'}.\label{le1fb}\ee
 Since $\tau e^{\b' j_1}=V_x(R_{j_1})\geq V_x(r_{j_1}^{})\ge 1$,  from \eqref{seq2}, \eqref{seq5}, and  \eqref{le1fb} (with $j=j_1+n_0$)  we have
 \be \sum_{k=j_1+n_0+2}^{J+2}|a_k|^{\b'}=\frac{A_\infty}{A_0}\sum_{j=j_1+n_0}^{J}\a_j^{\b'}\le\frac{C_2 A_\infty}{ A_0 e^{\b'n_0}}\|(Tf)\chi^{}_{E_\tau^c}\|^{\b'}_{\b'}\le \k\|(Tf)\chi^{}_{E_\tau^c}\|^{\b'}_{\b'},\label{le1fc}\ee
 provided $n_0$ is chosen sufficiently large, independently of $x,z,f.$
 Once again,  we can choose suitable ${\la_k}$ similar to \eqref{lambda1}, with at most $n_0+2$ of the $\la_k$ less than $1$. Indeed, we can write 
 \be h_1=\sum_{k=0}^{J+2}a_k= C_1\ab_J+C_1\b_J+\sum_{j=0}^{j_1-1}\a_j+\dfrac{A_\infty^{\frac1 {\b}}}{A_0^{\frac1 {\b}}}\sum_{j=j_1}^{j_1+n_0-1}\a_j+\dfrac{A_\infty^{\frac1 {\b}}}{A_0^{\frac1 {\b}}}\sum_{j=j_1+n_0}^J \a_j\ee  so that we  choose $\lambda_k=1$ for $2\le k\le j_1+1$ and $j_1+n_0+2\le k\le J+2$, and the other $n_0+2$ values of $\lambda_k$ chosen as  in \eqref{lambda1}. 
 With this choice we have $L_\lambda\le Cn_0\le C$ and  
 \be\ba\aa \lambda^{\b'}_{k}|a_k|^{\b'}=\sum_{k=0}^{j_1+n_0+1}+\sum_{k=j_1+n_0+2}^{J+2}\le\|f_\tau'\|_{\b'}^{\b'}+\k\|(Tf)\chi^{}_{E_\tau^c}\|_{\b'}^{\b'}
 \le 1-\e.\ea\label{growth4}\ee
Therefore \eqref{le1d01} holds and \eqref{l1a02} in Proposition  \ref{le102} follows. 

 \hfill\QEDopen
 
\bigskip
\ni{\bf Remark 3.}  In the above proof the reader  can appreciate the reason why under $\|f\|_{\b'}\le~1$ one cannot use the constant $A_0$ in the main inequality of Proposition 2, and hence of Theorem~1. Indeed, if $A_\infty>A_0$ one can choose  $h_1$ as  \eqref{h1} so that \eqref{h11} holds. The  corresponding sequence $\lambda=\{\lambda_k\}$ defined in \eqref{lambda1} will statisfy  $\|\lambda a\|_{\b'}^{\b'}\le \|f\|_{\b'}^{\b'}(1-\e)$,  however the number of $\lambda_k$  strictly less than 1 can grow arbitrarily, as the support of $f$ gets larger,  so that  the inequality $L_\lambda\le C$ can fail. 
\hfill\QEDopen

\bigskip
\ni
{\underline{\bf {Proof of \eqref{rufnew2}}.}}

\medskip 
In this case we are assuming that \eqref{1b1} and \eqref{lambda11} hold  with $\kappa>0$. 
It is enough to show that  for all $x\in E_\tau$ with $(Tf)^\circ(\tau,x)> C^*$ (with $C^*$ as in Prop. \ref{cl}) we have 
\be {\exp\bigg[\dfrac{(1-\e)^{1-\b}}{A_0}\big(\Tfc(\tau,x)+B\|f_\tau'\chi_{B_\tau(z)}\|_{\b'}\big)^{\b}\bigg]}\leq \frac{C}{\tau}.\label{rufnew3}\ee
We first consider the case where $\ab_J=\a'$ and $\ J-j_1\geq 2$. By \eqref{lambda11} and \eqref{cl2} we have
\be \ba &\sum_{j=0}^{J-1}\a_{j}^{\b'}+(\a')^{\b'}+\b_J^{\b'}+C_2^{-1}\k\tau\sum_{j=0}^{J-1}\a^{\b'}_j e^{\b'j}\cr
&\le \|f_\tau'\|^{\b'}_{\b'}+\k\|(Tf)\chi^{}_{E^c_\tau}\|\le 1-\e
\ea\label{rufnew4}\ee
By \eqref{cl1} in proposition \ref{cl}, we  get
\be \ba &A_0^{-1/\b}\big(\Tfc(\tau,x)+B\|f_\tau'\chi_{B_\tau(z)}\|_{\b'}\big)\le (\b')^{1/\b}\bigg(\sum_{j=0}^{j_1-1}\a_j+\big(\frac{A_\infty}{A_0}\big)^{1/\b}\sum_{j=j_1}^{J-1}\a_j+C_1\a'+(1+C_1)\b_J\bigg)
\cr&\le (\b')^{1/\b}\bigg(\sum_{j=0}^{J-1}(1+C_2^{-1}\k\tau e^{\b'j})\a^{\b'}_j+(\a')^{\b'}+\b_J^{\b'}\bigg)^{1/\b'}\times
\cr&\times \bigg(\sum_{j=0}^{j_1-1}(1+C_2^{-1}\k\tau e^{\b'j})^{-\b/\b'}+\frac{A_\infty}{A_0}\sum_{j=j_1}^{J-1}(1+C_2^{-1}\k\tau e^{\b'j})^{-\b/\b'} +(2^{\b})(C_1^{\b}+1)\bigg)^{1/\b}.
\ea\label{rufnew5}\ee
By definition of $j_1$ in \eqref{j1} and \eqref{vols}, \eqref{vols1} we have that 
\be \tau e^{\b'j_1}\geq 1,\qquad j_1\le \frac{1}{\b'}\log\frac{1}{\tau}+C.\nt\ee
Hence by \eqref{rufnew4} we get
\be\ba A^{-1}_0\big(\Tfc(\tau,x)+B\|f_\tau'\chi_{B_\tau(z)}\|_{\b'}\big)^\b&\le \b'(1-\e)^{\b-1}(j_1+C\sum_{j=j_1}^{J-1}e^{-\b(j-j_1)}+C)\cr& \le (1-\e)^{\b-1}(\log\frac{1}{\tau}+C).\label{rufnew6}\ea\ee

For the other cases such as $\ab_J=\a_j$ for some $j\le J$, or $J\le j_1-1$, the proof is completely similar, so we omit it here. 
\hfill\QEDopen

%###########################################################
%\input{5-prop3-pf}

 \section{Proof of Proposition \ref{cl}}\label{prop1}
Recall that $R_j=r_x(\tau e^{\b'j})$, $V_x(R_j)=\mu(B(x,R_j))=\tau e^{\b'j}$, and  $r_j$ is defined in \eqref{defr}. The integer $N$ was defined in \eqref{nn}, and we are assuming $\Tfc(\tau,x)>Q_2$, where $Q_2$ is as in \eqref{intk}, which implies  $N\ge 2$ and \eqref{vols2}.

 \bigskip
For the rest of this section we will  set for any measurable function $\phi:M\to \R^m$
   \be S_j \phi=\phi\chi_{D_j^c\setminus D_0'}^{}, \nt\ee 
   where we defined $D_j,\ D_0'$ in \eqref{Dj}.
 With this notation we then have
   \be( S_j-S_{j+1})f_\tau'=f_\tau'\chi_{(D_{j+1}\setminus D_j)\setminus D_0'}^{} ,\qquad \a_j=\|( S_j-S_{j+1})f_\tau'\|_{\b'}, \quad \b_j=\|S_j f_\tau'\|_{\b'},\;\;\; j\ge0. \nt\ee 
\def\supp{{\rm supp }}\def\N{{\mathbb N}}

%and from \eqref{fde1} 

%\be \ba f_\tau'&=S_0f_\tau'+f_\tau' \chi_{D_0\setminus D_0'}^{}+f_\tau' \chi_{D_0'}^{}.\cr\ea\label {de22}\ee

  We first give some preliminary estimates on $(Tf)^\circ(\tau,x).$ 
  Recall that  \be \Tfc(\tau,x)=|T(f_\tau'\chi^{}_{D_0^c})(x)|\le|Tf_\tau'\chi^{}_{D_0^c\setminus D_0'}(x)|+|Tf_\tau'\chi^{}_{ D_0^c\cap D_0'}(x)|. \label{mtau1}\ee
   By H\"older and  \eqref{k3} we have
  \be\ba |Tf_\tau'\chi^{}_{ D_0^c\cap D_0'}(x)|&\le \Big(\int_{ \{d(x,y)\ge r_0\}\cap D_0'}|k(x,y)|^\b d\mu(y)\Big)^{1/\b}\Big(\int_{D_0'}|f'_\tau(y)|^{\b'}d\mu(y)\Big)^{1/\b'}\cr&\le B^{1/\b}V_x(r_0)^{-1/\b}\mu(D_0')^{1/\b}\|f_\tau'\chi_{D_0'}\|_{\b'}= B^{1/\b}\a'\le B\a'.\ea\label{fla8}\ee
  
 By \eqref{mtau1} we can write, for any $J\in\{0,1,...,N-1\}$
  \be \ba \Tfc(\tau,x)+B\|f_\tau'\chi_{B_\tau(z)}\|_{\b'}&\le |{T}S_0f_\tau'(x)|+2B\a'=\bigg|  \sum_{j=0}^{J}T(S_j-S_{j+1})f_\tau'+TS_{J+1}f_\tau'\bigg|+2B\a'\cr&\le\sum_{j=0}^{J}|{T}\big(S_jf_\tau'-S_{j+1}f_\tau'\big)(x)|+|{T}S_{J+1}f_\tau'(x)|+2B\a'.
  \ea
  \label{fla7}\ee
  
For any integer $j$, by \eqref{defr}, \eqref{intk} we have the estimate
   \be\ba |T(S_jf_\tau'-S_{j+1}f_\tau')(x)|&\le \bigg(\mathop\int\limits_{D_{j+1}\setminus D_j}|k(x,y)|^{\b}d\mu(y)\bigg)^{1/\b}\|S_jf_\tau'-S_{j+1}f_\tau'\|_{\b'}\cr
   &\le\bc (\b' A_0)^{1/\b}\a_j,\ \ \ &\text{if}\ 1\le j\le j_1-1,\; j\neq N-1\cr(\b' A_\infty)^{1/\b}\a_j,\ \ \ &\text{if}\ j_1\leq j\le N-1,\;  j\neq 0\cr (Q_2)^{1/\b} \a_j &\text{if}\ j=0,\; j=N-1.\ec\ea\label{la3} \ee 

   We then get that there is $Q_3$ such that  for any $J=0,1,...,N-1$
    \be \ba &\Tfc(\tau,x)+B\|f_\tau'\chi_{B_\tau(z)}\|_{\b'}\cr&\leq\bc (\b' A_0)^{1/\b}\sum_ {j=0}^{J}\a_j+ Q_3\ab_J + |TS_{J+1}f_\tau'(x)|,&\text{if}\ J\le j_1-1,\cr(\b' A_0)^{1/\b}\sum_ {j=0}^{j_1-1}\a_j+(\b' A_\infty)^{1/\b}\sum_ {j=j_1}^{J}\a_j+ Q_3\ab_J + |TS_{J+1}f_\tau'(x)|,\ &\text{if}\ J\ge j_1 .\ec
  \ea
  \label{la1}\ee

  \par
 \ni  What remains to be proved is that for $(Tf)^\circ(\tau,x)>C^*$ large enough,  there is   $J\in\{0,1,...,N-~1\}$ (depending on $x$) and $Q^*$ (depending on $k,\b$) such that 
\be |TS_{J+1}f_\tau'(x)|\le  Q^*\ab_J +Q^*\b_J \label{clcl2},\ee
(so that \eqref{cl1} holds for some $C_1$) and \eqref{cl2} holds i.e.
\be \sum_{j=0}^{J} \a_j^{p}e^{\b'j}+C_1\ab_J^{p} +C_1\b_J^{p}\leq \frac{C_2}{\tau}\|(Tf)\chi^{}_{E_\tau^c}\|_{p}^{p},\label{cl2b}\ee
for some $C_2>0$.
%  Recall that $N$ is an integer such that $\supp f\subseteq D_N$.  

 For each $j=0,1,...,N-1$ consider the estimates

\be \b_{j-1}\le \avgdj |Tf|d\mu.  \label{j1a}\ee
\be|TS_jf_\tau'(x)|\le \frac{2e^{\b'-1}}{e^{\b'-1}+1}\, |TS_{j+1}f_\tau'(x)|.\label {j2a}\ee
If \eqref{j1a} holds for  $j=0$ let 
\be J_1=\max\big\{k\in\{1,...,N\}: \; {\text{\eqref{j1a} holds for }} 0\le j\le k-1\,\big\}   \ee
whereas if \eqref{j1a}  does not hold    for $j=0$ let $J_1=0$.

\smallskip
\ni Similarly, if \eqref{j2a} holds for  $j=0$ let 
\be J_2=\max\big\{k\in\{1,...,N\}: \; {\text{\eqref{j2a} holds for }} 0\le j\le k-1\,\big\}   \ee
whereas if \eqref{j2a}  does not hold    for $j=0$ let $J_2=0$.

 The integers $J_1, J_2$ could be considered as  stopping times: they represent the first times after 0 that \eqref{j1a}, \eqref{j2a}  do not occur, if  they do occur at $0$.

  \bigskip
 \ni {\underline{\emph{Proof of \eqref{clcl2}, \eqref{cl2b} in  case $0\le J_2\leq J_1\leq N$.}}}
 
  \bigskip
First assume that $J_2\neq N$. From the definition of $J_2$ we have 
 \be \ba |{T}S_{J_2+1}f_\tau'(x)|&< \frac{e^{{\b'}-1}+1}{2e^{{\b'}-1}} \; |{T}S_{J_2}f_\tau'(x)| \cr
% &= \left(\frac{e^{{\b'}-1}+1}{2e^{{\b'}-1}}  \right) \big({T}S_{J_2+1}f_\tau'(x)+{T}(S_{J_2}-S_{J_2+1})f_\tau'(x)\big)\cr
 &\leq \frac{e^{{\b'}-1}+1}{2e^{{\b'}-1}}  \; \big(|{T}S_{J_2+1}f_\tau'(x)|+\big|{T}(S_{J_2}-S_{J_2+1})f_\tau'(x)\big|\big)\cr
 &\le \frac{e^{{\b'}-1}+1}{2e^{{\b'}-1}}\;   \big(|{T}S_{J_2+1}f_\tau'(x)|+Q_4\a_{J_2}\big)\cr
 \ea\label{cl3}\ee
 where the last inequality is by \eqref{la3}. This implies 
 \be |{T}S_{J_2+1}f_\tau'(x)|<Q_5\a_{J_2}\le Q_5 \ab_{J_2},\label{cl4}\ee
 which is   \eqref{clcl2} with  $J=J_2$.  
 
 To show \eqref{cl2b}, note that if $J_1=0$, then $J_2=0$ and from \eqref{la1} there is $Q_6$ such that 
 \be \Tfc(\tau,x)\le Q_6 \ab_{0}\le Q_6.\ee
 Therefore, assuming $\Tfc(\tau,x)> \max\{Q_6,Q_2\}$ 
we have $J_1\neq 0$ and \eqref{j1a} is true for all $j\leq J-1$, so
 \be  \ba &\sum_{j=0}^{J} \a_j^{p}e^{{\b'}j}+C_1\ab_J^{p} +C_1\b_J^{p}\leq Q_7 \sum_{j=0}^{J-1} \b_{j-1}^{p}e^{{\b'}j}\le Q_7\sum_{j=0}^{J-1} \left(\avgdj |Tf|d\mu\right)^{p}e^{{\b'}j}\cr&  \le Q_7\sum_{j=0}^{J-1} e^{{\b'}j} \avgdj |Tf|^{p}d\mu
 \le\frac {2Q_7}\tau \sum_{j=0}^{J-1} \frac{V_x(R_j)}{V_x(r_{j+1})-V_x(r_j)}\dj |Tf|^{p}d\mu\cr&\le \frac{C_2}{\tau}\sum_{j=0}^{J-1}\dj |Tf|^{p}d\mu\le \frac{C_2}{\tau}||(Tf)\chi^{}_{E_\tau^c}||_{p}^{p}
 \ea\label{cl6}\ee
where the first inequality is due to the fact that $\b_j$ is decreasing,  the third last inequality is by that $\mu(E_\tau)= \tau$ and \eqref{vols2}, and whereas the second last inequality is by \eqref{vols1}.  

\smallskip
If $J_2 = N$, then $J_1 = N$, and if we take $J=N-1$ it is clear that \eqref{clcl2} holds, since $TS_Nf_\tau'(x) = 0$.
Note that \eqref{j2a} is true for $j = 0,....,N-1$, therefore \eqref{cl2b} still follows by the same
calculations as in \eqref{cl6}. \\
% If $J_1=0$, then we just need to check \eqref{cl2} for $J=0$:
% \be  \a_0^{\si\b'}+\b_{0}^{\si\b'}\leq 2\b_0^{\si\b'}\le C=C\frac{\tau}{\tau}\le \frac{C}{\tau}||Tf||_{\si\b'}^{\si\b'}. \label{cl666}\ee
 %Case 2
 
 \bigskip
  \ni {\underline{\emph{Proof of \eqref{clcl2}, \eqref{cl2b}  in case $N\ge J_2\geq J_1+1$.}}}
  
   \bigskip If $J_1= N-1$ then the proof is the same as in the case $J_1 = J_2 = N$ given above, so
 we can assume $J_1 <  N - 1$.

 We will need the following lemma to handle this case. Let us state it here, and its proof will be postponed to the Appendix.
 \begin{lemma}\label{lb}
  There exist  constants $C_3,C>0$, depending only on $k,\b$,  such that for any $j\le N-2$ 
  \be \avgdj |Tf_\tau|d\mu\leq C_3 e^{-(\b'-1)j},\label{lb1}\ee
  
\be \bigg| {T}S_{j+1}f_\tau'(x)- \avgdj {T}S_{j+2}f_\tau'd\mu\bigg|\le C\b_{j+1},\label{lc1}\ee

 \be \bigg|\avgdj {T}(S_0-S_{j+2})f_\tau'd\mu\bigg|
\le C\ab_{j+1},\label{ld1}\ee

\be \bigg|\avgdj T(f_\tau'\chi^{}_{D_0\setminus D_0'})d\mu\bigg|
\le C\a_{-1},\label{ld2}\ee
and
\be\bigg|\avgdj T(f_\tau'\chi^{}_{D_0'})d\mu\bigg|
\le C\a' .\label{ldd1}\ee
\end{lemma} 
 %
 %reduction
  %Assuming Lemma \ref{lb}, let us first make a reduction. 
  By \eqref{lc1} in Lemma \ref{lb} and noticing  that 
 \be f=f_\tau+f_\tau'=f_\tau+f_\tau'\chi^{}_{D_0\setminus D_0'}+f_\tau'\chi^{}_{(D_{J_1+2}\setminus D_0)\setminus D_0'}+f_\tau'\chi^{}_{D_{J_1+2}^c\setminus D_0'}+f_\tau'\chi^{}_{D_0'},\nt\ee
  we have, using \eqref{lc1}-\eqref{ldd1},
 \be \ba & |{T}S_{J_1+1}f_\tau'(x)|\leq  \bigg|\avgdJo {T}S_{J_1+2}f_\tau'd\mu\bigg|+C\b_{J_1+1}\cr
 &=\bigg|\avgdJo\big(Tf-Tf_\tau-T(f_\tau'\chi^{}_{D_0\setminus D_0'}) -T(S_0-S_{J_1+2})f_\tau'-T(f_\tau'\chi^{}_{D_0'})\big) d\mu\bigg|+C\b_{J_1+1}\cr
% &\le\bigg|\avgdJo Tfd\mu\bigg|+\bigg|\avgdJo Tf_\tau d\mu\bigg|+\bigg|\avgdJo T(f_\tau '\chi^{}_{D_0\setminus D_0'})d\mu\bigg|\cr&+\bigg|\avgdJo T(S_0-S_{J_1+2})f_\tau'd\mu\bigg|+\bigg|\avgdJo T(f_\tau '\chi^{}_{D_0'})d\mu\bigg|+C\b_{J_1+1}\cr
 &\le \avgdJo |Tf|d\mu+\avgdJo |Tf_\tau|d\mu+C\a_{-1}+C\ab_{J_1+1}+C\b_{J_1+1}.
 \ea\label{cl8}\ee

 To estimate the second integral, note first that by \eqref{mtau1}, \eqref{fla8}
  \be |{T}S_0f_\tau'(x)|\geq  \Tfc(\tau,x)-B. \label{cl7}\ee
  Hence, if 
  \be \Tfc(\tau,x)> \max\{B+4C_3,Q_2\},\label{mtest3}\ee
  where $C_3$ is the constant in \eqref{lb1}, then $C_3\le \frac1 4 |TS_0f_\tau'(x)|$.
 Using \eqref{lb1} in Lemma \ref{lb}, and condition \eqref{j2a} applied $J_1+1$ times (which is possible since $J_1+1\le J_2$), we get
 \be \ba\avgdJo |Tf_\tau|d\mu&\leq C_3e^{-(\b'-1)J_1} \leq \frac{1}{4}e^{-(\b'-1)J_1}|{T}S_0f_\tau'(x)|\cr&\le \frac{1}{4}e^{-(\b'-1)J_1}\left(\frac{2e^{{\b'}-1}}{e^{{\b'}-1}+1}\right)^{J_1+1}|{T}S_{J_1+1}f_\tau'(x)|\cr
 &= \frac{1}{4}\left(\frac{2}{e^{{\b'}-1}+1}\right)^{J_1}\left(\frac{2e^{{\b'}-1}}{e^{{\b'}-1}+1}\right)|{T}S_{J_1+1}f_\tau'(x)|\le \frac{1}{2}|{T}S_{J_1+1}f_\tau'(x)|,
 \ea\label{cl7}\ee
 and
 \be \ba  |{T}S_{J_1+1}f_\tau'(x)|\leq  \avgdJo |Tf|d\mu+\frac{1}{2}|{T}S_{J_1+1}f_\tau'(x)|+C\ab_{J_1+1}+C\b_{J_1}.
 \ea\nt\ee
 So the above inequality, along with the fact that  \eqref{j1a} is false for $j=J_1$, give us 
 \be \ba |{T}S_{J_1+1}f_\tau'(x)|&\leq 2\avgdJo |Tf|d\mu+2C\ab_{J_1+1}+2C\b_{J_1+1}\cr
 &\le 2\b_{J_1-1}+ 2C\ab_{J_1+1}+2C\b_{J_1+1}\le Q_8\ab_{J_1}+Q_8\b_{J_1}.\ea\label{cl9}\ee
which is \eqref{clcl2}, with $J=J_1+1$ (in the last inequality above we used   from \eqref{abj} that for any $j$ we have $\b_{j-1}\le \b_j+\a_{j-1}\le \b_j+\ab_j$ and that $\ab_{j+1}\le \max\{\ab_j,\a_{j+1}\}\le\ab_j+\beta_{j+1}\le \ab_j+\b_j.$)

By the same argument as in the previous  case, if $J_1=0$ then we have, from \eqref{la1} that  $\Tfc(\tau,x)<Q_9$ for some $Q_9>0$, so that if $\Tfc(\tau,x)>\max\{Q_9,Q_2,B+4C_3\}$ then $J_1 \neq 0$, and estimate \eqref{cl2b} for this choice of $J$ follows exactly in the same way as in \eqref{cl6}.

  \bigskip
Finally, we can take $C^*=\max\{Q_2,Q_6,Q_9,B+4C_3\}$ and $Q^*=\max\{Q_5,Q_8\}$, to cover all of the above cases.
 %Case 3

\hfill\QEDopen \\
 
% \bigskip
% \underline{\bf Proof of Proposition \ref{clprop}:} The proof of Proposition \ref{clprop} is completely similar, just replace $A_2$ by $C_0$. So we omit the proof here.
%pf of lemmas 9,10,11

%#########################################################
%\input{6-Ainfty0-pf}

\section{Proof of Theorem \ref{m1} when $A_\infty=0$  }\label{Ainfinity}
%\centerline{\bf \Large{Sharp Adams inequalities with exact growth condition }}
%\centerline{\bf\Large{on metric measure spaces and applications}}
%\font\twelvemi=cmmi12 at 13pt
%\renewcommand{\chi}{\raisebox{.13\baselineskip}{\hbox{\twelvemi\char31}}}
%\vskip 0.1in
%\centerline{\bf{Carlo Morpurgo,\ Liuyu Qin}}

%
%
%
%%
%% \chapter*{Acknowledgements}
%\begin{abstract}
%
%\end{abstract}
%\numberwithin{equation}{section}

Since \eqref{k3} is a more restrictive condition than \eqref{k3} with small $V_x$ and \eqref{kk12}, it is enough to prove the inequalities under the conditions $A_\infty=0$, the first estimate in \eqref{k3} and \eqref{k4} for $V_x(d(x,y))\le 1$, and \eqref{kk12}.   

%Therefore, it also proves the inequalities in the case when $A_\infty=0$ and \eqref{k1}-\eqref{k4} hold.

\smallskip

The proof in this case is similar to the proof under conditions \eqref{k1}-\eqref{k4} and $A_\infty\neq 0$ but with some appropriate modifications. Most parts of the proof follow immediately given the above assumptions. We will mention the main changes as well as the related proof.  

First we prove that \eqref{k1p} still holds, since it is needed to get \eqref{Adams}, the Adams inequality on measure spaces with finite measure.  Let $0\le V_x(r_1)\le V_x(r_2)\le 1$, then by \eqref{k1}, \eqref{k3} for small balls, \eqref{kk12} and \eqref{k2} with $A_\infty=0$, 

\be \ba  &\int_{V_x(r_1)}^{V_x(r_2)}  (k^*(x,u))^\b du\le \mathop\int\limits_{r_1\le d(x,y)\le r_2}|k(x,y)|^\b d\mu(y)+\mathop\int\limits_{d(x,y)\le r_1,\ |k(x,y)|\le k^*(x,V_x(r_1))}|k(x,y)|^\b d\mu(y)\cr
&+V_x(r_2)\supess_{r_2\le d(x,y)\le r_x(1)} |k(x,y)|^\b+\mathop\int\limits_{ d(x,y)\geq r_x(1)} |k(x,y)|^\b d\mu(y)\le \mathop\int\limits_{r_1\le d(x,y)\le r_2}|k(x,y)|^\b d\mu(y)+C.\ea\ee

%It is also clear that (3.17) holds by (K1) and (K2) with $A_\infty=0$. 

Next, the definition of $r_j$ in \eqref{defr} should be modified since $A_\infty=0$. If $V_x(R)\le 1$, let $N$ be the same as in definition \eqref{nn} . Otherwise $N$ is defined by replacing all the $V_x(R)$ by $1$ in \eqref{nn}, that is 

\be N=\bc0&\text{if}\ 1\le\tau \cr 1&\text{if}\ \tau<1\le \tau e^{2\b'}\cr\Big\lceil \dfrac{1}{\b'}\log\dfrac{1}{\tau}\Big\rceil-1\ \ \ &\text{if}\ 1>\tau e^{2\b'}\cr\ec\label{hy5}\ee

 Let 

\be r_j=\bc R_0,\ \ \ & j=0\cr 
\sup\{ r\le R_j\ :\ \mathop\int\limits_{r_{j-1}\le d(x,y)\le r} |k(x,y)|^\b d\mu(y)\le  \b'A_0\}, &j=1,...,N-1\ \text{if}\ N\geq 2\cr
R, &j=N,\ec\label{hy2}\ee
where $R_j$ is defined in \eqref{Rj}, that is, $R_j=r_x(\tau e^{\b'j})$. Note that 
\be V_x(r_j)\le 1,\ \ \text{if}\ \  0\le j\le N-1,\ N\geq 1,\label{hy7}\ee
and $V_x(r_N)>1$ if and only if $\tau>1$ or $V_x(R)>1$.

Since we have \eqref{k1} and \eqref{k3} for small balls, by the same argument in Appendix \ref{Vols}, we get \eqref{vols} for $0\le j\le N-1$, and \eqref{vols1} for $0\le j\le N-2$ with slightly different constants depending only on $\b',\ B,\ A_0$. On the other hand, by the definition of $R_j, r_j$ for $j=N-1,N$, we have
\be V_x(r_N)\geq CV_x(R_N),\qquad  V_x(r_N)- V_x(r_{N-1})\geq CV_x(R_{N-1}).\label{hy6}\ee
Recall that $m$ is defined in \eqref{vols2}, and we will assume $m=1$ as before.
The inequalities that were derived by using \eqref{vols},\eqref{vols1}, \eqref{hy6} are \eqref{cl6},\eqref{lem2} and \eqref{lem6}. It is easy to verify that they still hold.

The inequalities \eqref{circ0},\eqref{fla8}, which were proved using \eqref{k3}, still hold by \eqref{k1} and \eqref{k2} with $A_\infty=0$.

For the rest of the proof in sections \ref{lemma2}, \ref{prop1}, we always take $j_1=N$ by definition and $ A_\infty=0$. 
%Recall the definition of $j_1$
%\be j_1=\bc \min\{j\ :\ V_x(r_j)\geq 1\}\ \ \ &\text{if}\ V_x(r_j)\geq 1\ \text{for some}\ 0\le j\le N\cr
%N &\text{otherwise}.\ec\label{hy8}\ee
%Therefore by \eqref{hy7} we get $j_1=N$.

Lastly, we mention that \eqref{lc4} still holds. We can apply \eqref{k4} for $V_x(r)\le 1$ since by \eqref{hy7} we have $V_x(r_{N-1})\le 1$ and clearly $\xi\in B(x,r_{N-1}) $.

\hfill\qed

%############################################################
%\input{7-thm2-pf}

 \section{Proof of Theorem 2}\label{sharpness} For simplicity, in this proof we will let $V(r)=V_{x_0}(r)$. Let $B(x_0,\e_0)$ be the largest ball with zero volume, and  given any $\e>\e_0$ we will let $\e'=r_{x_0}^{}(\frac 1 2 V(\e))$, i.e. $B(x_0,\e')$ is the smallest ball inside $B(x_0,\e)$  with
 $$V(\e)=2 V(\e').$$
 
  In the notation of \eqref{yreg2}, \eqref{yreg3} let
 \be P(x,y)=\sum_{j=0}^{m-1} k_j(x,x_0)p_j(y,x_0),\qquad Q(y,z)=\sum_{j=0}^{m-1}{v_j(y,x_0)}{v_j(z,x_0)}\label {pq}\ee
 i.e. $Q(y,z)$ is the kernel of the orthogonal projection  of $L^2(B(x_0,r))$ to the space spanned by the $p_j(\cdot ,x_0)$ restricted to the ball $B(x_0,r)$. Using \eqref{yreg3} we get
\be |Q(y,z)|\le \frac{C}{V(r)},\qquad y,z\in B(x_0,r).\label{Q}\ee

It's easy to check that \eqref{xreg}  implies

% (1) Suppose there exists $x_0\in M$ and $P(x,y)\ :\ B^c(x_0,R)\times B(x_0,r)\to \R$ such that for all $0<r<R$ with $V(R)= 2V(r)$
% \be\mathop\int\limits_{r_1\le d(x,y)\le r_2} |k(x_0,y)|^\b d\mu(y)\geq  A_0\log{\frac{V(r_2)}{V(r_1)}}-C,  \quad  0<V(r_1)<V(r_2)\le1,\label{k5}\ee

\be \int_{d(x_0,y)>\e}|(k(x,y)-k(x_0,y))||k(x_0,y)|^{\b-1}d\mu(y)\le C,\ \ \ x\in B(x_0,\e') \label{k6}\ee
%and assume there exists $P(x,y)\ :\ B^c(x_0,R)\times B(x_0,r)\to \R$ such that
(note that the above estimate implies  (K4), if valid for all $x_0$).

%\be C_1|k(x_0,x)|\le |P(x,y)|\le C_2|k(x_0,x)|, \qquad \text{for}\ x\in B^c(x_0,R),\ y\in B(x_0,r)\label{k10}\ee
%\be \int_{d(x_0,x)>R}|P(x)|^{\si\b'}d\mu(x)\le CV^{1-\si\b'/\b},\label{k9}\ee

\smallskip
Given $\e_0<\e<r$, let
\be \phi_{\e,r}(y)=k(x_0,y)|k(x_0,y)|^{\b-2}\chi^{}_{B(x_0,r)\setminus B(x_0,\e)}(y),\label{s0}\ee
%where $0<\e<r$. 
and let 
\be \wt{\phi}_{\e,r}(y)=\phi_{\e,r}(y)-Q\phi_{\e,r}(y)\nt\ee
where 
\be Q\phi_{\e,r}(y)=\int_{B(x_0,r)}Q(y,z)\phi_{\e,r}(z)d\mu(z).\ee
Clearly $\wt\phi_{\e,r}$ is orthogonal to $P(x,\cdot)$, for any $x\in B(x_0,r)^c.$
By \eqref{k3}, \eqref{Q}, and \eqref{vstar} we have 
\be| Q\phi_{\e,r}(y)|\le \frac{C}{V(r)}\int_{B(x_0,r)} |k(x_0,z)|^{\b-1}d\mu(z)\le CV(r)^{-1/\b'}.\label{s22}\ee

We now prove the following estimates:
\be A\log\frac{V(r)}{V(\e)}-C\le\|\wt{\phi}_{\e,r}\|^{\b'}_{\b'}\le A\log\frac{V(r)}{V(\e)}+C.\label{s2}\ee

  \be|T\wt{\phi}_{\e,r}(x)| \geq A\log\frac{V(r)}{V(\e)}-C,\qquad x\in B(x_0,\e')\label{s3}\ee
  \be \|T\wt{\phi}_{\e,r}\|_{p}^{p}\le CV(r).\label{s18}\ee
  
\ni  where in \eqref{s2}, \eqref{s3} we have $A=A_0$ if $k(x,y)$ is proper on the diagonal and $0<V(\e)<V(r)\le1$, and $A=A_\infty$ if $k(x,y)$ is proper at infinity, with $V(r)>V(\e)\ge1.$  

To prove \eqref{s2}, using \eqref{s22}
\be |\phi_{\e,r^{}}|-CV(r)^{-1/\b'}\le |\wt{\phi}_{\e,r^{}}|\le |\phi_{\e,r^{}}|+CV(r)^{-1/\b'},\nt\ee
hence 
\be |\phi_{\e,r^{}}|^{\b'}-C|\phi_{\e,r^{}}|^{\b'-1}V(r)^{-1/\b'}\le |\wt{\phi}_{\e,r}|^{\b'}\le |\phi_{\e,r}|^{\b'}+C|\phi_{\e,r}|^{\b'-1}V(r)^{-1/\b'}+CV(r)^{-1}.\nt\ee
(for the first inequality use that if $c>a-b$ then $c^p>a^p-pba^{p-1}$, for $a,b,c\ge0$, $p>1$, for the second, use $(a+b)^p\le a^p+p2^{p-1}(a^{p-1}b+b^p)$).
%Let \be \phi_{\e,r_0^{}}_\e(y)=k(x_0,y)|k(x_0,y)|^{\b-2}\chi^{}_{B\setminus B(x_0,\e)}(y).\label{s1}\ee
So by \eqref{k1}, \eqref{k2},   \eqref{k5}, \eqref{k7} and \eqref{k3} we obtain \eqref{s2}.
 
Note that if $x\in B(x_0,R)$, then by \eqref{s22}, \eqref{kk12} (or\eqref{vstar}) 
 \be |TQ\phi_{\e,r}(x)|\le CV(r)^{-1/\b'}\int_{B(x_0,r)}|k(x,y)|d\mu(y)\le C\label{s23}\ee
 (recall that $Q(y,z)$ is supported in $B(x_0,r)\times B(x_0,r)$).

Estimate \eqref{s3} follows from

  \be\ba &|T\wt{\phi}_{\e,r}(x)|\geq \bigg|\int_{\e\le d(x_0,y)\le r}k(x,y) k(x_0,y)|k(x_0,y)|^{\b-2}d\mu(y)\bigg|-C\cr&\geq\int_{\e\le d(x_0,y)\le r}|k(x_0,y)|^{\b}d\mu(y)-\bigg|\int_{\e\le d(x_0,y)\le r}(k(x,y)-k(x_0,y)) k(x_0,y)|k(x_0,y)|^{\b-2}d\mu(y)\bigg|-C\cr
 & \geq A\log\frac{V(r)}{V(\e)}-C\ea\ee
where the first inequality is by \eqref{s23} and in the last inequality we used \eqref{k6}.
 
Finally, to prove \eqref{s18},  take any $R>0$ such that $V(R)=2V(r)$. For any  $x\in B(x_0,R)^c$ using the orthogonality property of $\wt \phi_{\e,r}$ we have 
 $$T\wt{\phi}_{\e,r}(x)=\int_{B(x_0,r)}(k(x,y)-P(x,y))\wt\phi_{\e,r}(y)d\mu(y).$$
 Assuming  $p>(1+p/\b') (1+\eta)^{-1}$ and using \eqref{yreg2}, \eqref{vstar} we get, for $y\in B(x_0,r)$,
\be\ba \int_{d(x_0,x)>R}|k(x,y)-P(x,y)|^{p}d\mu(x)&\le  C (V(r))^{p\eta}\int_{B(x_0,r)^c}\big(V(d(x_0,x)\big)^{-p(\eta+1/\b)}d\mu(x)\\ &\le C(V(r))^{1-p/\b}.\ea\label{k8}\ee
  Therefore by Minkowski's inequality, \eqref{k8}, \eqref{s22}  
 \be\ba &\int_{d(x_0,x)>R}|T\wt{\phi}_{\e,r}(x)|^{p}d\mu(x)\cr
 &\le \int_{d(x_0,x)>R}\bigg(\int_{{B(x_0,r)}}|k(x,y)-P(x,y)|\,|\wt\phi_{\e,r}(y)|d\mu(y)\bigg)^{p}d\mu(x)   \cr 
 &\le \int_{d(x_0,x)>R}|k(x,y)-P(x,y)|^{p}d\mu(x)\,\bigg(\int_{B(x_0,r)}|\wt\phi_{\e,r}(y)|d\mu(y)\bigg)^{p}\cr
 &\le C(V(r))^{1-p/\b}\int_{B(x_0,r)} \big(|k(x_0,y)|^{\b-1}+|Q\phi_{\e,r}(y)|\big) d\mu(y)\le CV(r).
 \ea\label{s19}\ee 
 %where the last inequality is by \eqref{kk12}.
 
 On the other hand, if $x\in B(x_0,R)$, by \eqref{kk12}, $\|\phi_{\e,r}\|_q\le CV(r)^{-1/\b'+1/q}$ if $q<\b'$. So by taking $q'\geq p,\ q<\b'$ such that $(q')^{-1}=q^{-1}-(\b')^{-1}$ and applying the Hardy-Littlewood-Sobolev inequality for $T\phi_{\e,r}$, (see [A2])  we have
 \be\ba &\bigg( \int_{d(x_0,x)\le R}|T\phi_{\e,r}(x)|^{p}d\mu(x)\bigg)^{1/p}\le CV(R)^{1/p-1/q'}\bigg( \int_{d(x_0,x)\le R}|T\phi_{\e,r}(x)|^{q'}d\mu(x)\bigg)^{1/q'}\cr
 &\le CV(r)^{1/p-1/q'}\|\phi_{\e,r}\|_q\le CV(r)^{1/p}.
 \ea\label{s20}\ee
 Therefore 
 \be \int_{d(x_0,x)<R}|T\wt{\phi}_{\e,r}(x)|^{p}d\mu(x)\le C\int_{d(x_0,x)<R}(|T{\phi_{\e,r}}(x)|^{p}+|TQ\phi_{\e,r}(x)|^{p})d\mu(x)\le V(r).\nt\ee
 
 \bigskip
 
\ni{\udl{Proof of (a)}}. Suppose that $k(x,y)$ is proper on the diagonal.
Let $r_0=r_{x_0}(1)$, i.e. $B(x_0,r_0)$ is the smallest ball with  volume  $V(r_0)= 1$, and recall that $\e_0\ge0$ is the largest $\e$ with $V(\e)=0$.  
For $\e_0<\e<r_0$ let \be\psi_\e=\frac{\wt{\phi}_{\e,r_0^{}}}{\|\wt{\phi}_{\e,r_0^{}}\|_{\b'}}.\label{psie}\ee Clearly $\|\psi_\e\|_{\b'}=1$ and from \eqref{s2}, \eqref{s3} and \eqref{s18}
\be |T{\psi}_\e(x)|^\b\geq A_0\log\frac{1}{V(\e)}-C,\ \ \ x\in B(x_0,\e'),\label{s5}\ee
(recall: $V(\e)=2V(\e'))$) and 
\be \|T\psi_\e\|_{p}^{p}\le C\bigg(\!\!\log\frac{1}{V(\e)}\bigg)^{-p/\b'}\to0,\qquad \e\to\e_0.\label{s7}\ee 
For the sharpness in \eqref{1a}, note that for $\theta>1$
 \be \ba\int_{B(x_0,\e')}\frac{\exp\bigg[\dfrac{\theta}{A_0}|T\psi_\e|^\b\bigg]}{1+|T\psi_\e|^{p\b/\b'}}d\mu(x)&\geq CV(\e')V(\e)^{-\theta}\!\left(\!\log\frac{1}{V(\e)}\right)^{-p/\b'}\cr&\qquad =CV(\e)^{1-\theta}\!\left(\!\log\frac{1}{V(\e)}\right)^{-p/\b'}\to+\infty,\qquad \e\to\e_0\ea\label{s8a}\ee
 whereas for $0<\theta<1$
  \be   \frac{1}{\|T\psi_\e\|_{p}^{p}}\int_{B(x_0,\e')}\frac{\exp\bigg[\dfrac{1}{A_0}|T\psi_\e|^\b\bigg]}{1+|T\psi_\e|^{\theta p\b/\b'}}d\mu(x)\ge C \left(\!\log\frac{1}{V(\e)}\right)^{(1-\theta)p/\b'}\to+\infty,\qquad \e\to\e_0.\label{s8}\ee

To prove the sharpness for \eqref{1a1}, let 
\be \psi_\e=\frac{\wt\phi_\e}{(\|\wt\phi_\e\|_{\b'}^{\b'}+\k\|T\wt\phi_\e\|_{\b'}^{\b'})^{1/\b'}},\label{s16}\ee
then it's easy to check that  \eqref{s5}-\eqref{s8} still hold. 

For the sharpness of \eqref{1a2}, let $\psi_\e$ be defined as in \eqref{s16}.  If $\theta>1$, by \eqref{s5}, \eqref{s7} we get
\be  \int_{B(x_0,\e')}\exp\bigg[\dfrac{\theta}{A_0}|T\psi_\e|^\b\bigg]d\mu(x)\geq CV(\e')V(\e)^{-\theta}= CV(\e)^{1-\theta}\to\infty,\ \ \e\to \e_0.\label{s17}\ee

\bigskip

\ni{\udl{Proof of (b)}}.  Suppose that $k(x,y)$ is proper at infinity, and  let $r_0^{}$ be as above.  Let $0<r_1<r_0^{}<r $ with $1=V(r_0)=2V(r_1)$. 
 Let \be \psi_r=\frac{\wt{\phi}_{r_0,r}}{\|\wt{\phi}_{r_0,r}\|_{\b'}},\nt\ee 
then from \eqref{s2}, \eqref{s3}, \eqref{s18}
\be |T{\psi}_r(x)|^\b\geq A_\infty\log V(r)-C\qquad x\in B(x_0,r_1),\label{s12}\ee
and 
\be \|T\psi_r\|_{p}^{p}\le CV(r)\big(\log{V(r)}\big)^{-p/\b'}.\label{s13}\ee 
 Therefore 
 \be \frac{1}{\|T\psi_r\|_{p}^{p}} \int_{B(x_0,r_1)}\frac{\exp\bigg[\dfrac{\theta_1}{A_\infty}|T\psi_r|^\b\bigg]}{1+|T\psi_r|^{\theta_2 p\b/\b'}}d\mu(x)\geq CV(r)^{\theta_1-1}\big(\log{V(r)}\big)^{(1-\theta_2)p/\b'}.\label{s14}\ee

 Hence if either $\theta_1>1,\theta_2=1$ or $\theta_1=1,\theta_2<1$, then we have 
  \be   \frac{1}{\|T\psi_r\|_{p}^{p}}\int_{B(x_0,r_1)}\frac{\exp\bigg[\dfrac{\theta_1}{A_\infty}|T\psi_r|^\b\bigg]}{1+|T\psi_r|^{\theta_2 p\b/\b'}}d\mu(x)\to\infty,\ \ r\to \infty.\label{s15}\ee

\hfill\qed

%, since \eqref{k} is satisfied by condition \eqref{kk12}. In general $P(x)$ could be the Taylor expansion of $k(x_0,y)$ in the $y$ variable at $x$ if it exists.

\begin{rk}\label{ss1}
An example where $A_\infty>A_0$ and inequality \eqref{1a} is sharp. 
\end{rk}

Consider a Riesz-like kernel on $\R^n$ such that $k\in C^n(\R^n\setminus 0)$ and
\be k(x,y)=k(|x-y|)=\bc |x-y|^{\a-n}  \ \ \ &\text{if}\ |x-y|\le 1\cr 2|x-y|^{\a-n} &\text{if}\ |x-y|\geq 2.\ec\nt\ee

We will verify that the above kernel satisfies the conditions of Theorem \ref{s} with $\b=n/(n-\a),\ \b'=n/\a,\ A_0=|B_1|,\ A_\infty=2^\b|B_1|.$

  It is clear that \eqref{k1}, \eqref{k2} and \eqref{k3} are satisfied with the given $A_0,\ A_\infty$. To verify \eqref{k4}, note that $V_x(R)\geq (1+\delta)V_x(r)$ implies $R\geq (1+\delta')r$, and hence $c_1R\le |x'-y|\le c_2R$ for all $x'\in B(x,r),\ y\in B(x,R)^c$. So we have 
 $||x'-y|^{\a-n}-|x-y|^{\a-n}|\le C|x'-x||y|^{\a-n-1}\le Cr|y|^{\a-n-1}$ and \eqref{k4} follows.
 
 Let $x_0=0$. It is clear that \eqref{xreg} and  \eqref{k7} hold. Since $k(x,y)$ is a $C^n$ function, we can take its Taylor expansion of order $n-1$ so that \eqref{yreg2} holds with $\eta=1,\ m=n$. 
 
%  and \eqref{k6} also holds by similar argument for verifying \eqref{k4}. Let $P(x)=k(|x|)$, then \eqref{k10} holds and 
% \be\ba &\int_{|x|>R}(|k(x,y)-P(x)|)^{\si\b'}dx\le C\int_{|x|>R}(|x-y|^{\a-n}-|x|^{\a-n})^{\si\b'}dx\cr
% &\le C\int_{|x|>R}(|y||x|^{\a-n-1})^{\si\b'}dx= CR^{n(1-\si\b'/\b)}\ea\label{s21}\ee
% which is \eqref{k8}. 

%##########################################################
%\input{8-pdeapplication}

\section{Existence of ground state solutions for quasilinear equations}\label{pdeapp}
\ \ \

 In [MS2] Masmoudi and Sani showed that their inequality \eqref{MS1} in the case $\a=1$ can be applied to show the existence of a radial ground state solution for the following quasilinear equation in $\R^n,\ n\geq 2$:
\be -\Delta_n u +V_0|u|^{n-2}u=f(u) \label{app1}\ee
where \be \Delta_n u=\text{div}(|\nabla u|^{n-2}\nabla u)\nt\ee is the $n$-Laplacian operator,  $V_0$ is a positive constant, and the allowed maximal growth on the nonlinear term $f$ is exponential.

In this section, by making use of an improved version of \eqref{MS1} we prove an existence result for the following more general quasilinear elliptic equation:
\be \bc -\Delta_n u +V_0|u|^{p-2}u=f(u)\qquad & p>1,\ n\geq2\cr
||\nabla u||_n<\infty, \ ||u||_{p}<\infty.\ec\label{app2}\ee

First, let's observe that as an immediate consequence  of Theorem 1, we obtain the improved MSI given in \eqref{MS3}, simply by writing $u\in C_c^\infty(\R^n)$ in terms of the Riesz potential $I_\a f=|y|^{\a-n}*f$. Specifically,  $u=c_\a |y|^{\a-n}*(-\Delta)^{\frac\a2}u$ for $\alpha$ even, and $u=\tilde c_\a|y|^{\a+1-n}\nabla y*\nabla(-\Delta)^{\frac{\a-1}2}u$ for $\a$ odd, where 
\be c_\a=\frac{2^{-\a}\Gamma\big({\frac{n-\a}2}\big)}{\pi^{n/2}\Gamma\big(\frac\a 2\big)},\qquad \tilde c_\a=\frac{2^{-\a}\Gamma\big({\frac{n-\a+1}2}\big)}{\pi^{n/2}\Gamma\big(\frac{\a+1} 2\big)}.\ee
 Theorem 1 then yields \eqref{MS3} with 
 \be \gam_{n,\a}=
\begin{cases}
c_\a^{-\frac{n}{n-\a}}|B_1|^{-1}\ \ \ &\text{if} \ \a \ \text{even}\\
\tilde c_\a^{\;-\frac{n}{n-\a}}|B_1|^{-1}\ \ &\text{if} \  \a\ \text{odd},
\end{cases}
\label{mf4}
\ee
and where $|B_1|$ is  the measure of the unit ball. 

In order to extend the validity of \eqref{MS3} to a wider space, we  define the {\it mixed Sobolev space}

\be W^{\a,q,p}(\R^n):=\{u\in L^p(\R^n):\; |\nabla^\a u|\in L^q(\R^n)\}\label {mix}\ee
which can be easily seen to be a reflexive Banach space, under the norm $\|u\|_p+\|\nabla^\a  u\|_q$, and obviously $W^{\a,q,q}(\R^n)$ coincides with the usual Sobolev space $W^{\a,q}(\R^n)$. It is also easy to check using standard arguments (for ex. as in [AF, Thm. 3.22])  that $C_c^\infty(\R^n)$ is a dense subspace of $W^{\a,q,p}(\R^n)$. 

We record the full statement of the improved MSI,  equivalent to  \eqref{MS3},  here:

\begin{theorem}\label{msip}
Let $0<\a<n$ be an integer, and let $p\ge1$. There exists $C$ such that for every $u\in  W^{\a,\frac{n}{\a},p}(\R^n)$ with $\; \|\nabla^\a u\|_{n/\a}\le 1$
we have
\be \int_{\R^n}\frac{\exp_{\left\lceil p\frac{n-\a}{n}-1\right\rceil}\left[\gam_{n,\a}|u|^{\frac{n}{n-\alpha}}\right]}{1+|u|^{\frac{ p\a}{n-\a}}}dx\leq C||u||_{p}^{ p} \label{msip1}\ee
where $\gam_{n,\a}$ is given in \eqref{mf4} and it is sharp.\end{theorem}

In particular, we have, for all $u\in W^{1,n,p}(\R^n)$ such that $\|\nabla u\|_n\le 1$ and $p\ge 1$
\be\int_{|u|\geq 1} \frac{\exp[\gam_{n,1} |u|^{n/(n-1)}]}{ |u|^{p/(n-1)}}dx\le C||u||_p^p.\label {MSIMP}\ee
and using this we obtain that  $W^{1,n,p}$ is continuously embedded in  $W^{1,n}$ for $1\le p<n$, and $W^{1,n}$ continuously embedded in $W^{1,n,p}$ if $p>n$ (consider the cases $|u|<1$ and $|u|\ge 1$). Similar such embeddings can be obtained for general $\alpha$, using \eqref{msip1}.

\smallskip
The main theorem of this section is the following: 

\begin{theorem}\label{pde}
Let $f : \R\to\R$ be a continuous function satisfying $f(0)=0$ and  such that the following hold:
\be \text{There exists}\ \mu>p\ \text{such that }\ 0< F(t) := \int_0^t f(s)ds\le \frac{t}{\mu}\, f(t)\ \ \forall t\in\R\setminus \{0\},\label{app3}\tag{$f_1$}\ee
\be \text{There exist}\ t_0,M_0>0\ \text{such that }\ F(t)\le M_0f(t)\ \ \forall t\geq t_0,\label{app4}\tag{$f_2$}\ee
\be \text{There exist}\ \gam_0\geq 0\ \text{such that} \ \lim_{t\to\infty}{}{\exp[-{\gam\, t^{n/(n-1)}}]}f(t)=\bc 0\ \ \ &\text{if}\ \gam>\gam_0,\cr \infty  \ \ &\text{if}\ \gam<\gam_0.\ec\label{app5}\tag{$f_3$}\ee
\be  \lim_{t\to\infty}{ t^{p/(n-1)}}\exp[-\gam_0 t^{n/(n-1)}]F(t)=\infty.\label{pde6}\tag{$f_4$}\ee
Then for each $V_0>0$, \eqref{app2} has a positive radial solution $u_0\in W^{1,n,p}(\R^n)$
which is also a ground state solution.
\end{theorem}
In [MS2] the above theorem was proved for $p=n$. 
A sufficient condition for \eqref{pde6} (from l'Hospital rule) is 
\be  \lim_{t\to\infty}{ t^{(p-1)/(n-1)}}\exp[-\gam_0 t^{n/(n-1)}]f(t)=\infty\label{pde6f5}\tag{$f_5$}.\ee
Given any $\gam_0>0$, a typical function $f$ for which  all of the above conditions are satisfied is given by  
\be f(t)=\bc t^{-\e}  e^{\sgam_0 t^{\frac{n}{n-1}}} & \text {if }\; t\ge t_0\cr t^{\mu-1} & \text {if } \;0<t<t_0\cr\ec\ee
where $\mu>p$, $\;0\le \e<\frac{p-1}{n-1}$, and where $t_0$ is any positive number.

In general conditions \eqref{app3}-\eqref{app5} imply
\be F(t)\le C t^\mu,\qquad  0<t<1\label{F1}\ee
\be F(t)\le C_\sgam \, t^{-\frac{p}{n-1}} e^{\sgam t^{\frac{n}{n-1}}}, \qquad \gam>\gam_0, \quad  t\ge 1
\label{F2}\ee
while \eqref{pde6} can be rewritten as
\be F(t)=C_0(t) t^{-\frac{p}{n-1}} e^{\sgam_0 t^{\frac{n}{n-1}}},\qquad \lim_{t\to+\infty}C_0(t)=+\infty.\label{F3}\ee

\begin{rk}\label{cgamma} The estimate in  \eqref{F2} can be refined  a little: $F(t)\le C_\sgam(t) t^{-\frac{p}{n-1}} e^{\sgam\, t^{\frac{n}{n-1}}}$, where $C_\sgam(t)\to 0$ as $t\to+\infty$.
\end{rk}

The proof of Theorem \ref{pde} is essentially the same as the proof of [MS2], with some small changes. We will therefore  outline the main steps of the  argument, highlighting the changes and the role of the crucial inequality in Theorem \ref{MSIMP}.\par

\medskip
 \ni {\bf Step 1.}  Consider the functional associated to the variational formulation of \eqref{app2} i.e.
\be I(u)=\frac{1}{n}||\nabla u||_n^n+ \frac {V_0}{p}||u||_{p}^{p}-\int_{\R^n}F(u)dx,\qquad u\in W^{1,n,p}.\nt\ee
Any  critical point of $I(u)$ is a solution to the quasilinear equation \eqref{app2}.\par
A solution $u_0$ is called a \emph{ground state solution} if it satisfies 
\be I(u_0)=m:=\inf\{I(u)\ :\ u\in W^{1,n,p}\setminus\{0\},\, u\ \text{a solution of \eqref{app2}}\},\nt\ee
i.e. $u_0$ is a  a solution with the least energy among all solutions of \eqref{app2}.

\smallskip

 \ni {\bf Step 2.} Observe that if $u$ is a solution of \eqref{app2}, then it must satisfy the so-called {\it Pohozaev identity}, which in this case is given as 
 \be V_0||u||_{p}^{p}-p\int_{\R^n}F(u)dx=0.\label{pohoz}\ee
The identity can be proved by arguing along the same lines as  in   [BL, Proposition 1]. In the case $u$ radial (which is the case of interest here) the proof is easier:  multiply \eqref{app2} by $r^nu'(r)$ and integrate by parts. This step is not necessary, strictly speaking, but it provides a motivation for the next steps.
%Define the functional 
% \be G(u)= V_0||u||_{p}^{p}-p\int_{\R^n}F(u)dx,\nt\ee
%which is the Pohozaev identity of the equation \eqref{app2}. Then we have that any solution to \eqref{app2} must satisfy the condition $G(u)=0$. 

\bigskip
 \ni {\bf Step 3.} Observe that there exists at least one $u\in W^{1,n,p}(\R^n)\setminus \{0\}$ such that the Pohozaev identiy is verified and that if  there exists one such  $u_0$ solving the constrained minimization problem
\be \ba I(u_0)&=m_0 :=\inf\big\{I(u)\ :\ u\in W^{1,n,p}\setminus\{0\},\, \text {and \eqref{pohoz} holds} \big\}\cr& =\inf\bigg\{\frac{1}{n}||\nabla u||_n^n\ :\ u\in W^{1,n,p}\setminus\{0\},\, \frac{p}{||u||_{p}^{p}}\int_{\R^n}F(u)dx=V_0\bigg\}\le m
\ea\label{min}\ee
then 
$ \wt u_0(x)=u_0\big(x(1-p\theta)^{-1/n}\big)$ is a ground state solution, where $\theta$ is a Lagrange multiplier for the constrained problem. Hence the problem is reduced to showing that  $m_0$ is attained at a radial function. Note that using the P\'olya-Szeg\fH{o} inequality, in the above constrained minimization problem
 it is enough to consider those $u$ which are radially decreasing.\bigskip

\ni{\bf Step 4.} 

Define for $D>0$
\be Q(D):=\sup\bigg\{\frac{p}{||u||_{p}^{p}}\int_{\R^n}F(u)dx\ :\ ||\nabla u||_n\le D, \ ||u||_{p}<\infty,\ u\neq 0\bigg\}\nt\ee
and show that for $\gam_0>0$
\be \sup \{D>0\ :\ Q(D)<\infty \} =\bigg(\frac{\gam_{n,1}}{\gam_0}\bigg)^{\frac{n-1}{n}}:=D^*.\label{app9}\ee
%(note that if $\gam_0=0$, then $Q_D(F)<\infty$ for any $D>0$. Hence $C(F)=\infty$ and $m_0$ is always attained for any $V_0>0).$

To see this, note that from \eqref{F1}, \eqref{F2} and the
 improved Moser-Trudinger inequality inequality with exact growth condition \eqref{MSIMP}, we have for each fixed $D<D^*=(\gam_{n,1}/\gam_0)^{(n-1)/n}$
\be \int_{\R^n} F(u)dx \le C\int_{|u|\le 1}|u|^p+C_D \int_{|u|\geq 1} \frac{\exp\left[\gam_{n,1}D^{-n/(n-1)} |u|^{n/(n-1)}\right]}{ |u|^{p/(n-1)}}dx\le C||u||_p^p\nt\ee whenever  $ ||\nabla u||_n\le D$ and $ ||u||_p<\infty$.
Hence we can deduce that  $Q(D)<\infty$ when $D<D^*.$
On the other hand, if $\ds D=D^*$ , using \eqref{F3} and the usual Moser extremal sequence  $u_{\e}$ (which is constant on $B(0,\e)$ and goes to $+\infty$ as $\e\to0$), we have that 
\be \frac{p}{||u_{\e}||_{p}^{p}}\int_{\R^n}F(u_{\e})dx\ge\frac{p}{||u_{\e}||_{p}^{p}}\int_{|x|<\e}C_0(|u_\e|)\frac{\exp[\gam_0 |u_\e|^{n/(n-1)}]}{ |u_\e|^{p/(n-1)}}dx \to \infty\qquad \text{as}\ \e\to0^+\nt\ee
with $||u_{\e}||_{p}\le 1$ and $\ds ||\nabla u_{\e}||_n=D^*$.
Therefore $Q(D^*)=+\infty\label{qdfcf2}$ and  \eqref{app9} holds.

\bigskip

\ni{\bf Step 5.} Show that for any $V_0>0$
\be m_0< \frac{(D^*)^n}{n}.\label{app12}\ee

\ni To prove \eqref{app12}, first note that since $Q(D^*)=\infty$, then there exists $u_0\neq 0$ with $||u_0||_{p}<\infty$ and $||\nabla u||_n\le D^*$ such that 
\be V_0<\frac{p}{||u_0||_{p}^{p}}\int_{\R^n}F(u_0)dx,\label{app13}\ee
and using this estimate the proof of \eqref{app12} is exactly the same as in [MS2, Proposition 7.2] 

\bigskip
\ni{\bf Step 6.} Take a sequence of radially decreasing functions  $\{u_k\}$ which is minimizing for \eqref{min}, i.e. 
$$\lim_{ k\to+\infty}\frac 1{n}\|\nabla u_k\|_n^n=m_0,\qquad\|u_k\|_p=1,\qquad  p\int_{\R^n} F(u_k)dx=V_0$$
and such that, moreover, $u_k$ converges weakly to $u_0$ in $W^{1,n,p}(\R^n)$. Then prove that $m_0>0$ and 
\be\frac{1}{n}\|\nabla u_0\|_n^n=m_0,\qquad p\int_{\R^n} F(u_0)=V_0\|u_0\|_p^p.\label{u0}\ee

\ni The key estimate here is again based on the improved MSI  \eqref{MSIMP}. In particular,  since 
$$\|\nabla u_k\|_n^n\to n m_0< (D^*)^n=\left(\frac{\gam_{n,1}}{\gam_0}\right)^{n-1}$$
then we can assume that for some $D>nm_0$ and for all $k$  $$nm_0\le \|\nabla u_k\|_n^n<D^n<\left(\frac{\gam_{n,1}}{\gam_0}\right)^{n-1}.$$
 Hence, Remark \ref{cgamma}   applied to $\gam=\gam_{n,1}/D^{n/(n-1)}>\gam_0$,  gives that for any $L>0$

\be\int_{|u_k|>L} F(u_k)dx\le \int_{|u_k|>L}C_D(|u_k|) \;\frac{\exp\left[\gam_{n,1}D^{-n/(n-1)} |u_k|^{n/(n-1)}\right]}{ |u_k|^{p/(n-1)}}dx\ee
where $C_D(t)\to 0$ as $t\to+\infty$. This means, using  \eqref{MSIMP}, that for each $\e>0$
$$\int_{|u_k|>L} F(u_k)dx\le C\e$$
if $L$ is chosen large enough, independently of $k$, and the same is true for $u$. Arguing as in [MS2, Sect. 5], one can then show that $\int_{\R^n} F(u_k)\to \int_{\R^n} F(u)=V_0/n$ and  complete the proof arguing as in  [MS2, Prop. 7.1].

Note that in the above proof we assumed $\gam_0>0$. If $\gam_0=0$ then $Q(D)<\infty$ for all $D>0$ hence the supremum in \eqref{app9} is $D^*=+\infty$ and the proof still works.
Also note that the condition \eqref{pde6} is only needed to guarantee $Q(D^*)=+\infty$, in case $\gam_0>0$. Without assuming \eqref{pde6} one can still conclude that a radial ground state solution exists when $0<V_0<Q(D^*)$, when $Q(D^*)<\infty$, where $D^*$ is given still as in \eqref{app9} (see [MS2, Thm. 7.4]).

%###########################################################
%\input{9-heisenberg}

\section{Moser-Trudinger inequalities with exact growth conditions on Heisenberg group}\label{hsbg}
%The Lie algebra of the Heisenberg group is a vector space defined as follows. 
%\be T=\frac{\partial}{\partial t},\qquad X_j=\frac{\partial}{\partial x_j}+2y_j\frac{\partial}{\partial t},\qquad Y_j=\frac{\partial}{\partial y_j}-2x_j\frac{\partial}{\partial t}\nt\ee
%for $j=1,...,n$. The above $2n+1$ vector fields form a basis for the vector space of left-invariant vector fields on $\H^n$. Note we have the commutation relations
%\be [Y_j,X_k]=4\delta_{jk}T\nt\ee
%and all other commutators vanish.\\
%The sub-Laplacian on $\H^n$ is 
%\be \mathcal{L}_0=\frac{1}{4}\sum_{j=1}^n(X_j^2+Y_j^2).\nt\ee
%And the subelliptic gradient is the $2n$ dimensional vector
%\be \nabla_{\H^n} u=(X_1u, Y_1u,...,X_nu,Y_nu).\nt\ee

The Heisenberg group will be denoted as $\H^n$, with elements $x=(z,t)\in   \C^n\times \R$ and endowed with the usual group law, dilation, and Haar measure $dx=dzdt.$ The norm in $\H^n$ is denoted as 
$ |(z,t)|=(|z|^4+t^2)^{\frac{1}{4}}$ and the homogeneous dimension is $Q=2n+2.$
The unit sphere is denoted as  $\Sigma=\{x\in \H^n\ :\ |x|=1\}.$ Given any $x\in\H^n$ we will let 
$$x^*=(z^*,t^*)=\frac x{|x|}=\bigg(\ds\frac{z}{|x|},\frac{t}{|x|^2}\bigg)\in\Sigma.$$

The sublaplacian on $\H^n$ is defined as 
\be \mathcal{L}=-\frac{1}{4}\sum_{j=1}^n(X_j^2+Y_j^2),\nt\ee
where  $X_1,..X_n,Y_1,...,Y_n,T$ is the usual basis of the left-invariant vector fields on $\H^n$, and
the horizontal gradient is 
\be \nabla_{\H^n} =(X_1, Y_1,...,X_n,Y_n).\nt\ee\def\L{\mathcal {L}}
For $1\le p,q<\infty,$ and $\a\ge0$ an even integer we define {\it mixed Sobolev space} on $\H^n$ as 
\be W^{\a,q,p}(\H^n):=\big\{u\in L^p(\H^n)\ :\ \mathcal{L}^{\frac{\a}{2}} u\in L^q(\H^n)\big\},\nt\ee
with the graph norm $||u||_p+||\mathcal{L}^{\frac{\a}{2}} u||_q$. This defines a Banach space, which coincides with the usual Sobolev space $W^{\a,q}(\H^n)$ if $p=q.$ Also, it's possible to check that $C_c^\infty(\H^n)$ is dense in $W^{\a,q,p}(\H^n)$, by the same arguments one can use for $W^{\a,q,p}(\R^n)$ (see Section \ref{pdeapp}).

\def\am{\frac{\a}{2}}

 If $\a<Q$  the fundamental solution of the powers of the $\H^n$ sublaplacian, is homogeneous of degree $\a-Q$, smooth away from 0, and  can be written as 
\be \mathcal{L}^{-\am}(x,0)=g_\a(\theta)|x|^{\a-Q} \label{lap1}\ee
where  \be \theta=\theta(t^*)=\textup{arg} \frac{|z|^2+i t}{|x|^2}\in\bigg[-\frac{\pi}{2},\frac{\pi}{2}\bigg]. \nt\ee
The function $g_\a(\theta)$, which is bounded,  was computed explicitly  (in integral form except a few cases)  in [BDR], [CT], and when $\alpha=2$ reduces to  $g_2(\theta)=2^{n-2}\Gamma(n/2)^2\pi^{-n-1}$, derived by Folland in [Fol].

\begin{theorem}\label{hs2}
Let $0<\a<Q$ be an even integer, and let $p\ge1$. There exists $C$ such that for every $u\in  W^{\a,\frac{Q}{\a},p}(\H^n)$ with $\; \|\mathcal{L}^{\frac{\a}{2}}u\|_{Q/\a}\le 1$
we have
\be \int_{\H^n}\frac{\exp_{\left\lceil p\frac{Q-\a}{Q}-1\right\rceil}\left[\dfrac{1}{A_\a}|u|^{\frac{Q}{Q-\alpha}}\right]}{1+|u|^{\frac{ p\a}{Q-\a}}}dx\leq C\|u\|_{p}^{ p} \label{hs2b}\ee
where
\be A_\a=\frac{1}{Q}{\int_{\Sigma}|g_\a(\theta)|^{\frac{Q}{Q-\alpha}}dx^*}.\label{hs2ag}\ee
Moreover, for every $u\in  W^{\a,\frac{Q}{\a}}(\H^n)$ with \be \|u\|^{Q/\a}_{Q/\a}+\|\mathcal{L}^{\frac{\a}{2}}u\|^{Q/\a}_{Q/\a}\le 1\ee
we have 
\be \int_{\H^n}{\exp_{[\frac{Q}{\alpha}-2]}\left[\frac{1}{A_\a}|u|^{\frac{Q}{Q-\alpha}}\right]}dx\leq C \label{rufhsp1}\ee
The exponential constant $A_\a^{-1}$ in \eqref{hs2b} and \eqref{rufhsp1}  is sharp, i.e. it cannot be replaced by any larger number.
\end{theorem} 
Note that the relation between the surface measure $dx^*$ in the above \eqref{hs2ag} and $\theta$ is as follows [BFM, (2.4)]
\be \int_{\Sigma}\phi dx^*=\omega_{2n-1}\int_{-\pi/2}^{\pi/2}\phi(\theta)(\cos\theta)^{n-1}d\theta, \nt\ee
where $\phi\ :\ [-\frac{\pi}{2},\frac{\pi}{2}]\to\R$ is a measurable function.\par

\bigskip
In the case of the horizontal gradient, let

\be W^{1,q,p}(\H^n):=\big\{u\in L^p(\H^n)\ :\ |\nabla_{\H^n} u|\in L^q(\H^n)\big\},\nt\ee
with the graph norm $||f||_p+||\nabla_{\H^n}  u ||_q$. Once again this is a Banach space, which coincides with the usual Sobolev space $W^{1,q}(\H^n)$ if $p=q,$ and which contains $C_c^\infty(\H^n)$ as a dense subspace.
  
\medskip
\begin{theorem}\label{hs3}
For every $p\ge1 $, there exists $C$ such that for every $u\in W^{1,Q,p}(\H^n)$ with $\;\|\nabla_{\H^n} u||_{Q}\le 1$
we have 
\be \int_{\H^n}\frac{\exp_{\llceil p\frac{Q-1}{Q}-1\rrceil}\left[\dfrac{1}{A}|u|^{\frac{Q}{Q-1}}\right]}{1+|u|^{\frac{ p}{Q-1}}}dx\leq C||u||_{ p}^{p} \label{hs3b}\ee
where
\be A=\frac{1}{Q}\left(\int_{\Sigma} |z^*|^Q dx^* \right)^{-\frac{1}{Q-1}}=\frac{1}{Q}\left(\frac{\Gamma{(\frac{Q}{2})}\Gamma{(n)}}{2\pi^{n+1/2}\Gamma{(\frac{Q-1}{2})}}\right)^{\frac{1}{Q-1}},\label{A}\ee
and the exponential constant $A^{-1}$ is sharp,i.e. it cannot be replaced by any larger number. 
\end{theorem} 
For the computation of the constant $A$ in \eqref{A},  see [CL].\\

\bigskip

\ni{\bf Proof of Theorem \ref{hs2}.} The inequalities are immediate consequences of Theorem 1 as applied to the operator $T$ with kernel $k(x,y)=K(y^{-1}x)$, where $K(x)=g_\a(\theta)|x|^{\a-Q}$, and with $\beta=\frac{Q}{Q-\a}$. Indeed, it is enough to prove the inequalities for $u$ smooth and compactly supported, in which case we can write  $u=Tf$ with $f=\L^{\am}u$ in $L^{Q/\a}$ and compactly supported.
All one needs to check is that the conditions (K1)-(K4) are verified. Estimate (K3) is obvious since $g_\a(\theta)$ is bounded and $|B(0,r)|=C r^Q$. Using polar coordinates one readily gets
\be \int_{r_1<|x|<r_2} |K(x)|^{\frac{Q}{Q-\a}}dx=A_\a\log{ \frac{|B(0,r_2)|}{|B(0,r_1)|}},\qquad r_1<r_2\label{KH}\ee
which gives (K1) and (K2). Regarding (K4) we in fact have both \eqref{xreg} and \eqref{yreg2} (and also \eqref{yreg3}) simply from Taylor's formula and the fact that 
$$|D_z^h D_t^\ell K(x)|\le C|x|^{\a-Q-|h|-2\ell}$$
for any multi-index $h=(h_1,...,h_{2n})$ with $|h|=h_1+...+h_n$ (and the obvious notation for $D_z^\ell$). 

The sharpness statements follow from Theorem 2 and the pointwise  regularity estimates on $K$. Indeed, one just has to  consider the family of functions  $u_\e=T\psi_\e$, where $\psi_\e$ is defined as in \eqref{psie}.\hfill\QEDopen
\bigskip

\ni{\bf Proof of Theorem \ref{hs3}.} 
Let $x=(z,t),\ y=(z',t')$. In the case of subgradient, the following representation formula was derived by [CL]:
\be u(x)=-\frac{1}{4c_Q^{}}\int_{\H^n}\frac{|z'|^{Q-2}}{|y|^{2Q}}\nabla_{\H^n}(|y|^4)\cdot \nabla_{\H^n}F(y^{-1}x) dy\label{hs2c}\ee
where ``$\ \cdot\ $'' denotes the standard inner product of two vectors in $\R^{2n}$, $c_Q=\int_{\Sigma}|z^*|^Q dx^*$ and $u\in C_0^\infty(\H^n)$. Hence, we let
\be K(x)=-\frac{1}{4c_Q^{}}\frac{|z|^{Q-2}}{|x|^{2Q}}\nabla_{\H^n}(|x|^4),\ \ x\ne 0,\nt\ee
 and since $|\nabla_{\H^n}(|x|^4)|=4|z||x|^2$ (see [CL, Lemma 2.1] then
\be |K(x)|=c_Q^{-1}|z|^{Q-2}|x|^{2-2Q}=c_Q^{-1}|z^*|^{Q-2}|x|^{1-Q}:=g(x^*)|x|^{1-Q}.\nt\ee
% First, we compute
%\be \nabla_{\H^n}(|u|^4)=4(|z|^2x_1+y_1t, |z|^2y_1-x_1t,...,|z|^2x_n+y_nt, |z|^2y_n-x_nt), \nt\ee
%and $|\nabla_{\H^n}(|u|^4)|=4|z||u|^2.$ Then we can rewrite $K(u)$ as follows
%\be \ba K(u)&=\frac{-(c_Q)^{-1}}{4}\frac{|z|^{Q-2}}{|u|^{Q+1}}4(|z|^2x_1+y_1t, |z|^2y_1-x_1t,...,|z|^2x_n+y_nt, |z|^2y_n-x_nt)|u|^{1-Q}\cr
%&=-(c_Q)^{-1}\frac{|z|^{Q-2}}{|u|^{Q-2}}\bigg(\frac{|z|^2x_1+y_1t}{|u|^3}, \frac{|z|^2y_1-x_1t}{|u|^3},...,\frac{|z|^2x_n+y_nt}{|u|^3}, \frac{|z|^2y_n-x_nt}{|u|^3}\bigg)|u|^{1-Q}\cr
%&=-(c_Q)^{-1}|z^*|^{Q-2}(|z^*|^2x^*_1+y^*_1t^*, |z^*|^2y^*_1-x^*_1t^*,...,|z^*|^2x^*_n+y^*_nt^*, |z^*|^2y^*_n-x^*_nt^*)\cr
%&=:g(u^*)|u|^{1-Q}.\ea\nt\ee
Therefore, arguing as in \eqref{KH}  via polar coordinates
\be \int_{r_1<|x|<r_2} |K(x)|^{\frac{Q}{Q-\a}}dx\le A \log{ \frac{|B(0,r_2)|}{|B(0,r_1)|}},\qquad r_1<r_2\label{KH1}\ee
with 
\be A=\frac{1}{Q}\int_\Sigma |g(x^*)|^{\frac{Q}{Q-1}}dx^*=\frac{1}{Q}(c_Q)^{-\frac{Q}{Q-1}}\int_\Sigma|z^*|^{(Q-1)\frac{Q}{Q-1}}dx^*=\frac{1}{Q}(c_Q)^{-\frac{1}{Q-1}}.\nt\ee
The inequality in \eqref{hs3b} follows now from Theorem 1, applied to the  kernel $|K(x)|$, which also satisfies (K4) since it is of the same type discussed in the previous proof.

 To prove the sharpness, we use the same extremal family of functions as in [CL]. We consider the function
\be v_\e(y)=\bc \log(1/|y|)\qquad&\text{for}\ \e\le|y|\le 1\cr \log(1/\e)&\text{for}\ |y|\le \e.\ec\nt\ee
Then [CL] computed that
\be ||\nabla_{\H^n}v_\e||_Q^Q=c_Q^{}\log\frac{1}{\e},\ee
and it is clear that $ ||v_\e||_{\s Q}\le C$. Let
\be u_\e=\frac{v_\e}{||\nabla_{\H^n}v_\e||_Q},\nt\ee
so that $ ||\nabla_{\H^n}u_\e||_Q= 1$
\be ||u_\e||_{p}^{p}\le C\left(\log\frac{1}{\e}\right)^{p/Q}\to0,\qquad \e\to0^+,\nt\ee
and \be |u_\e|^{\frac{Q}{Q-1}}\geq c_Q^{-\frac{1}{Q-1}}\log\frac{1}{\e}\qquad\text{for}\ |y|\le\e.\nt\ee
If  $\theta >1$ then
%%\be \ba &\int_{\H^n}\frac{\exp_{[\frac{n}{\alpha}-2]}\big[\theta A^{-1}|u_\e|^{\frac{Q}{Q-1}}\big]}{1+|u_\e|^{\frac{q}{Q-1}}}dy\cr
%&\geq \int_{|y|\le\e}\frac{\exp\big[\theta A^{-1}|u_\e|^{\frac{Q}{Q-1}}\big]}{1+|u_\e|^{\frac{q}{Q-1}}}dy\cr
%&\geq C\e^Q\frac{\exp\bigg[\theta\log\dfrac{1}{\e^Q}+C\bigg]}{1+C\Big(\log\dfrac{1}{\e}\Big)^{q/Q}\cr
%&\ge C\e^{(1-\theta)Q}(\log\dfrac{1}{\e})^{-q/Q}\to\infty \ea\nt\ee
%%
%%\hfill\QEDopen
\be \ba \int_{\H^n}\frac{\exp_{[Q-2]}\big[\theta A^{-1}|u_\e|^{\frac{Q}{Q-1}}\big]}{1+|u_\e|^{\frac{p}{Q-1}}}dy&\geq \int_{|y|\le\e}\frac{\exp\big[\theta A^{-1}|u_\e|^{\frac{Q}{Q-1}}\big]}{1+|u_\e|^{\frac{p}{Q-1}}}dy\geq C\e^Q\,\frac{\exp\bigg[\theta\log\dfrac{1}{\e^Q}+C\bigg]}{1+C\Big(\log\dfrac{1}{\e}\Big)^{p/Q}}\cr&\ge C\e^{(1-\theta)Q}\left(\log\dfrac{1}{\e}\right)^{-p/Q}\to\infty \ea\nt\ee
as $\e\to0^+.$

\hfill\QEDopen

\par

%##########################################################
%\input{10-manifolds}

\section{Masmoudi-Sani inequalities on Riemannian manifolds with negative curvature}\label{mf}
 Let $M$ be a complete, connected, smooth Riemannian manifold with metric tensor $g$. The geodesic distance between two points $x,y\in M$ is denoted by $d(x,y)$, which gives a metric in $M$. The Riemannian measure $\mu$ associated to $g$ is defined in local coordinates as $d\mu=\sqrt{|g|} dm$, where $|g|=\det(g_{ij})$, $g=(g_{ij})$, and $dm$ is the  Lebesgue measure. We denote the inverse matrix of $(g_{ij})$ by $(g^{ij})$.
 
For vector fields $Z,W$ we will let 
 \be \lan Z,W\ran=g(Z,W),\qquad |Z|=\lan Z,Z\ran^{1/2}.\label{norm}\ee
 \def\ric{{\rm{Ric } }} Covariant derivatives of order $k$ of a smooth function $u$ are denoted as $\nabla^k u$, and their components in a local chart are denoted as $(\nabla^k u)_{j_1...j_k}$. In particular $\nabla^0 u=u$, $(\nabla^1 u)_j=\partial_j u$. Define
 $$|\nabla^k u|^2=g^{i_1j_1}....g^{i_kj_k}(\nabla^k u)_{i_1...i_k}(\nabla^k u)_{j_1...j_k}.$$

The gradient and the Laplacian of  a  smooth function $u$ on $M$ are given as  
\be \nabla u=g^{ij}\partial_j u \,\partial_i,\qquad \Delta u=g^{ij}(\nabla^2 u)_{ij}=\frac{1}{\sqrt{|g|}}\,{\partial_i}\left(\sqrt{|g|}\,g^{ij}\partial_j u  \right).\nt\ee
and 
$$|\nabla u|^2=g(\nabla u, \nabla u)=g_{ij}(\nabla u)^i(\nabla u)^j=g^{ij}\partial_i u\, \partial_j u=|\nabla^1 u|^2.$$
Let us define  
\be D^\a=\bc (-\Delta)^{\frac{\a}{2}}  &\text{  for $\alpha$ even}\cr\nabla(-\Delta)^{\frac{\a-1}{2}} & \text{ for $\alpha$ odd.} 
\ec\ee

 \ni The sectional curvature in the metric $g$ will be denoted as $K$ and the Ricci curvature as $\ric$.

The Bessel potential space is defined for $\a>0,p\ge1$  as 
\be H^{\a,p}(M)=\{u\in L^p(M): \; u=(I-\Delta)^{-\frac\a2}f,\; f\in L^p(M)\}\label {bessel1}\ee
endowed with the norm $\|(I-\Delta)^{\frac\a2}u\|_p$. These spaces were studied, on general noncompact manifolds, by Strichartz [St] who also showed  that $C_c^\infty(M)$ is dense in  $H^{\a,p}(M)$ and that $\|u\|_p+\|(-\Delta)^{\frac\a2}u\|_p$ is an equivalent norm if $p>1$. On Hadamard manifolds, work of Auscher et al [ACDH] shows  the Riesz transform is continuous in $L^p$, for $1<p<\infty$,  and this implies that $\|u\|_p+\|D^\a u\|_p$
is an equivalent norm in $H^{\a,p}$, for $\a\ge1.$  Further, the Poincar\'e inequality on such manifolds implies that $\|D^\a u\|_p$ is also an equivalent norm in $H^{\a,p}$.

\begin{rk}\label {sob} The classical Sobolev spaces $W^{\a,p}(M)$ and $W_0^{\a,p}(M)$ for $\alpha$ integer are  defined in terms of norm $\|u \|_{W^{\a,q}}=\sum_0^\a\|\nabla^ju\|_p$ as follows. If $C^{\a,p}(M)=\big\{u\in C^\infty(M): \|u\|_{W^{\a,p}}<~\infty\big\}$, then $W^{\a,p}(M)$ is the closure of $C^{\a,p}(M)$ relative to the norm $\|\cdot\|_{W^{\a,p}}$, and  $W_0^{\a,p}(M)$ is the closure of $C_c^\infty(M)$ in $W^{k,q}(M)$. As it turns out if $p\ge1$ then $W_0^{\a,p}(M)=W^{\a,p}(M)$ for $\a=1$,  and for $\a\ge 2$ if $\inj(M)>0$ and $|\nabla^j \ric|\le C$ for all  $j\le \a-2$ (see e.g. [H]). Also, for $p>1$, $\a\in\N$ we have $W^{\a,p}(M)\subseteq H^{\a,p}(M)$ with equality if $\a=1$, and for all $\alpha\in\N$, $\alpha\ge2$ if $M$ satisfies $\inj(M)>0$ and $|\nabla^j R|\le C$ for all $j$, where $R$ is the Riemann curvature tensor  (see e.g. [T, p. 320]). It is conjectured that $W^{\a,p}=H^{\a,p}$ provided $\inj(M)>0$ and $|\nabla^j\ric|\le C$ for $j\le \a-2.$

\end{rk}

\begin{theorem}\label{mf1}
Let $(M,g)$ be a complete and simply connected Riemannian manifold of dimension $n\ge3$, such that $-b^2\le K\le -a^2$  for some $a,b>0$. For any integer $0<\a<n$, there exists a constant $C$ such that for all $u\in H^{\a,n/\a}(M)$ with 
\be \|D^\a u\|_{n/\a}\le 1\ee we have 
\be  \int_M \frac{\exp_{\llceil\frac{n-\a}{\a}-1\rrceil}{\left[\gam_{n,\a}|u|^{\frac{n}{n-\a}}\right]}}{1+|u|^{{\frac{n}{n-\a}}}}d\mu\le C\|u\|_{{ n/\a}}^{{n/\a}}, \label{mf3}\ee 
where $\gam_{n,\a}$ is defined in \eqref{mf4} and it is sharp.

\end{theorem}

The above theorem in the special case $K=-1$ (hyperbolic space) was derived in [LT] when $\alpha=1$, and in [NN] when $\alpha=2$, for $n\ge2$. The classical version of the above inequality, without the norm on the right-hand side, and without denominator, was derived in [BS], extending earlier results  in [FM3] on the hyperbolic space. Observe that in [BS] the result is stated and proved for $n\ge 3$, under $\ric\ge -(n-1)b^2,\; K\le -a^2$, for some $a,b>0$, which is however equivalent to $-b^2\le K\le -a^2$ for some $a,b>0$. 

\begin{rk} The proof of Theorem \ref{mf1} relies on Green function estimates that were proved in [BS] for $n\ge3$. We believe that both Theorem \ref{mf1} and the main result in [BS] are still valid in the case $n=2,\; \alpha=1$, with suitable modifications of the current proofs.

\end{rk}

\begin{rk} A version of Theorem \ref{mf1}  can be given with an inequality like  \eqref{msip1}, for any $p\ge1$. For simplicity we stated it only in the case $p=n/\a$.\end{rk}
\smallskip

\ni {\bf Proof.} It suffices to prove \eqref{mf3} for $u\in C_c^\infty(M)$. If $G_\a(x,y)$ denotes the Green function of $D^\a=(-\Delta)^{\a/2}$ for $\alpha$ even, then for each fixed $x\in M$
\be u(x)=\int_{M}G_{\a}(x,y)D^\a u(y)d\mu(y),\label{mf5}\ee
where $G_\a$ satisfies 
\be G_\a(x,y)=\int_M G_{\a-2}(x,z)G_2(z,y)d\mu(z),\qquad \alpha {\hbox{ even,}}\,\; 0<\alpha\le n.\label{itG}\ee
If $\alpha$ is odd, then applying the above to $\alpha+1$ and integrating by parts gives
\be u(x)=\int_{M}d\mu(z)\int_M G_{\a-1}(x,z)\lan\nabla_y G_2(z,y),\nabla_y D^{\a-1}u(y)\ran d\mu(y),\label{mf5a}\ee
for $\a$ odd, where $\nabla_y$ denotes the gradient in the $y$ variable.

We then  let
\be K_\a(x,y)=\bc G_\a(x,y) &\text{  for $\alpha$ even}\cr\nabla_y G_{2}(x,y) & \text{ for $\alpha=1$} \cr \int_M G_{\a-1}(x,z)\nabla_y G_2(z,y)d\mu(z)  & \text{ for $\alpha$ odd and $1<\a<n$.}
\ec\label{Ka}\ee
Note that $|K_\a(x,y)|$ is symmetric in $x,y$ for $\a$ even, but not necessarily so for $\a$ odd.

We would like to apply Theorem 1, where $T$ is the integral operator defined as in \eqref{mf5}, \eqref{mf5a}. One observation to make here is that when $\alpha$ is odd the inner product in \eqref{mf5a}  and the norms $|K_\a(x,y)|,\, |D^\a u(y)|$ are with respect to the metric $g$, and not the standard Euclidean metric, in any given chart, so Remark \ref{rkvector} following Theorem 1, does not apply directly. However, on a Hadamard manifold we can use a single coordinate chart $(M,\phi)$,  and in these coordinates we can write the metric tensor as $\big(g_{ij}(y)\big)=R_y^TR_y$, for some invertible matrix $R_y$, for each $y\in M$. Hence the operator in \eqref{mf5a} for $\alpha$ odd  can be written as
\be u(x)=Tf(x)=\int_{\R^n} R_y K_\a(x,y)\cdot R_y f(y) \sqrt{|g(y)|}dm(\xi),\label{mf5b}\ee
where $ y=\phi^{-1}\xi$ and $f=\nabla_y D^{\a-1}u(y)=D^\a u(y)$, and Remark 2 can be applied to the kernel $R_y K_\a$ and the function $R_y f$, which clearly satisfy $|K_\a|=\big((R_y K_\a)\cdot (R_y K_\a)\big)^{1/2}$, and $|f|=|D^\a u|=\big((R_y f)\cdot (R_y f)\big)^{1/2}.$ With this observation in mind we will then just check all the Riesz-like  kernel conditions using the norm in the metric $g$, when $\a$ is odd, and the usual absolute value when $\alpha$ is even. We will show  that \eqref{k1}, \eqref{k2} hold, with $A_0=\gamma_{n,\a}^{-1}$ and $A_\infty=0$. We will also show that the first estimate of \eqref{k3}  holds for $V_x(d(x,y)\le 1$ and that \eqref{kk12} holds. Finally, we will show that \eqref{k4} holds whenever $V(r)\le 1$.

\smallskip

Bertrand and Sandeep [BS, Lemma 4.1, Theorem 4.2] proved that, given the assumptions on the  curvature
\begin{numcases}{|K_\a(x,y)|\le} \big(|B_1|\gam_{n,\a}\big)^{-\frac{n-\a}{n}} d(x,y)^{\a-n}+Cd(x,y)^{\a-n+\frac{1}{2}} &\text{if} $d(x,y)<1$ \label{mf6a}\\ Ce^{-\delta d(x,y)}&\text{if} $d(x,y)\geq 1$,\label{mf6b}\end{numcases}
for some $C>0$ and $\delta>0$, and that there are $C,\delta'>0$ such that for any $x\in M$, and any $\sigma\in (0,1)$, and  for any $x\in M$
\begin{numcases}{|K_\a|^*(x,t)\le}
{\gam_{n,\a}}^{\hskip-1em-\frac{n-\a}{n}}\,t^{-\frac{n-\a}{n}}+Ct^{-\frac{n-\a}{n}+\delta'}&\text{for}  $0<t\le 1$\label{mf7a}\\
 C(\sigma)t^{-\sigma} &\text{for} $t>1$,\label{mf7b}\end{numcases}
where $C(\sigma)$ is a positive constant depending on $\sigma$ (see Remark \ref{green} below regarding these estimates).

By the classical comparison theorems, denoting $V^\lambda(r)$ to be the volume of (any) ball of radius $r$ in the $n-$dimensional space form  of constant sectional curvature $\lambda$, we have (since $\ric\ge -(n-1)b^2$)
\be V_x^0(r)=|B_1|r^n\le V_x(r)\le V^{-b^2}(r)=n |B_1| b^{-n}\int_0^{br}\big(\sinh u)^{n-1} du\ee
so that if $ r_0>0$ there is $C(r_0)$ such that for all $r\in [0,r_0]$ we have 
\be |B_1|r^n\le V_x(r)\le C(r_0) r^n \label {volest}\ee
where the left inequality holds for all $r$.
Hence for $V_x(d(x,y))\le 1$ we have $d(x,y)\le r_0:=|B_1|^{-1/n}$ and $d(x,y)\ge d_0 V_x(r)^{1/n}$, some $d_0>0$ (independent of $x$), and 
\be |K_\a(x,y)|\le C d(x,y)^{\a-n}\le C V_x(d(x,y))^{-\frac{n-\a}{n}}\label{k3small}\ee
which is the first condition in \eqref{k3}, for balls with small volumes. Also, with the definition of $r_x(\tau)$ given in \eqref{rtau} we have $0<d_0\le r_x(1)\le r_0$, and \eqref{mf6a}, \eqref{mf6b} imply that $|K_\a(x,y)|\le C$  if $V(d(x,y))\ge 1$ (since in that case $d(x,y)\ge r_x(1)\ge d_0$).

Clearly \eqref{mf7a} implies \eqref{k1p} and hence, together with \eqref{k3small}, it also implies \eqref{k1} with $A_0=\gam_{n,\a}^{-1}$  (see Remark \ref{k1prime1} in Appendix A). Additionally,  \eqref{k2p} holds with $A_\infty=0$ (critical integrability), then \eqref{k2} also holds with $A_\infty=0$, since $|K_\a(x,y)|$ bounded when $V_x(d(x,y))\ge 1$ (again by Remark \ref{k1prime1}). 

 The first estimate in \eqref{kk12}, for $\sup_x \big(|K_\a(x,\cdot)|^*\big)(t)$, follows directly from \eqref{mf7a}, \eqref{mf7b}. The second estimate in \eqref{kk12}, for $\sup_y  \big(|K_\a(\cdot,y)|^*\big)(t)$, follows from the symmetry in $x,y$ of $|K_\a(x,y)|$, when $\alpha$ is even, however it needs to be justified when $\alpha$ is odd, since symmetry is lost, in general. It's enough to argue that the estimates \eqref{mf7a},  \eqref{mf7b} are also valid for $\big(|K_\a(\cdot,y)|^*\big)(t)$, following the proof of \eqref{mf7a},  \eqref{mf7b}  given  in [BS].
Estimate  \eqref{mf7a} follows directly  from  \eqref{mf6a}, however the same method does not work for  \eqref{mf7b}  (see Remark \ref{green}). The way this was done in [BS] was to first prove 
\be|G_2|^*(x,t)\le 
C \bc t^{-\frac{n-2}{n}}&\text{for}  \ 0<t\le 1 \\
 t^{-1} &\text{for} \ t>1,\ec\label{G2}\ee
  \be|\nabla_y G_2|^*(x,t)\le 
C \bc t^{-\frac{n-1}{n}}&\text{for}  \ 0<t\le 1 \\
 t^{-1} &\text{for} \ t>1,\ec\label {GG2}\ee
 and then prove \eqref{mf7b} for all $\alpha$ by induction, via O'Neil's Lemma.
 
 Estimate \eqref{G2} for $t\le 1$ follows clearly from  \eqref{mf7a}, whereas the one for $t>1$ is a consequence of the following  distribution function estimate   ([BS, eq. (3.25]):
\be\mu\{y:G_2(x,y)>s\}\le \frac{C_{n,a}}{s},\qquad s>0\label{distrG}\ee
%\be\mu\{y:\|\nabla_y G_2(x,y)|>s\}\le \frac{C_{n,a}}{s},\qquad s>0\label{distrGG}\ee
(where $C_{n,a}$ is some explicit constant). Estimate \eqref{GG2} follows from the classical gradient estimate for  positive harmonic functions (see e.g. [Li, Thm. 6.1]) applied to  $G_2(x,y)$, in the $y$ variable inside $B(y,d(x,y))$
\be |\nabla_y G_2(x,y)|\le C\left(1+\frac 1{d(x,y)}\right) G_2(x,y)\le C\bc d(x,y)^{1-n} & \text{if }\, d(x,y)\le 1 \cr G_2(x,y) &\text {if }\, d(x,y)>1\ec\label {gradest}\ee
combined with \eqref{volest} and \eqref{distrG}. It is clear then that since the bound in \eqref{gradest} is symmetric in $x$ and $y$ then one can get the same estimate in \eqref{GG2}, and hence \eqref{kk12}, for $\big(|\nabla_y G_2(\cdot,y)|\big)^*(t)$, for any $y\in M$.

Now we will  prove that  \eqref{k4}  holds for $V_x(r)\le 1$.
Let us then assume $V_{x}(r)\le 1$ and $V_{x}(R)\ge (1+\delta)V_{x}(r)$, $x'\in B(x,r)$ and $y\in B(x,R)^c.$ We then have $r\le r_x(1)$ since on a Riemannian manifold  $V_x(r)$ is strictly increasing, and we know that $r_x(1)\le r_0=|B_1|^{1/n}$. Assume first $R>2r_0$, in which case $d(x',y)\ge d(x,y)-d(x,x')\ge R-r\ge  r_0\ge r_x(1)$, and since \eqref{k2} holds (with $A_\infty=0$) we have
\be\ba\int_{d(x,y)\ge R}&|K_\a(x',y)-K_\a(x,y)|^\nna d\mu(y) \cr &\le C \int_{d(x,y)\ge R}|K_\a(x',y)|^\nna d\mu(y) +C\int_{d(x,y)\ge R}|K_\a(x,y)|^\nna d\mu(y)\cr& \le  C \int_{d(x',y)\ge r_x(1)}|K_\a(x',y)|^\nna d\mu(y) +C\int_{d(x,y)\ge r_x(1)}|K_\a(x,y)|^\nna d\mu(y)\le C.\ea\ee
To settle the case  $R\le 2r_0$, it is enough to prove that  the following pointwise condition holds  for all $x\in M$,  if $r$ and $R$  are such that $V_{x}(r)\le 1$ and $V_{x}(R)\ge (1+\delta)V_{x}(r)$:
\be|K_\a(x',y)-K_\a(x,y)|\le C_\delta d(x,x') d(x,y)^{\a-n-1},\qquad \forall x'\in B(x,r),\;\;\forall  y\notin B(x,R),\label{regK}\ee
since it  implies  \eqref{xreg} when $d(x,y)\le 2r_0$ (due to \eqref{volest}), and hence \eqref{k4} follows.
We will prove \eqref{regK} only for $\alpha=1,2$; the proof for general $\alpha$ is obtained easily using the integral formulas in \eqref{itG} and \eqref{Ka}; the details are left to the reader.

First, note that the condition   $V_{x}(R)\ge (1+\delta)V_{x}(r)$ implies the existence of $c_\delta>0$, independent of $x,r,R$ such that $R\ge (1+c_\delta)r$. Indeed, from the standard comparison theorems $V_x(r)/V^{-b^2}(r)$ is decreasing in $r$ (and tends to 1 as $r\to0$). Hence, since $R>r$,   $V_x(R)/V^{-b^2}(R)\le V_x(r)/V^{-b^2}(r)$ and from this inequality we get 
\be V_x(R)-V_x(r)\le V^{-b^2}(R)-V^{-b^2}(r)=n|B_1|b^{-n}\int_{br}^{bR}(\sinh u)^{n-1}du\le C_0 (R^n-r^n)\ee
with $C_0$  independent of $x$, where in the last inequality we used $R\le 2r_0$. The result follows since $\delta|B_1|r^n\le \delta V_x(r)\le V_x(R)-V_x(r)$.

Integrating along the geodesic from $x$ to $x'$, inside $B(x,r)$, we get
\be\ba |G_2(x',y)-G_2(x,y)|&\le d(x,x')\sup_{z\in B(x,r)}|\nabla_z G_2(z,y)|\le C d(x,x')\sup_{z\in B(x,r)}d(z,y)^{1-n}\cr&\le C d(x,x')\bigg(\frac{c_\delta}{1+ c_\delta} d(x,y)\bigg)^{1-n}=C_\delta d(x,x')d(x,y)^{1-n},\ea\ee
which gives \eqref{regK} for $\alpha=2.$ 

To prove \eqref{regK} for $\alpha=1$, for each $y$ in a given coordinate chart let $\big(g_{ij}(y)\big)=R^T_y R_y$ for some invertible matrix $R_y$, and let $U_j(z,y)=\big(R_y\nabla_y G_2(z,y)\big)^j$, so that 
\be |\nabla_y G_2(x',y)-\nabla_y G_2(x,y)|^2=\sum_{j=1}^n \big(U_j(x',y)-U_j(x,y)\big)^2\le d(x,x')^2\sum_{j=1}^n \sup_{z\in B(x,r)}|\nabla_z U_j(z,y)|^2. \label{gradlip}\ee
Now note that for each $y$ in the chart and any $j$,  the function $U_j(z,y)$ is harmonic  in $z$ inside $B(x,d(x,y))$, in particular if $\sigma\in (0,1)$ is so that $r<\sigma d(x,y)<d(x,y)$ then
\be \min_{z\in B(x,\sigma d(x,y))} U_j(z,y)=U_j(z_y,y),\qquad {\text{some }} \; z_y\in \partial B(x,\sigma d(x,y)).\ee
It's easy to see, using the definition of Green function,  that $U_j(\cdot, y)$ cannot be constant in $B(x,\sigma d(x,y))$, hence one can apply the gradient estimate to the positive harmonic function  $U_j(z,y)-U_j(z_y,y)$ on the ball $B(x,\sigma d(x,y))$ and obtain
\be\ba|\nabla_z U_j(z,y)|&=|\nabla_z(R_y\nabla_y G_2(z,y)^j|\le C\bigg(1+{1\over d(x,y)}\bigg) |U_j(z,y)-U_j(z_y,y)|\cr & \le Cd(x,y)^{-1}\big(d(z,y)^{1-n}+d(z_y,y)^{1-n}\big)\le C d(x,y)^{-n},\ea\label{mixed}\ee
where in the second inequality we used $|U_j(z,y)|\le |\nabla_y G_2(z,y)|\le C d(z,y)^{1-n}$ (for example from \eqref{mf6a}, \eqref{mf6b}).
Combining \eqref{gradlip} and \eqref{mixed} gives \eqref{regK} for $\alpha=1.$\smallskip

To prove the sharpness, consider a smoothing of the following function
\be v_\e(y)=\bc \log\frac{1}{d(x_0,y)}\qquad&\text{for}\ \e\le d(x_0,y)\le 1\cr \log\frac{1}{\e}&\text{for}\ d(x_0,y)\le \e,\ec\nt\ee
where $x_0$ is a fixed point in $M$.
Then  by [F, Proposition 3.6], we have
\be \|\nabla^{\a}v_\e\|_{n/\a}^{\frac{n}{n-\a}}=\frac{\gam_{n,\a}}{n}(\log\frac{1}{\e})^{\frac{\a}{n-\a}}+O(1).\ee
Also, note that $ \|v_\e\|_{n/\a}\le C$. Take $u_\e={v_\e}({\|\nabla^{\a}v_\e\|_{n/\a}})^{-1}$.
We have \be \|\nabla^\a u_\e\|_{n/\a}\le 1,\qquad \|u_\e\|_{n/\a}^{n/\a}\le C(\log\frac{1}{\e})^{-1},\nt\ee
and \be |u_\e|^{\frac{n}{n-\a}}\geq \gam_{n,\a}^{-1}\log\frac{1}{\e}\qquad\text{for}\ d(x_0,y)\le\e.\nt\ee
For the sharpness of the exponential constant, we take $\theta >1$ and estimate
\be \ba &\int_{M}\frac{\exp_{\lceil\frac{n-\a}{\a}-1\rceil}\big[\theta \gam_{n,\a}|u_\e|^{\frac{n}{n-\a}}\big]}{1+|u_\e|^{\frac{n}{n-\a}}}d\mu(y)\geq \int_{d(x_0,y)\le\e}\frac{\exp\big[\theta \gamma_{n,\a}|u_\e|^{\frac{n}{n-\a}}\big]}{1+|u_\e|^{\frac{n}{n-\a}}}dy\cr &\geq C\e^n\frac{\exp\bigg[\theta\log\dfrac{1}{\e^n}+C\bigg]}{1+C\log\dfrac{1}{\e}}= \frac{C\e^{(1-\theta)n}}{1+C\log\dfrac{1}{\e}}\to\infty \ea\nt\ee
as $\e\to0^+.$\hfill\QEDopen
\par
%For the sharpness of the power of the denominator, we take $\theta <1$ and get 
%\be \ba &\|u_\e\|_{n/\a}^{-n/\a}\int_{M}\frac{\exp_{\lceil\frac{\sigma n}{n-\a}-\si-1\rceil}\big[\ga(n,\a)|u_\e|^{\frac{n}{n-\a}}\big]}{1+|u_\e|^{\frac{\theta \si n}{n-\a}}}d\mu(y)\geq \|u_\e\|_{n/\a}^{-n/\a}\int_{d(x_0,y)\le\e}\frac{\exp\big[ \ga(n,\a)|u_\e|^{\frac{n}{n-\a}}\big]}{1+|u_\e|^{\frac{\theta\si n}{n-\a}}}d\mu(y)\cr&
%\geq C(\log\dfrac{1}{\e})^{\si}\e^n\frac{\exp\bigg[\log\dfrac{1}{\e^n}+C\bigg]}{1+C(\log\dfrac{1}{\e})^{\theta\si}}\geq C(\log\dfrac{1}{\e})^{\si(1-\theta)}\to\infty \ea\nt\ee
%as $\e\to0^+$.
%\hfill\qed

\begin{rk}\label{green} Estimate \eqref{mf6b} can be refined as
\be C_{b,n}' r^{\a/2-1}e^{-b(n-1)d(x,y)}\le |K_\a(x,y)|\le C_{a,n} r^{\a/2-1}e^{-a(n-1)d(x,y)},\qquad d(x,y)\ge 1\label {refined}\ee
where $C_{a,n},C_{b,n}'$ are two positive constants depending on $b,n$. This bound can be obtained  directly from the heat kernel comparison theorem due to Debiard-Gaveau-Mazet [DGM, Th\'eor\`eme~1])
\be p^{(-b^2)}(t,x,y)\le p(t,x,y)\le p^{(-a^2)}(t,x,y),\qquad t>0,\; x,y\in M\ee
where $p(t,x,y)$ is the heat kernel for the Laplace-Beltrami operator on $M$ and $p^{(\lambda)}(t,x,y)$ the one on the space form of constant sectional curvature $\lambda.$

We note that when $\a=2$ the above refined bound is not enough to prove \eqref{GG2}, unless $a$ and $b$ are close enough. In fact \eqref{GG2} suggests a bound of type 
\be |G_2(x,y)|\le CV_x\big(d(x,y)\big)^{-1},\ee
which is certainly true when $K$ is constant, but  not known in general.

Also, it is possible to refine \eqref{mf7b} as follows
\be |K_\a|^*(x,t)\le C{(\log t)^{{\a\over2}-1}\over t},\qquad t>1\ee
with the same strategy as in [BS, Thm. 4.2], using induction and O'Neil's Lemma.
\end{rk}

\bigskip\medskip

%###########################################################
%\input{Appendix}

\newpage
\appendix
\ni{\LARGE \bf Appendices}\bigskip
\section{Proof of the equivalency of \eqref{k1},\eqref{k2} and  \eqref{k1p},\eqref{k2p}}\label{app-equivalency}

 \bigskip
It is enough to prove that assuming \eqref{k3}, there exists $C$ such that 
 \be\bigg|\mathop\int\limits_{r_1\le d(x,y)\le r_2} |k(x,y)|^\b d\mu(y)-\mathop\int_{V_x(r_1)}^{V_x(r_2)} (k^*(x,t))^\b dt\bigg|\le C \quad \forall x\in M ,\label{k11}\ee  
where $0<V_x(r_1)<V_x(r_2).$ 

Fix  $x\in M$. By the same argument used to find $E_\tau, F_\tau$ in \eqref{d3},\eqref{d2}, we can find two sets $H_i,\ i=1,2$ such that  
\be\bc &\{y\ :\ |k(x,y)|>k^*(x,V_x(r_i))\}\subseteq H_i\subseteq \{y\ :\ |k(x,y)|\geq k^*(x,V_x(r_i))\}\cr & \mu(H_i)=V_x(r_i),\ec
\nt\ee
Note that $H_1\subseteq H_2$ by construction and the fact that $k^*(x,V_x(r_1))\geq k^*(x,V_x(r_2))$. 

 We can write
\be\ba& \int_{V_x(r_1)}^{V_x(r_2)} (k^*(x,t))^\b dt=\int_{H_2\setminus H_1}|k(x,y)|^{\b}d\mu(y)\le\mathop\int\limits_{r_1\le d(x,y)\le r_2} |k(x,y)|^\b d\mu(y) \cr
&\hskip3em +\mathop\int\limits_{(H_2\setminus H_1 )\cap B(x,r_1)} |k(x,y)|^\b d\mu(y)+\mathop\int\limits_{(H_2\setminus H_1 )\cap B(x,r_2)^c} |k(x,y)|^\b d\mu(y)\cr&
%+\mu(H_2\setminus H_1)\supess_{y:\ d(x,y)\geq r_2}|k(x,y)|^\b\cr& 
 \le\mathop\int\limits_{r_1\le d(x,y)\le r_2} |k(x,y)|^\b d\mu(y)+V(r_1)\big(k^*(x,V_x(r_1)\big)^\b+CV_x(r_2)\supess_{ d(x,y)\geq r_2}\frac{1}{V_x\big(d(x,y)\big)}\cr
%&\le CV_x(r_1)(V_x(r_1))^{-1}+\mathop\int\limits_{r_1\le d(x,y)\le r_2} |k(x,y)|^\b d\mu(y)+CV_x(r_2)(V_x(r_2))^{-1}\cr
&\le \mathop\int\limits_{r_1\le d(x,y)\le r_2} |k(x,y)|^\b d\mu(y)+C,
\ea\label{rk31}\ee
where for the last inequality we used  \eqref{k3} and \eqref{kk12} (which follows from \eqref{k3}).

\begin{rk} \label{k1prime} If $A_\infty=0$ then \eqref{k1} implies \eqref{k1p} provided that the first inequality in \eqref{k3} holds for $V_x(d(x,y))\le 1$, and that \eqref{kk12} holds. Indeed, we argue as in \eqref{rk31} but we estimate

\be\ba  \mathop\int\limits_{(H_2\setminus H_1 )\cap B(x,r_2)^c} |k(x,y)|^\b d\mu(y)\le \mu(H_2\setminus H_1)&\supess_{r_2\le d(x,y)\le r_x(1)}|k(x,y)|^\b\cr&+\mathop\int\limits_{d(x,y)\ge r_x(1)} |k(x,y)|^\b d\mu(y)\le C\ea\ee

Note that  if $k$ is symmetric, then we do  not even  need to assume \eqref{kk12} since it follows from the above hypothesis. Indeed, critical integrability implies, for all $s>0$,
\be\mu\{y:\;|k(x,y)|>s,\, d(x,y)\ge r_x(1)\}\le Cs^{-\b}\ee
and the validity of \eqref{k3} for volumes less then or equal  1 implies
\be\mu\{y:\;|k(x,y)|>s,\; d(x,y)<r_x(1)\}\le Cs^{-\b}\ee
from which it follows $k_1^*(t)=k_2^*(t)\le Ct^{-1/\b}.$

\end{rk}

Next let  $\tilde k(x,y)=k(x,y)\chi_{\{r_1\le d(x,y)\le r_{2}\}}$. Then clearly $\tilde k^*(x,t)\le k^*(x,t)$, and by \eqref{k3} we have  $|\tilde k(x,y)|^\b\le CV_x(r_1)^{-1}$, hence
 \be \tilde k^*(x,t)^\b\le CV_x(r_1)^{-1}. \nt\ee
 So we have that
 \be \ba \mathop\int\limits_{r_1\le d(x,y)\le r_2} |k(x,y)|^\b d\mu(y)&\le \int_{0}^{V_x(r_2)}| \tilde k^*(x,t)|^\b dt\cr
 &=\int_{V_x(r_1)}^{V_x(r_2)}| \tilde k^*(x,t)|^\b dt+\int_0^{V_x(r_1)}| \tilde k^*(x,t)|^\b dt\cr
 &\le\int_{V_x(r_1)}^{V_x(r_2)}| k^*(x,t)|^\b dt+C\int_0^{V_x(r_1)}(V_x(r_1))^{-1}dt\cr
 &\le\int_{V_x(r_1)}^{V_x(r_2)}| k^*(x,t)|^\b dt+C.
 \ea\label{rk32}\ee
 which shows that \eqref{k1} and \eqref{k2} follows from \eqref{k1p} and \eqref{k2p}, assuming \eqref{k3}.
  \hfill\QEDopen
 
\begin{rk}\label{k1prime1} Under the assumption that the first inequality in  \eqref{k3} holds when $V_x(d(x,y))\le~1$ the argument in \eqref{rk32} is still valid, to show how \eqref{k1} follows from \eqref{k1p}. If, additionally, \eqref{k2p} holds with  $A_\infty=0$ and $k$
is bounded for large volumes, then \eqref{k2} holds, with $A_\infty=0$. Indeed, if $|k(x,y)|\le C_0$ for $V_x(d(x,y))\ge 1$ then taking $r_1=r_x(1)$ we get  $\tilde k^*(x,t)^\b\le C_0$ and as in \eqref{rk32} we deduce
\be \ba \mathop\int\limits_{ d(x,y)\ge r_x(1) } |k(x,y)|^\b d\mu(y)
\le  \int_1^{\infty}  (k^*(x,t))^\b dt+\int_0^1 C_0^\b du\le C.\ea\ee

\end{rk}
%\be\bigg|\mathop\int\limits_{r_1\le d(x,y)\le r_2}|k(x,y)|^\b d\mu(y)-  A\log{\frac{V_x(r_2)}{V_x(r_1)}}\,\bigg|\le B,  \qquad  0<r_1<r_2,\;\;\forall x\in M\label{A1p}\tag{A1'}\ee
 
\bigskip\section{Proof of \eqref{vols}, \eqref{vols1}}\label{Vols}
First note that  $V_x(R_{j+1})=\tau e^{\b'(j+1)}=e^{\b'}V_x(R_j).$ 
If $j\le N-1$ then  $r_j\le R_j$, hence $V_x(r_j)\le V_x(R_j)$.
From \eqref{RN} we have
\be V_x(R_N)\le V_x(r_N)\le e^{\b'}V_x(r_N)\nt\ee and
\be \ba (e^{\b'}-1)V_x(R_{N-1})\le V_x(R_N)-V_x(R_{N-1})\le V_x(r_N)-V_x(r_{N-1})&\le V_x(r_N)\cr&\le e^{\b'}V_x(R_N)=e^{2\b'}V_x(R_{N-1}).\ea\nt\ee

\ni Now let us  prove that for any $j=0,1,...,N-1$ there is $C>0$ such that $V_x(r_j)\ge CV_x(R_j)$. Clearly this is true if $r_j=R_j$. If $r_j<R_j$ then, since $r_0=R_0$,  there is $i$ with $0\le i<j$, $r_i=R_i$ and $r_k<R_k$ for $i<k\le j$.
  
  If $i<j_1<j$ then from Remark \ref{rk4} following the definition of $r_j$ in \eqref{defr}, and using conditions  \eqref{k1}, \eqref{k2} we get 
  
 \be (j_1-i)\b'A_0=\int_{r_{i}\le d(x,y)\le r_{j_1}} |k(x,y)|^\b d\mu(y)\le A_0\log{\frac{V_x(r_{j_1})}{V_x(r_i)}}+B,\label{ctau2}\ee
 and
  \be  (j-j_1)\b'A_\infty=\int_{r_{j_1}\le d(x,y)\le r_{j}} |k(x,y)|^\b d\mu(y)\le A_\infty\log{\frac{V_x(r_{j})}{V_x(r_{j_1})}}+B,\label{ctau3}\ee
therefore  
\be V_x(r_j)\ge e^{-\frac{B}{A_0}-\frac{B}{A_\infty}} e^{\b'(j-i)} V_x(r_i)= e^{-\frac{B}{A_0}-\frac{B}{A_\infty}} e^{\b'(j-i)} V_x(R_i)= e^{-\frac{B}{A_0}-\frac{B}{A_\infty}} V_x(R_j).\ee

If $j_1\le i$ then the proof is identical except we only use \eqref{ctau3} with $i$ in place of $j_1$, and if $j_1\ge j$ we only use \eqref{ctau2} with $j$ in place of $j_1$.
 
\medskip
We need to show now that there is $C>0$ such that $V_x(r_{j+1})-V_x(r_j)\ge C V_x(R_j)$ for $j=0,1,2,...,N-2.$ If $r_{j+1}=R_{j+1}$ then
 
  \be V_x(r_{j+1})-V_x(r_j)=V_{x}(R_{j+1})-V_x(r_j)\geq V_{x}(R_{j+1})-V_x(R_j)=(e^{\b'}-1)V_x(R_j). \nt\ee

 If $r_{j+1}<R_{j+1}$, then by definition of $r_j$ and Remark \ref{rk4} we have that
 \be \int_{r^{}_{j}\le d(x,y)\le r^{}_{j+1}} |k(x,y)|^\b d\mu(y)\geq\b'\min\{A_0,A_\infty\}.\nt\ee
 On the other hand, we can use condition \eqref{k3} to get
 \be\ba \int_{r^{}_{j}\le d(x,y)\le r^{}_{j+1}} |k(x,y)|^\b d\mu(y)\le \frac{B}{ V_{x}(r_j)}\big(V_x(r_{j+1})-V_x(r_j)\big)\ea\nt\ee
 and therefore, 
 \be\ba V_x(r_{j+1})-V_x(r_j)&\geq \frac{\b'\min\{A_0,A_\infty\}}{B}V_{x}(r_j)\geq \frac{\b'\min\{A_0,A_\infty\}}{B}C V_x(R_j).\ea\label{ctau5}\ee
 \ni \hfill\QEDopen
 
\section{Proof of Lemma \ref{lo}}\label{app-lo}

\vskip 0.2in  Suppose first that $\lambda_k=1$ for all $k.$ We can assume WLOG that the sequences $a=\{a_k\}$ in the definition of $\mu(h)$ in \eqref{muh1} are nonnegative.

It is enough to show the result for $h\ge1$. Indeed, if  $0<h_0\le h\le1$ and $\|a\|_1=h$, then we have $\|a\|_{\b'}\le \|a\|_1=h$, and $\mu(h)=h^{p}\mu(1)$ so \eqref{lo2} follows in this case. 

\smallskip
By dilation it is  easy to see that $\mu(h)$ increases as $h>0$ increases, so it is enough to show that for any integer $N\geq 1$, the following holds 
\be C\frac{e^{\b'N}}{N^{p/\b' }}\leq \mu(N^{1/\b})\leq C' \frac{e^{\b'N}}{N^{p/\b' }}.\label{ap1}\ee

%Let $\|a\|_{(e)}=\big(\sum_{k=0}^\infty |a_k|^{\si\b'}e^{\b'k}\big)^{1/\b'}$.
Consider the sequence $a_k=N^{-1/\b'}$ for $k=0,...,N-1$, and $a_k=0$ otherwise. Then $ \|a\|_1=N^{1/\b}, \|a\|_{\b'}=1$ and 
\be\mu(N^{1/\b})\le\sum_{k=0}^\infty a_k^{p}e^{\b'k}=N^{-{p/\b'}}\frac{e^{N\b'}-e^{\b'}}{e^{\b'}-1}\le C\frac{e^{N\b'}}{N^{p/\b'}}\nt\ee
So we are left to prove that  for any nonnegative  sequence such that 

\be \|a\|_1=N^{1/\b} \qquad\ \|a\|_{\b'}\le 1\qquad\ \label{ap2}\ee  we have 
\be  \|a_k^{p/\b'} e^k\|_{\b'}^{\b'}=\sum_{k=0}^\infty a_k^{p} e^{\b'k}\ge C \frac{e^{\b'N}}{N^{{p/\b'} }}.\label{ap0}\ee 

If 
\be a_{N-1}\ge \frac{1}{2\beta}N^{-1/\b'}\ee
then \eqref{ap0} follows using
\be  \|a_k^{p/\b'} e^k\|_{\b'}^{\b'}\geq a_{N-1}^{p}e^{\b'(N-1)}\ge({2\b})^{-p}e^{-\b'}\,\frac{e^{\b'N}}{N^{p/\b'}}\ee

Assume then that
\be a_{N-1}< \frac{1}{2\beta}N^{-1/\b'}\label{ap00}\ee
and let\def\lla{\Lambda} $\lla$  be defined so that
\be \lla^{p}= \|a_k^{p/\b'} e^k\|_{\b'}^{\b'}\frac{N^{p/\b'}}{ e^{\b' N}}.\ee

\ni Let us also define $a_{-1}=0$ so that the proof below works out also for $N=1$. 

On one hand we have, from H\"{o}lder's inequality and \eqref{ap2}
\be\bigg(\sum_{j=0}^{N-2} a_j\bigg)^\b\le (N-1)\bigg(\sum_{j=0}^{N-2} a_j^{\b'}\bigg)^{\b/\b'}\le N-1\label{ap3}\ee

On the other hand, 
\be\sum_{j=0}^{N-2}a_j=N^{1/\b}-a_{N-1}-\sum_{j=N}^\infty a_j\label{ap4}\ee

since $a_j\le \|a_k^{p/\b'} e^k\|_{\b'}^{\b'/p}e^{-j\b'/p}$ we have
\be\sum_{j=N}^\infty a_j\le \sum_{j=N}^\infty \|a_k^{p/\b'} e^k\|_{\b'}^{\b'/p}e^{-j\b'/p}=\|a_k^{p/\b'} e^k\|_{\b'}^{\b'/p}\frac{e^{-N\b'/p}}{1-e^{-\b'/p}}=C_0\lla N^{-1/\b'}\ee
so that from \eqref{ap4} and \eqref{ap00}  we get
\be \sum_{j=0}^{N-2}a_j\ge N^{1/\b}-\frac1{2\b}N^{-1/\b'}-C_0\lla N^{-1/\b'}=N^{1/\b}\bigg(1-\Big(\frac1{2\b}+C_0\lla\Big) N^{-1}\bigg).\ee
If $\frac1{2\b}+C_0\lla\ge1$ then $\lla\ge C_0^{-1}(1-\frac1{2\b})>0$ and we are done. If instead $\frac1{2\b}+C_0\lla<1$ then using 
the inequality $(1-x)^\beta> 1-\beta x$ for $0<x<1$, we get
\be \bigg(\sum_{j=0}^{N-2}a_j\bigg)^\b> N-\Big(\frac1 2+\b C_0\lla\Big)\ee
which together with \eqref{ap3} yields $N-1>N-\frac12-\b C_0\lla$, i.e.  
\be \lla>\frac1 {2\b C_0}.\ee
\bigskip

To settle the more general case, let  $\lambda=\{\lambda_k\}$ be so that  $\lambda_k\in(0,1]$ with $\lambda_k<1$, for finitely many $k$, and let 
\be \mu'(h)=\inf \Big\{\aa a_k^{p}e^{\b'k}: \ \|a\|_1=h, \|\lambda a\|_{\b'}\le 1\Big\}. \label{muh2}\ee
 We will prove \eqref{lo2} by showing that for some $C_3>0$
\be C_3 e^{-L\b'}\mu(h)\le \mu'(h)\le \mu(h).\label{muh3}\ee
It is clear that $\mu'(h)\le \mu(h)$, since any sequence satisfying $\|a\|_{\b'} \leq 1$ also satisfies $\|\la a\|_{\b'}\le1$. \par
It is then enough to show that for any sequence $a$ such that
$ \|a\|_1=h,$  $\|\la a\|_{\b'}\le1$ there is a  nonnegative sequence $b$ with
  $ \|b\|_1=h,$ $\|b\|_{\b'}\le1$
and 
\be C_3 e^{-L\b'}\aa b_k^{p}e^{\b'k}\le \aa a_k^{p}e^{\b'k}.\label{muh6}\ee
 Define
\be J_k=\bc \lceil\la_k^{-\b}\rceil \ &\text{if}\ k\ge0 \cr 0 & \text {if} \ k=-1,\ec\nt\ee
and
\be b_j=\bc a_0/J_0,\  & 0\le j\le J_0-1\cr
a_1/J_1, &J_0\le j\le J_0+J_1-1\cr
......\cr
a_k/J_k, &\sum_{\ell=0}^{k-1}J_\ell\le j\le \sum_{\ell=0}^k J_\ell -1\cr
......\ec\nt\ee

Then we have 
\be \|b\|_1=\frac{a_0}{J_0}J_0+...+\frac{a_k}{J_k}J_k+...=\|a\|_1=h\nt\ee
\be \ba\|b\|_{\b'}^{\b'}=\bigg(\frac{a_0}{J_0}\bigg)^{\b'}J_0+...+\bigg(\frac{a_k}{J_k}\bigg)^{\b'}J_k+...=\sum_{k=0}^\infty J_k^{1-\b'}a_k^{\b'}\le \sum_{k=0}^\infty \la_k^{\b'}a_k^{\b'}\le 1,\ea\nt\ee
and 
\be\ba \sum_{j=0}^\infty & b_j^{p}e^{\b'j}=\bigg(\frac{a_0}{J_0}\bigg)^{p}(1+...+e^{\b'(J_0-1)})+...\cr&\hskip4em +\bigg(\frac{a_k}{J_k}\bigg)^{p}(e^{\b'(J_0+...+J_{k-1})}+...+e^{\b'(J_0+...+J_{i}-1)})+...\cr
&\le \sum_{k=0}^\infty a_k^{p}e^{\b'\sum_{\ell=0}^{k-1}J_\ell}\sum_{j=0}^{J_k-1} e^{\b'k}\le\sum_{k=0}^\infty a_k^{p}e^{\b'\sum_{\ell=0}^{k-1}J_\ell}\sum_{j=0}^{J_k-1} e^{\b'k}\le\frac{1}{e^{\b'}-1}\sum_{k=0}^\infty a_k^{p} e^{\b'\sum_{\ell=0}^{k}J_\ell}\cr&=
\frac{e^{\b'}}{e^{\b'}-1}\sum_{k=0}^\infty a_k^{p}e^{\b'k} e^{\b'\sum_{\ell=0}^{k}(J_\ell-1)}\le\frac{e^{\b'}e^{\b'L}}{e^{\b'}-1}\sum_{k=0}^\infty a_k^{p}e^{\b'k},\ea\nt\ee
which gives \eqref{muh6}.

\hfill\QEDopen

\bigskip\section{Proof of Lemma \ref{lb}}\label{app-lb}

\vskip 0.2in

\ni{\bf Proof of \eqref{lb1}:} 
Recall \eqref{vols}, \eqref{vols1}, \eqref{vols2} and that $\mu(E_\tau)\le\tau.$ Using \eqref{kk12} we get, for fixed $y\in M$, 
\be \ba \int_{(D_{j+1}\setminus D_j)\setminus E_\tau}|k(\xi,y)|d\mu(\xi)&\le \int_0^{\mu(D_{j+1}\setminus D_j)}k_2^*(t,y)dt\le C\int_0^{\mu(D_{j+1}\setminus D_j)} t^{-1/\b}dt
&=C\mu(D_{j+1}\setminus D_j)^{1-1/\b}.\ea\label{lem1}\ee
Hence, we have
\be \ba \avgdj& |Tf_\tau|d\mu
  \leq \frac{C}{\mu(D_{j+1}\setminus D_j)}\int _{M}|f_\tau(y)|\int_{(D_{j+1}\setminus D_j)\setminus E_\tau}|k(\xi,y)|d\mu(\xi)d\mu(y)\cr&\leq  C\mu(D_{j+1}\setminus D_j) ^{-1/\b}\int _{M}|f_\tau(y)|d\mu(y)\le CV_x(R_j)^{-1/\b}\int _{M}|f|\chi_{F_\tau}^{}d\mu\cr&\le CV_x(R_j)^{-1/\b}\mu(F_\tau)^{1/\b}\|f\|_{\b'}=Ce^{-j\b'/\b}\|f\|_{\b'}\le C_3e^{-(\b'-1)j}. \ea\label{lem2}\ee
  
  \bigskip
 \ni{\bf Proof of \eqref{lc1}:} We will show that in fact for every $\xi\in D_{j+1}\setminus D_j$ we have
 \be | {T}S_{j+2}f_\tau'(\xi)-{T}S_{j+1}f_\tau'(x)|\le C\b_{j+1},\label{ineqD}\ee
 which obviously implies \eqref{lc1}. First write 
   \be \ba &| {T}S_{j+2}f_\tau'(\xi)- {T}S_{j+1}f_\tau'(x)|\le|{T}S_{j+2}f_\tau'(\xi)-{T}S_{j+2}f_\tau'(x)|+|{T}\big(S_{j+1}-S_{j+2})f_\tau'(x)|.
 \ea\label{lc2}\ee
 From \eqref{la3}
 \be \ba |{T}\big(S_{j+1}-S_{j+2})f_\tau'(x)|\le C\alpha_{j+1}\le C\b_{j+1}.
 \ea\label{lem4}\ee

\smallskip
Now let us  estimate the first term in \eqref{lc2}:
 \be \ba  &|{T}S_{j+2}f_\tau'(\xi)-{T}S_{j+2} f_\tau'(x)|
 =\bigg|\int_{D^c_{j+2}\setminus D_0'}(K(\xi,y)-K(x,y)) f_\tau'(y)d\mu(y)\bigg|\cr
&\leq \bigg(\int_{D^c_{j+2}\setminus D_0'}|k(\xi,y)-k(x,y)|^\b d\mu(y)\bigg)^{1/\b}\bigg(\int_{D^c_{j+2}\setminus D_0'}| f_\tau'(y)|^{\b'}d\mu(y)\bigg)^{1/\b'}
\cr& \leq \bigg(\int_{d(x,y)\ge r_{m+j+1}}|k(\xi,y)-k(x,y)|^\b d\mu(y)\bigg)^{1/\b}\bigg(\int_{D^c_{j+1}\setminus D_0'}| f_\tau'(y)|^{\b'}d\mu(y)\bigg)^{1/\b'}\le C\b_{j+1}.\ea\label{lc4}\ee
where the last inequality follows from \eqref{k4} since $x,\xi\in B(x,r_{j+1})$ and there is $\delta>0$ (depending only from $k,\b$) such that
\be V_x(r_{j+2})\ge (1+\delta)V_x(r_{j+1}),\ee
which can be deduced  from \eqref{vols}, \eqref{vols1}.

\bigskip
\ni {\bf Proof of \eqref{ld1}:} We use \eqref{lem1} to get, for any $k$
\be \ba \bigg|\avgdj T(S_k&-S_{k+1})f_\tau'(\xi)d\mu(\xi)\bigg|=\bigg|\avgdj \int_{(D_{k+1}\setminus D_k)\setminus D_0'} k(\xi,y) f_\tau'(y)d\mu(y)d\mu(\xi)\bigg|\cr
&\le C\mu(D_{j+1}\setminus D_j)^{-1/\b}\int_{(D_{k+1}\setminus D_k)\setminus D_0'}|f_\tau'(y)|d\mu(y)\cr& 
\le CV_x(R_j)^{-1/\b} V_x(R_k)^{1/\b}\a_k=Ce^{-\b'(j-k)}\a_k.\ea
\label{lem6}\ee

Therefore,
\be \ba \bigg|&\avgdj{T}(S_0-S_{j+2})f_\tau'(\xi)d\mu(\xi)\bigg|=\bigg|\avgdj \sum_{k=0}^{j+1}T(S_k-S_{k+1})f_\tau'(\xi)d\mu(\xi)\bigg|\cr
&\le \sum_{k=0}^{j+1}\bigg|\avgdj T(S_k-S_{k+1})f_\tau'(\xi)d\mu(\xi)\bigg|\le  C\sum_{k=0}^{j+1} e^{-\b'(j-k)}\a_k\le C\ab_{j+1}.
\ea
\label{lem8}\ee

\bigskip
{\bf Proof of \eqref{ld2} and \eqref{ldd1}:} Both \eqref{ld2} and \eqref{ldd1} can be derived exactly the same way as we did in proof of \eqref{lb1}, using  $f_\tau'\chi^{}_{D_0\setminus D_0'}$ and $f_\tau'\chi^{}_{D_0'}$ in place of $f_\tau$. \hfill\QEDopen

\section{Some results on $(Tf)^\circ$ for the Riesz potential on  $\R^n$}\label{tf}
%udl{(1) An example of $Tf^\circ(\tau)=u(r)$:}\par

\bigskip

First let us remark that for an operator of type 
$$Tf(x)=\int_{\R^n} K(|x-y|)f(y)dy$$
if both $K$ and $f$ are nonnegative and  radially decreasing, then $Tf$ is also nonnegative and radially decreasing, hence $E_\tau=F_\tau=B(0,r)$, up to a set of 0 measure, where $|B(0,r)|=\tau.$
Also, we have that $(Tf)^\circ(\tau)=T\big(f\chi_{B(0,r)^c}\big)(0)$, i.e. the $\supess$ in the definition of $(Tf)^\circ$ is a max attained at $x=0$. This is because  for $|x|\le r$
\be \ba& T(f\chi_{B(0,r)^c  \cap B(x,r)^c})(x)=\int_{\R^n}\big( K\chi_{B(0,r)^c}\big)(x-y) \big(f\chi_{B(0,r)^c}\big)(y)dy\cr&\hskip2em\le \int_0^\infty \big( K\chi_{B(0,r)^c}\big)^*(t)\big(f\chi_{B(0,r)^c}\big)^*(t)\cr&\hskip2em =\int_\tau^\infty K^*(t) f^*(t)dt=\int_{B(0,r)^c}K(y)f(y)dy =T(f\chi_{B(0,r)^c})(0).\ea\ee

Let $I_\alpha f$ be the Riesz potential on $\R^n$ i.e.
$$I_\a f(x)=\int_{\R^n} |x-y|^{\a-n} f(y)dy$$
and let, for a vector-valued function $f$,
$$\wt I_1 f(x)=\int_{\R^n} \nabla_y |x-y|^{2-n} f(y)=(n-2)\int_{\R^n}|x-y|^{-n}(x-y)\cdot f(y)dy.$$

\begin{prop}\label{circstar} Let $f:\R^n\to\R$ be nonnegative, radially decreasing, and with compact support. For any $\tau>0 $ the following hold:
\smallskip \item{a)}   If $0<\a\le 2$, then $(I_\a f)^\circ(\tau)\le (I_\a f)^* (\tau)$ with strict inequality if $f>0$ near $0$.\smallskip
\smallskip\item{b)} If $\a>2$, $f$ bounded and positive near 0, then $(I_\a f)^\circ(\tau)>(I_\a f)^* (\tau)$ for all $\tau$ small enough.

\smallskip\ni
 If $f$ is vector-valued and $|f|$ is radially decreasing with compact support, then for every $\tau>0$ we have $(\wt I_1  f)^\circ (\tau)=(\wt I_1  f)^*(\tau).$ \end{prop}

\bigskip
\ni{\bf Proof.}  We have
\be I_\a f(x)=\int_0^{\infty}f(\rho)\rho^{\a-1}d\rho\int_{\S^{n-1}}|\rho^{-1}x-y^*|^{\a-n}dy^*\label {X1}\ee
and due to invariance under rotation we can assume $x=r e_1$. We will write $f(r)=f(re_1)$ and $I_\a f (r)=I_\a f(re_1)$.

 Note that we have
\be\avint_{\S^{n-1}}|z-y^*|^{2-n}dy^*=\bc |z|^{2-n} & \text {if}\, |z|>1\cr 1 & \text {if}\, |z|<1 \ec\label{form1}\ee
and for $0<\alpha<2$
\be\avint_{\S^{n-1}}|z-y^*|^{\a-n}dy^*\ge\bc |z|^{2-n} & \text {if}\, |z|>1\cr 1 & \text {if}\, |z|<1 \ec\label{form2}\ee
while  the  inequality in \eqref{form2} is reversed if $\alpha>2$.
The above are due to the fact that for fixed $z$ with $|z|>1$ the function $|z-y|^{\a-n}$ in the unit ball is harmonic for $\alpha=2$, subharmonic for $0<\alpha<2$, and superharmonic for $\alpha>2$. If $|z|<1$ then $|z-y^*|=\big|z^*-|z|y^*\big|$ and the function $\big|z^*-|z|y\big|^{\a-n}$ defined for $y$ in the unit ball, is harmonic for $\alpha=2$, subharmonic for $0<\alpha<2$ and superharmonic for $\alpha>2$, and continuous up to the boundary.

When $0<\alpha\le 2$ we then have 
\be\ba(I_\alpha f)^\circ(\tau)&=I_\alpha(f\chi_{B(0,r)^c})(0)=\omega_{n-1}\int_r^\infty\!\! f(\rho)\rho^{\a-1}d\rho \le \int_r^\infty\!\! f(\rho)\rho^{\a-1}d\rho\int_{\S^{n-1}}|r\rho^{-1}e_1-y^*|^{\a-n}dy^*\cr&=\int_{|y|\ge r} |re_1-y^*|^{\a-n}f(y)dy\le I_\a f(r)=(I_\a f)^*(\tau)\ea\ee

with strict inequality in the last inequality if $f>0$ around $0$.

\medskip

Let $\a>2$. Assume further that $f$ is bounded. First we estimate
\be\ba   &\big|{\rho^{-1}}{re_1}-y^*\big|^{\a-n}=\big(1-2{\rho^{-1}}{r} y_1^*+{\rho^{-2}}{r^2}\big)^{(\a-n)/2}\cr
  &=1+\frac{n-\a}{2}\big({2\rho^{-1}}{r} y_1^*-{\rho^{-2}}{r^2}\big)+\frac{n-\a}{2}(n-\a+2){\rho^{-2}}{r^2}(y_1^*)^2+O({\rho^{-3}}{r^3}),
  \ea\label{tf12}\ee

  so that 
  \be \ba \int_{\S^{n-1}}\big|\rho^{-1}{r}e_1-y^*\big|^{\a-n}dy^*&=\omega_{n-1}-\frac{n-\a}{2}\rho^{-2}{r^2} \int_{\S^{n-1}}(1-(n-\a+2)(y_1^*)^2)dy^*+O(\rho^{-3}{r^3})\cr
  &=\omega_{n-1}-\omega_{n-1}(\rho^{-2}{r^2})\frac{(n-\a)(\a-2)}{2n}+O(\rho^{-3}{r^3}),\ea\nt\ee
  
%   \be \ba& \int_{\S^{n-1}}\big|\frac{r}{\rho}e_1-y^*\big|^{\a-n}dy^*=\omega_{n-1}-\frac{n-\a}{2}(\frac{r^2}{\rho^2}) \int_{\S^{n-1}}(1-(n-\a+2)(y_1^*)^2)dy^*+O(\frac{r^3}{\rho^3})\cr
%  &=\omega_{n-1}-\omega_{n-1}(\frac{r^2}{\rho^2})\frac{(n-\a)(\a-2)}{2n}+O(\frac{r^3}{\rho^3}),\ea\nt\ee
  where we used the identities $\int_{\S^{n-1}} y_1^*dy^*=0$ and $\int_{\S^{n-1}}(y_1^*)^2dy^*=1/n$.
  Therefore, 
\be\ba (I_\a f)^*(\tau)&=I_\a f(r)=\int_{\R^n}|re_1-y|^{\a-n}f(y)dy=\int_0^\infty\rho^{\a-1}f(\rho)\int_{\S^{n-1}}|\rho^{-1}{r{ e_1}}-y^*|^{\a-n}dy^*d\rho\cr
&=\int_0^\infty\rho^{\a-1}f(\rho)(\omega_{n-1}-\omega_{n-1}(\rho^{-2}{r^2})\frac{(n-\a)(\a-2)}{2n}+O(\rho^{-3}{r^3}))d\rho\cr
&=\omega_{n-1}\int_0^\infty\rho^{\a-1}f(\rho)-\omega_{n-1}r^2\frac{(n-\a)(\a-2)}{2n}\int_0^\infty\rho^{\a-3}f(\rho)d\rho+O(r^3)\cr
&=\omega_{n-1}\int_0^r\rho^{\a-1}f(\rho)d\rho+\omega_{n-1}\int_r^\infty\rho^{\a-1}f(\rho)d\rho-Cr^2+O(r^3)\cr
&=\omega_{n-1}\int_r^\infty\rho^{\a-1}f(\rho)d\rho+C'r^{\a}-Cr^2+O(r^3),
\ea\label{tf11}\ee

where $C,C'>0$ and the second last equality is by the fact that $f$ is bounded and $\a>2$.
Hence for $r$ small enough, we have
\be (I_\a f)^\circ(\tau)-(I_\a f)^*(\tau)=Cr^2-C'r^\a>0.\nt\ee

Finally, let $|f|$ be radially decreasing with compact support,  and write 
 \be\wt I_1 f(x)=(n-2)\int_0^\infty |f(\rho)|d\rho \int_{\S^{n-1}}|\rho^{-1}x-y^*|^{-n}(\rho^{-1}x\cdot y^*-1)dy^*.\label{X2}\ee
 
 If $H(z)$  denotes the LHS of \eqref{form1} it is  easy to see that 

 \be H(z)-\frac{|z|}{2-n}\,\frac{\partial H}{\partial |z|}=\avint_{\S^{n-1}}|z-y^*|^{-n}(1-z\cdot y^*)dy^*=\bc 0 & \text {if}\, |z|>1\cr 1 & \text {if}\, |z|<1 \ec\label{form3}\ee

which gives  
\be\wt I_1 f(x)=-(n-2)\omega_{n-1}\int_{|x|}^\infty|f(\rho)|d\rho.\label{form4}\ee
Hence both $|f|$ and $ |\wt I_1 f|$ are radially decreasing so that $E_\tau=F_\tau=B(0,r)$.

 Next, note that 
 
\be (x-\rho y^*)\cdot |f(y)| =(x\cdot y^*-\rho)|f|(\rho)\leq 0,\ \ \ \text{for}\  |x|\le r,\ |y|=\rho\geq r,\nt\ee
so that 
\be 0\le -\wt I_1(f\chi_{B(0,r)^c\cap B(x,r)^c})(x)\le -\wt I_1(f\chi_{B(0,r)^c})(x),\ \ \forall x\in B(0,r).\label{tf6}\ee
%\be 0\le T(f\chi_{B^c(0,r)\cap B^c(x,r)})(x)\le T(f\chi_{B^c(0,r)})(x)= T(f\chi_{B^c(0,r)})(0)=u(r),\ \ \forall x\in B(0,r).\label{tf6}\ee
Using \eqref{form3} and \eqref{form4} we get that for every $x\in B(0,r)$
\be\ba & \big|\wt I_1(f\chi_{B(0,r)^c})(x)\big|=(n-2)\int_{|y|\geq r}|x-y|^{-n}(y-x)\cdot f(y)dy\cr
&=(n-2)\int_r^\infty| f(\rho)|d\rho\int_{\S^{n-1}}|\rho^{-1}x-y^*|^{-n}(1-\rho^{-1} x\cdot y^*)dy^* =|\wt I_1 f(r)|=(\wt I_1 f)^*(\tau)\ea\label{tf3}\ee

Therefore combining \eqref{tf6} and \eqref{tf3} we have
\be (\wt I_1f)^\circ(\tau)=\wt I_1(f\chi_{B(0,r)^c})(0)=(\wt I_1f)^*(\tau),\qquad\forall \tau>0.\label{tf2}\ee
\hfill\QEDopen
%which implies that $Tf^\circ(\tau)=T(f\chi_{B^c(0,r)})(0)=u(r)$.

\bigskip

\medskip\ni{\bf{\large{Conflict of interest statement}}} 

\medskip

\ni On behalf of all authors, the corresponding author states that there is no conflict of interest.

\medskip

\medskip\ni{\bf{\large{Data availability statement }}}

\medskip

\ni Data sharing not applicable to this article as no datasets were generated or analysed during the current study.

%#############################################################
%\input{references}

%\renewcommand\bibname{References}

\bigskip\bigskip
\ni Carlo Morpurgo \hskip13.55em Liuyu Qin

\ni Department of Mathematics \hskip 8.25em Department of Mathematics and statistics

\ni University of Missouri \hskip 11.2em  Hunan University of Finance and Economics

\ni Columbia, Missouri 65211 \hskip 9.25em Changsha, Hunan

\ni USA \hskip 18.55em China

\smallskip\ni morpurgoc@umsystem.edu\hskip 9.35em \text{Liuyu\_Qin@outlook.com}


\begin{thebibliography}{ABCD}  % 100 is a random guess of the total number of %references 
 
%\addtolength{\leftmargin}{0.2in} % sets up alignment with the following line. 
%\setlength{\itemindent}{-0.2in} 
 \addcontentsline{toc}{section}{References}


\bibitem[A1]{A} Adams D.R., \emph{
A sharp inequality of J. Moser for higher order derivatives},
Ann. of Math. {\bf128} (1988), no. 2, 385--398. 

\bibitem[A2]{a2}Adams D.R., \emph{A trace inequality for generalized potentials}, Studia Math. {\bf 48} (1973) 99–105.


%\bibitem[AH]{ah} Adams D.R., Hedberg L.I., \emph{Function spaces and potential theory}, Springer, 1996.

\bibitem[AF] {af} Adams R.A., Fournier J.J.F., \emph {Sobolev spaces},  Second edition. Pure and Applied Mathematics (Amsterdam), 140. Elsevier/Academic Press, Amsterdam, 2003.

%\bibitem[ASM]{aam}Alves C. O. , Souto M. A. S. , Montenegro M. , \emph{Existence of a ground state solution for a nonlinear scalar field equation with critical growth}, Calc. Var. Partial Differential Equations {\bf 43} (2012), No. 3-4, 537-554.


\bibitem[AT]{at} Adachi S., Tanaka K., \emph{Trudinger type inequalities in $\R^n$ and their best exponents},
Proc. Amer. Math. Soc. {\bf 128} (2000), 2051-2057.

\bibitem[ARSW] {arsw} Arcozzi N., Rochberg R., Sawyer E.T., Wick B.D.,  \emph{Potential theory on trees, graphs and Ahlfors-regular metric spaces},  Potential Anal. {\bf  41} (2014), no. 2, 317-366. 

%\bibitem[BDR]{bdr} C. Benson, A. H. Dooley, G. Ratcliff, Fundamental solutions for powers of the Heisenberg sub-Laplacian, Illinois J. Math. {\bf 37} (1993), 455-476. 
\bibitem[ACDH]{acdh} Auscher P., Coulhon, T., Duong X.T., Hofmann S.,
\emph{Riesz transform on manifolds and heat kernel regularity}, Ann. Sci. \'Ecole Norm. Sup. (4) {\bf37} (2004), 911–957.
%\bibitem[Bec]{bec} Beckner W., \emph{Sharp Sobolev inequalities on the sphere and the Moser-Trudinger inequality}, Annals of Mathematics, {\bf 138} (1993), 213-242.
\bibitem[BDR]{bdr} Benson C., Dooley A.H., Ratcliff G., \emph{Fundamental solutions for powers of the Heisenberg
sub-Laplacian}, Illinois J. Math. {\bf37} (1993), 455-476.
\bibitem[BL]{bl} Berestycki H., Lions P.-L., \emph{Nonlinear scalar field equations, I-existence of a ground state}. Arch. Rat. Mech. Anal. {\bf 82} (1983), 313–346 .
\bibitem[BFM]{bfm} Branson T.P. , Fontana L., Morpurgo C., \emph{Moser–Trudinger and Beckner–Onofri’s inequalities on the CR sphere}, Annals of Mathematics {\bf 177} (2013), 1-52.
\bibitem[BS]{bs} Bertrand J., Sandeep K., \emph{Sharp Green's function estimations on Hadamard manifolds and Adams inequality}, Int. Math. Res. Not. IMRN, {\bf 6} (2021),  4729-4767. 
\bibitem[BeSh]{besh} Bennett C., Sharpley R., \emph{Interpolation of Operators}, Pure Appl. Math., vol. 129, Academic Press Inc., Boston, MA, 1988.
%\bibitem[BL]{bl} Berestycki, H., Lions, P.-L.\emph{Nonlinear scalar field equations, I-existence of a ground state}. Arch. Rat. Mech. Anal. {\bf 82}(1983), 313–346 .

%\bibitem[Bon]{bon} Bonfiglioli A., \emph{Taylor formula for homogenous groups and applications},Math. Z. (2009) 262:255–279.

%\bibitem[Cao]{cao} Cao D. M., \emph{Nontrivial solution of semilinear elliptic equation with critical exponent
%in $\R^2$}, Comm. Partial Differential Equations {\bf17} (1992), 407-435.

%\bibitem[Ch]{c} S.Y. CHANG, The Moser-Trudinger inequality and applications to some problems in conformal geometry, \emph{Nonlinear partial differential equations in differential geometry}, Comm. on Pure and  Applied Math, {\bf LVI}(2003), 1135-1150.
\bibitem[Ch]{ch} Christ M., \emph{A $T(b)$ theorem with remarks on analytic capacity and the Cauchy integral}, 
Colloq. Math. {\bf 60/61} (1990), no. 2, 601-628.
%\bibitem[CLLY]{clly}  Cohn W., Lam N., Lu G., Yang Y., \emph{The Moser–Trudinger inequality in unbounded domains of Heisenberg
%group and sub-elliptic equations}, Nonlinear Analysis {\bf 75} (2012) 4483–4495.
\bibitem[Ci1]{ci1} Cianchi A., \emph{Moser-Trudinger inequalities without boundary conditions and isoperimetric problems}, Indiana Univ. Math. J. {\bf54} (2005) 669-705. 
\bibitem[Ci2]{ci2} Cianchi A., \emph{Moser-Trudinger trace inequalities},  Adv. Math. {\bf 217} (2008),  2005-2044.

\bibitem[CL]{cl} Cohn W.S., Lu G., \emph{Best constants for Moser-Trudinger inequalities on the Heisenberg
group}, Indiana Univ. Math. J. {\bf 50} (2001), 1567-1591.

%\bibitem[CL2]{cl2} Cohn W.S., Lu G., \emph{Sharp Constants for Moser-Trudinger Inequalities on Spheres in Complex Space $\C^n$}, Communications on Pure and Applied Mathematics,  {\bf LVII} (2004), 1458–1493.
\bibitem[CRT]{crt} Cassani D., Ruf B., Tarsi C., \emph{Best constants in a borderline case of second-order Moser type inequalities}, Ann. Inst. H. Poincar\'e Anal. Non Lin\'eaire {\bf 27} (2010), 73-93.
\bibitem[CT]{ct}  Chang D.C., Tie J., Estimates for powers of sub-Laplacian on
the non-isotropic Heisenberg group, J. Geom. Anal. {\bf 10} (2000), 653-678.

\bibitem[DGM]{dgm} Debiard A., Gaveau B., Mazet E.,
\emph{Th\'eor\`emes de comparaison en g\'eom\'etrie riemannienne}.
Publ. Res. Inst. Math. Sci. {\bf12} (1976/77), 391-425.

%\bibitem[Fol1]{fol} Folland G.B., \emph{Real analysis: Modern techniques and their applications}, Wiley, 1999.

%\bibitem[Fol2]{fol2} Folland G.B., \emph{Subelliptic estimates and function spaces
%on nilpotent Lie groups},  Ark. Mat. {\bf13}, (1975), 161–207. 

\bibitem[Fol]{fol3} Folland G.B., \emph{
A fundamental solution for a subelliptic operator}, Bull. Amer. Math.
Soc. {\bf 79} (1973), 373-376.

\bibitem[F]{f} Fontana L., \emph{Sharp borderline Sobolev inequalities on compact Riemannian manifolds}, Comment. Math. Helv. {\bf 68} (1993), 415-454.

\bibitem[FM1]{fm1} Fontana L., Morpurgo C., \emph{Adams inequalities on measure spaces}, Adv. Math. {\bf 226} (2011), 5066-5119.

\bibitem[FM2]{fm2} Fontana L., Morpurgo C., \emph{Sharp exponential integrability for critical Riesz potentials and fractional Laplacians on $\R^n$}, Nonlinear Anal. {\bf 167}(2018), 85-122.

\bibitem[FM3]{fm3} Fontana L., Morpurgo C., \emph{Adams inequalities for Riesz subcritical potentials}, Nonlinear Anal. {\bf 192}(2020), 111662.
\bibitem[FM4]{fm4} Fontana L., Morpurgo C., \emph{Optimal limiting embeddings for $\Delta$-reduced Sobolev spaces in $L^1$}, Ann. Inst. H. Poincar\'e C Anal. Non Lin\'eaire {\bf31} (2014), 217-230.
\bibitem[FM5]{fm5} Fontana L., Morpurgo C., \emph{Sharp Moser-Trudinger inequalities for the Laplacian without boundary conditions}, J. Funct. Anal. {\bf262} (2012), 2231-2271.
\bibitem[H]{h} Hebey E., 
\emph{Sobolev spaces on Riemannian manifolds}, 
Lecture Notes in Mathematics, 1635. Springer-Verlag, Berlin, 1996.
%\bibitem[FS]{fs} Folland G.B., Stein E.M., \emph{Hardy spaces on homogeneous groups}, Princeton University Press, 1982.

\bibitem[HK]{hk} Haj\l{}asz P., Koskela, P., \emph{Sobolev met Poincar\'e},  Mem. Am. Math. Soc. {\bf 145} (2000), no. 688. 

\bibitem[IMN]{imn} Ibrahim S., Masmoudi N., Nakanishi K., \emph{ Trudinger-Moser inequality on the whole plane with the exact growth condition}, J. Eur. Math. Soc. (JEMS) {\bf 17} (2015), 819-835.

\bibitem[LL1]{ll1} Lam N., Lu G., \emph{Sharp Moser-Trudinger inequality on the Heisenberg group at the critical case and applications}, Adv. Math. {\bf 231} (2012), 3259-3287.

\bibitem[LL2]{ll2} Lam N., Lu G., \emph{A new approach to sharp Moser-Trudinger and Adams type inequalities: a rearrangement-free argument}, J. Differential Equations {\bf 255} (2013), 298-325.

\bibitem[Li] {li} Li P., \emph{Geometric Analysis}, Cambridge Studies in Advanced Mathematics, 134. Cambridge University Press, Cambridge, 2012.

\bibitem[LR]{lr} Li Y., Ruf B., \emph{A sharp Trudinger-Moser type inequality for unbounded domains in $\R^n$}, Indiana Univ. Math. J. {\bf 57} (2008), 451-480.

%\bibitem[LLZ]{llz} Li, J., Lu, G., Zhu, M., \emph{ Concentration-compactness principle for Trudinger–Moser inequalities on Heisenberg groups and existence of ground state solutions}, Calc. Var. {\bf 57}, (2018),84. 

\bibitem[LT]{lt} Lu, G., Tang, H. \emph{Sharp Moser–Trudinger Inequalities on Hyperbolic Spaces with Exact Growth Condition}, J. Geom. Anal. {\bf 26} (2016), 837–857. 

\bibitem[LTZ]{ltz} Lu G., Tang H., Zhu M., \emph{Best constants for Adams' inequalities with the exact growth condition in $\R^n$}, Adv. Nonlinear Stud. {\bf 15} (2015), 763-788.

%\bibitem[LZ]{lz}Lu G., Zhu M., \emph{A sharp Trudinger-Moser type inequality involving $L^n$ norm in the entire space $\R^n$}, Journal of Differential Equations. {\bf 267}, Issue 5, (2019),  3046-3082.

\bibitem[MS1]{ms1} Masmoudi N., Sani F., \emph{Adams' inequality with the exact growth condition in $\R^4$}, Comm. Pure Appl. Math. {\bf 67} (2014), 1307-1335.

\bibitem[MS2]{ms2} Masmoudi N., Sani F., \emph{Trudinger-Moser inequalities with the exact growth condition in $\R^n$ and applications}, Comm. Partial Differential Equations {\bf 40} (2015), 1408-1440.

\bibitem[MS3]{ms3} Masmoudi N., Sani F., \emph{Higher order Adams' inequality with the exact growth condition}, Commun. Contemp. Math. {\bf 20}, No. 06(2018),1750072.

\bibitem[NN]{nn} Ng\^{o} Q., Nguyen V., \emph{Sharp Adams-Moser-Trudinger type
inequalities in the hyperbolic space}, Rev. Mat. Iberoam. {\bf36}, (2020), 1409–1467. 

%\bibitem[Mo]{mo} Moser J., \emph{A sharp form of an inequality by N. Trudinger}, Indiana Univ. Math. J. {\bf 20}(1970/71), 1077-1092.

%\bibitem[ON]{on}O’Neil R., \emph{ Convolution operators in $L(p, q)$ spaces}, Duke Math. J. {\bf 30} (1963),  129–142.

%\bibitem[Po]{po}Pohozhaev S.I. , \emph{On the imbedding Sobolev theorem for $pl = n$}, in: Doklady Conference, Section Math., Moscow Power Inst., 1965, 158–170.

%\bibitem[RS]{rs} Ruf B., Sani F., \emph{Sharp Adams-type inequalities in $\R^n$}, Trans. Amer. Math. Soc. {\bf 365}(2013), 645-670.


\bibitem[Qin1]{qin1} Qin L., \emph{Adams inequalities with exact growth condition on $\R^n$ and the Heisenberg group}, Ph.D. dissertation, University of Missouri, Columbia (2020).

\bibitem[Qin2]{qin2} Qin L., \emph{Adams inequalities with exact growth condition
for Riesz-like potentials on $\R^n$}, Adv.  Math., {\bf 397} (2022), 108195.

\bibitem[Ruf]{ruf} Ruf B., \emph{ A sharp Trudinger–Moser type inequality for unbounded domains in $\R^2$}, J. Funct. Anal. {\bf 219} (2005) , 340–367.

%\bibitem[S1]{s1}Stein E.M., \emph{Singular integrals and differentiability properties of functions}, Princeton Mathematical Series {\bf 30}, Princeton University Press (1970).

%\bibitem[S2]{s2} Stein E.M., \emph{Harmonic analysis: real-variable methods, orthogonality, and oscillatory integrals}, Princeton University Press (1993).


\bibitem[St]{st} Strichartz R.S.,
\emph{Analysis of the Laplacian on the complete Riemannian manifold}, J. Functional Analysis {\bf52} (1983),  48–79.

\bibitem[T]{t} Triebel H., 
\emph{Theory of function spaces. II}, Monogr. Math., 84
Birkh\"auser Verlag, Basel, 1992. 


 
\end{thebibliography}
\end{document}